\documentclass[11pt]{article}
\usepackage{amsfonts}
\usepackage{amsfonts,mathrsfs,amssymb,amsthm,mathptm}
\usepackage{amsmath,amscd}
\usepackage{mathptm,pslatex}
\usepackage{color}
\usepackage[dvips]{graphicx}
\usepackage[all]{xy}
\usepackage{graphicx}
\usepackage{fancyhdr}
\usepackage{hyperref}
\usepackage{appendix}

\begin{document}
\newtheorem{Def}{Definition}[section]
\newtheorem{Bsp}[Def]{Example}
\newtheorem{Prop}[Def]{Proposition}
\newtheorem{Theo}[Def]{Theorem}
\newtheorem{Lem}[Def]{Lemma}
\newtheorem{Koro}[Def]{Corollary}
\theoremstyle{definition}
\newtheorem{Rem}[Def]{Remark}

\newcommand{\add}{{\rm add}}
\newcommand{\con}{{\rm con}}
\newcommand{\gd}{{\rm gl.dim}}
\newcommand{\dd}{{\rm dom.dim}}
\newcommand{\cdm}{{\rm codom.dim}}
\newcommand{\tdim}{{\rm dim}}
\newcommand{\E}{{\rm E}}
\newcommand{\Mor}{{\rm Morph}}
\newcommand{\End}{{\rm End}}
\newcommand{\ind}{{\rm ind}}
\newcommand{\rsd}{{\rm res.dim}}
\newcommand{\rd} {{\rm rep.dim}}
\newcommand{\ol}{\overline}
\newcommand{\overpr}{$\hfill\square$}
\newcommand{\rad}{{\rm rad}}
\newcommand{\soc}{{\rm soc}}
\renewcommand{\top}{{\rm top}}
\newcommand{\pd}{{\rm pdim}}
\newcommand{\id}{{\rm idim}}
\newcommand{\fld}{{\rm fdim}}
\newcommand{\Fac}{{\rm Fac}}
\newcommand{\Gen}{{\rm Gen}}
\newcommand{\fd} {{\rm fin.dim}}
\newcommand{\Fd} {{\rm Fin.dim}}
\newcommand{\Pf}[1]{{\mathscr P}^{<\infty}(#1)}
\newcommand{\DTr}{{\rm DTr}}
\newcommand{\cpx}[1]{#1^{\bullet}}
\newcommand{\D}[1]{{\mathscr D}(#1)}
\newcommand{\Dz}[1]{{\mathscr D}^+(#1)}
\newcommand{\Df}[1]{{\mathscr D}^-(#1)}
\newcommand{\Db}[1]{{\mathscr D}^b(#1)}
\newcommand{\C}[1]{{\mathscr C}(#1)}
\newcommand{\Cz}[1]{{\mathscr C}^+(#1)}
\newcommand{\Cf}[1]{{\mathscr C}^-(#1)}
\newcommand{\Cb}[1]{{\mathscr C}^b(#1)}
\newcommand{\Dc}[1]{{\mathscr D}^c(#1)}
\newcommand{\K}[1]{{\mathscr K}(#1)}
\newcommand{\Kz}[1]{{\mathscr K}^+(#1)}
\newcommand{\Kf}[1]{{\mathscr  K}^-(#1)}
\newcommand{\Kb}[1]{{\mathscr K}^b(#1)}
\newcommand{\DF}[1]{{\mathscr D}_F(#1)}

\newcommand{\Kac}[1]{{\mathscr K}_{\rm ac}(#1)}
\newcommand{\Keac}[1]{{\mathscr K}_{\mbox{\rm e-ac}}(#1)}

\newcommand{\modcat}{\ensuremath{\mbox{{\rm -mod}}}}
\newcommand{\Modcat}{\ensuremath{\mbox{{\rm -Mod}}}}
\newcommand{\Spec}{{\rm Spec}}

\newcommand{\stmc}[1]{#1\mbox{{\rm -{\underline{mod}}}}}
\newcommand{\Stmc}[1]{#1\mbox{{\rm -{\underline{Mod}}}}}
\newcommand{\prj}[1]{#1\mbox{{\rm -proj}}}
\newcommand{\inj}[1]{#1\mbox{{\rm -inj}}}
\newcommand{\Prj}[1]{#1\mbox{{\rm -Proj}}}
\newcommand{\Inj}[1]{#1\mbox{{\rm -Inj}}}
\newcommand{\PI}[1]{#1\mbox{{\rm -Prinj}}}
\newcommand{\GP}[1]{#1\mbox{{\rm -GProj}}}
\newcommand{\GI}[1]{#1\mbox{{\rm -GInj}}}
\newcommand{\gp}[1]{#1\mbox{{\rm -Gproj}}}
\newcommand{\gi}[1]{#1\mbox{{\rm -Ginj}}}

\newcommand{\opp}{^{\rm op}}
\newcommand{\otimesL}{\otimes^{\rm\mathbb L}}
\newcommand{\rHom}{{\rm\mathbb R}{\rm Hom}\,}
\newcommand{\pdim}{\pd}
\newcommand{\Hom}{{\rm Hom}}
\newcommand{\Coker}{{\rm Coker}}
\newcommand{ \Ker  }{{\rm Ker}}
\newcommand{ \Cone }{{\rm Con}}
\newcommand{ \Img  }{{\rm Im}}
\newcommand{\Ext}{{\rm Ext}}
\newcommand{\StHom}{{\rm \underline{Hom}}}
\newcommand{\StEnd}{{\rm \underline{End}}}

\newcommand{\KK}{I\!\!K}
\newcommand{\gm}{{\rm _{\Gamma_M}}}
\newcommand{\gmr}{{\rm _{\Gamma_M^R}}}

\def\vez{\varepsilon}\def\bz{\bigoplus}  \def\sz {\oplus}
\def\epa{\xrightarrow} \def\inja{\hookrightarrow}

\newcommand{\lra}{\longrightarrow}
\newcommand{\llra}{\longleftarrow}
\newcommand{\lraf}[1]{\stackrel{#1}{\lra}}
\newcommand{\llaf}[1]{\stackrel{#1}{\llra}}
\newcommand{\ra}{\rightarrow}
\newcommand{\dk}{{\rm dim_{_{k}}}}
\newcommand{\holim}{{\rm Holim}}
\newcommand{\hocolim}{{\rm Hocolim}}
\newcommand{\colim}{{\rm colim\, }}
\newcommand{\limt}{{\rm lim\, }}
\newcommand{\Add}{{\rm Add }}
\newcommand{\Prod}{{\rm Prod }}
\newcommand{\Tor}{{\rm Tor}}
\newcommand{\Cogen}{{\rm Cogen}}
\newcommand{\Tria}{{\rm Tria}}
\newcommand{\Loc}{{\rm Loc}}
\newcommand{\Coloc}{{\rm Coloc}}
\newcommand{\tria}{{\rm tria}}
\newcommand{\Con}{{\rm Con}}
\newcommand{\Thick}{{\rm Thick}}
\newcommand{\thick}{{\rm thick}}
\newcommand{\Sum}{{\rm Sum}}

{\Large \bf
\begin{center}
Homological theory of self-orthogonal modules
\end{center}}

\medskip\centerline{\textbf{Hongxing Chen} and \textbf{Changchang Xi}$^*$}

\renewcommand{\thefootnote}{\alph{footnote}}
\setcounter{footnote}{-1} \footnote{ $^*$ Corresponding author.
Email: xicc@cnu.edu.cn; Fax: 0086 10 68903637.}
\renewcommand{\thefootnote}{\alph{footnote}}
\setcounter{footnote}{-1} 
\renewcommand{\thefootnote}{\alph{footnote}}
\setcounter{footnote}{-1} 

\begin{abstract}
Tachikawa's second conjecture predicts that a finitely generated, self-orthogonal module over a finite-dimensional self-injective algebra is projective. This conjecture is an important part of the Nakayama conjecture. Our principal motivation of this work is a systematic understanding of finitely generated, self-orthogonal generators over a self-injective Artin algebra from the view point of stable module categories. Consequently, we give equivalent characterizations of Tachikawa's second conjecture in terms of $M$-Gorenstein categories, and  establish a recollement of the $M$-relative stable categories for a self-orthogonal generator $M$. Further, we show that the Nakayama conjecture holds true for Gorenstein-Morita algebras.
\end{abstract}

{\footnotesize\tableofcontents\label{contents}}

\section{Introduction}\label{Introduction}
Since about half a century homological conjectures form a core set of problems in
representation theory and homological algebra of finite-dimensional algebras. One of the most
prominent conjectures in this system of closely related conjectures is the Nakayama conjecture posed by Nakayama in \cite{Nakayama}.

\smallskip
\textbf{(NC)} If a finite-dimensional algebra over a field
has infinite dominant dimension, then it is self-injective.

\smallskip
This conjecture has been verified
for a few classes of algebras, for which the conjecture could be checked more or less directly by
clever computations. Despite many efforts, taking various approaches, very little is known about the
homological conjectures and in particular about Nakayama's conjecture.

To deal with the conjecture, Tachikawa proposed another two homological conjectures (see \cite{Tac}), they are now called Tachikawa's first and second conjecture.

\smallskip
{\bf (TC1)} If a finite-dimensional algebra $A$ over a field $k$ satisfies $\Ext^n_A(D(A),A)=0$ for all $n\geq 1$, then $A$ is self-injective, where $D=\Hom_k(-,k)$ is the usual duality.

{\bf (TC2)}  Let $A$ be a finite-dimensional self-injective algebra and $M$ a finitely generated $A$-module. If $M$ is self-orthogonal, that is,  $\Ext_A^n(M, M)=0$ for all $n\geq 1$, then $M$ is projective.

\medskip
Concerning (TC1) and (TC2), there are only a few cases verified. For instance, (TC1) holds for special local algebras and commutative algebras (see \cite{Asashiba} and \cite{abs, GoTa}, respectively), while (TC2) holds true for self-injective algebras of finite representation type, symmetric algebras with radical cube zero, and local self-injective algebras with radical cube zero (see \cite{Tac, Hoshino, Sch}, respectively ). For further information on these conjectures, we refer to \cite{Yamagata}.

The validity of both (TC1) and (TC2) is equivalent to the one of (NC). Moreover, by a result of Mueller \cite{Muller}, a pair $(A, M)$, with $A$ a self-injective algebra and $M$ a finitely generated $A$-module, satisfies (TC2) if and only if the endomorphism algebra $\End_A(A\oplus M)$ satisfies (NC).

Thus it is of significant interest to understand self-orthogonal modules of the form $A\oplus M$ for a self-injective algebra $A$ and finitely generated $A$-modules $M$. Generally, a self-orthogonal module of the form $B\oplus Y$ over an Artin algebra $B$ with $Y$ a finitely generated $B$-module is termed \emph{self-orthogonal generators}.

In this paper, we investigate self-orthogonal generators over self-injective algebras from the point of view of stable categories. More precisely, we first establish a general theory for arbitrary (not necessarily self-orthogonal) generator by  constructing two pairs of triangle endofunctors for stable module categories, and then establish specially a recollement of the relative stable categories for a self-orthogonal generator. Finally, we describe compact objects of the right term of the recollement by the heart of a torsion pair in the stable module category. Based on these investigations, we give equivalent characterizations of (TC2) and show that the Nakayama conjecture holds for Gorenstein-Morita algebras.

\subsection{Equivalent characterizations of Tachikawa's second conjecture}

In this section, we present our equivalent characterizations of (TC2) in terms of  perpendicular categories or special modules associated with self-orthogonal generators. We then introduce the notion of Gorenstein-Morita algebras and state one of our main results, namely (NC) holds for Gorenstein-Morita algebras.

We begin with recalling a few notation and terminology.

Let $A$ be an Artin algebra. We denote by $A\Modcat$ (respectively, $A\modcat$) the category of (respectively, finitely generated) left $A$-modules, by $D$ the usual duality on $A\modcat$ and by $\nu_A$ the Nakayama functor ${_A}D(A)\otimes_A-$. For $M\in A\Modcat$, let $\Add(M)$ (respectively, $\add(M)$) be the full subcategory of $A\Modcat$ consisting of direct summands of (respectively, finite) direct sums of
copies of $M$, and let $M^{\bot 1}$ be the full subcategory of $A\Modcat$ consisting of modules $X$ with $\Ext^1_A(M,X)=0$. We say that $M$ is \emph{self-orthogonal} if $\Ext^i_A(M,M)=0$ for all $i\ge 1$; \emph{Nakayama-stable} if $\add(M)=\add(\nu_A(M))$; and a \emph{generator} if
$A\in\add(M)$. Clearly, $_AA$ is Nakayama-stable if $A$ is a self-injective algebra. Every module over a symmetric algebra $A$ (that is, $A\simeq D(A)$ as $A$-$A$-bimodules) is Nakayama-stable.

Let $A$ be a self-injective algebra, $M\in A\modcat$ a self-orthogonal generator and $\Lambda:=\End_A(M)$. Our strategy is to understand the relation between the category of Gorenstein-projective $\Lambda$-modules and the category of two-sided perpendicular category $\mathscr{G}$ of $M$, where
$$
\mathscr{G}:=\{X\in A\Modcat\mid \Ext_A^n(M, X)=0=\Ext_A^n(X,M)\;\;\mbox{for all}\;\; n\geq 1\}$$
has $A$-modules left and right orthogonal to ${_A}M$. These two categories are all Frobenius and equivalent (see Lemma \ref{Gorenstein}). This builds a new bridge between (TC2) and (NC).
Following \cite{SSW}, $\mathscr{G}$ is called an \emph{$M$-Gorenstein subcategory} in $A\Modcat$.
Then the quotient category of $\mathscr{G}$ modulo $\Add({_A}M)$, denoted by
$$ \mathscr{C}:=\mathscr{G}/[M],$$
is a triangulated category and equivalent to the stable category of Gorenstein-projective $\Lambda$-modules.
The category $\mathscr{C}$ is called an \emph{$M$-Gorenstein stable category}.
In particular, if $M=A$, then $\mathscr{C}$ is the usual \emph{stable module category} of $A$, denoted by
$\Stmc{A}$.

Next, we introduce two classes of $A$-modules determined by $M$.

\begin{Def} \label{Countably filtered} Let $X$ be an $A$-module.

$(i)$ $X$ is \emph{$M$-compact} if it is a compact object in the category $\mathscr{C}$, that is,
$X\in\mathscr{G}$ and the functor $\Hom_\mathscr{C}(X,-):\mathscr{C}\to\mathbb{Z}\Modcat$ commutes with coproducts.

$(ii)$ $X$ is \emph{$M$-filtered} if it has a filtration
$0=X_0\subseteq X_1\subseteq X_2\subseteq \cdots\subseteq  X_n\subseteq\cdots \subseteq X$
in $A\Modcat$ such that $X=\bigcup_{n=0}^{\infty}X_n$ and the subquotient $X_{n+1}/X_n$ of $X$ is isomorphic to a finite direct sum of $A$-modules in the set $\{_AA\}\cup\{\Omega_A^{-i}(M)\mid i\in\mathbb{N}\}$ for $n\in\mathbb{N}$. If $X=X_n$ for an integer $n$ in the filtration, then $X$ is said to be \emph{finitely $M$-filtered}.
\end{Def}

Finitely generated modules in $\mathscr{G}$ are $M$-compact, and finitely $M$-filtered $A$-modules are exactly finitely generated, $M$-filtered $A$-modules. Further, $M$-compact, finitely $M$-filtered $A$-modules lie in $\add(M)$ by Lemma \ref{FMF}. Clearly, $A$-compact modules are exactly $A$-modules that are isomorphic in $\Stmc{A}$ to finitely generated modules.

Now, our characterizations of (TC2) for Nakayama-stable generators read as follows.

\begin{Theo} \label{Conjecture}
Let $A$ be a self-injective Artin algebra and $M$ a self-orthogonal and Nakayama-stable generator for $A\modcat$.  The following are equivalent.

$(1)$ $M$ is a projective $A$-module.

$(2)$ $\mathscr{G}$ coincides with the full subcategory of $A\Modcat$ consisting of all filtered colimits of finitely generated modules in $\mathscr{G}$.

$(3)$ Any $M$-compact and $M$-filtered $A$-module lies in $\Add(M)$.

$(4)$ The minimal left $\mathscr{G}$-approximation $W$ of $\Omega_A^{-}(M)$ is a filtered colimit of finitely generated modules in $\mathscr{G}$.

$(5)$ The minimal left $\mathscr{G}$-approximation $W$ of $\Omega_A^{-}(M)$ has the property: the category $W^{\bot 1}$ is closed under countable direct sums in $A\Modcat$ of finitely $M$-filtered $A$-modules.
\end{Theo}

In Theorem \ref{Conjecture}, $(2)$ is equivalent to saying that the algebra $\Lambda$ is virtually Gorenstein (see Proposition \ref{Limits}) in the sense of Beligiannis; $\Add(M)$ contains finitely generated, $M$-compact and $M$-filtered $A$-modules; $(4)$ and $(5)$ hold true if the module $W$ is the direct sum of finitely generated $A$-modules. Moreover, Theorem \ref{Conjecture} implies that (TC2) holds for symmetric algebras of finite representation type because $(2)$-$(5)$ in Theorem \ref{Conjecture} are satisfied. This can be seen from a classical result, due to Auslander and Ringel-Tachikawa, that any module over an Artin algebra of finite representation type is a direct sum of finitely generated modules. Thus Theorem \ref{Conjecture} provides a different approach to studying (TC2).

\subsection{Nakayama conjecture for Gorenstein-Morita algebras}
As indicated by the relation between (TC2) and (NC), we can apply Theorem \ref{Conjecture} to discuss (NC) for \emph{strongly Morita} algebras, which are, by definition, the endomorphism algebras of Nakayama-stable generators over self-injective algebras. For this purpose, we focus on two classes of modules that are associated with compact objects in some stable categories. This leads to introducing the notions of compactly Gorenstein algebras and Gorenstein-Morita algebras in terms of these modules.

\begin{Def} \label{Gorenstein-Morita}
Let $B$ ba an Artin algebra and $Y$ a $B$-module.

$(i)$ The $B$-module $Y$ is \emph{compactly filtered} if it has a filtration
$0=Y_0\subseteq Y_1\subseteq Y_2\subseteq \cdots \subseteq Y_n\subseteq\cdots \subseteq Y$ of $B$-modules
such that $Y=\bigcup_{n=0}^{\infty}Y_n$ and the subquotient $Y_{n+1}/Y_n$ of $Y$ is isomorphic to a finitely generated $B$-module of finite projective dimension for $n\in\mathbb{N}$; and \emph{compactly Gorenstein-projective} if it is a compact object in the stable category of Gorenstein-projective $B$-modules.

$(ii)$ The algebra $B$ is \emph{compactly Gorenstein} if any compactly filtered, compactly Gorenstein-projective $B$-module is projective; and \emph{Gorenstein-Morita} if $B$ is both strongly Morita and compactly Gorenstein.
\end{Def}

Clearly, finitely generated and compactly filtered modules are exactly finitely generated modules of finite projective dimension, while finitely generated Gorenstein-projective modules are compactly Gorenstein-projective. Compactly Gorenstein algebras include virtually Gorenstein algebras (see \cite{Bel, Bel2}) and algebras of finite finitistic dimension (see Lemma \ref{Three conditions}). Moreover, compactly Gorenstein algebras over a field are closed under derived equivalences and stable equivalences of adjoint type by Corollary \ref{DM}. Since the Nakayama functor of a symmetric algebra is the identity functor, strongly Morita algebras capture \emph{gendo-symmetric algebras} which are, by definition, the endomorphism algebras of generators over symmetric algebras (see \cite{FK,xc6}). Examples of gendo-symmetric algebras include Hecke algebras, (quantized) Schur algebras, and blocks of the Bernstein-Gelfand-Gelfand category $\mathcal{O}$ of semisimple complex Lie algebras.

As a consequence of Theorem \ref{Conjecture}, we show that (NC) holds for Gorenstein-Morita algebras.

\begin{Koro}\label{Gendo-symmetric}
Let $B$ be a Gorenstein-Morita algebra. If $B$ has infinite dominant dimension, then it is self-injective.
In particular, any gendo-symmetric, virtually Gorenstein algebra with infinite dominant dimension
is symmetric.
\end{Koro}

As is known, \emph{not} all Artin algebras are virtually Gorenstein (see \cite{BK}), however, we would like to conjecture that \emph{all Artin algebras are compactly Gorenstein}. If this is true, then (NC) holds for strongly Morita algebras, and in particular for gendo-symmetric algebras.

In the preprint \cite{xcf}, we discuss (TC2) for symmetric algebras in terms of recollements of derived module categories and stratifying ideals of algebras. Consequently, it is shown that the validity of (TC2) for symmetric algebras is equivalent to saying that no indecomposable symmetric algebras  have stratifying ideals apart from themselves and $0$.

\subsection{Recollements of relative stable categories from self-orthogonal modules}
To prove Theorem \ref{Conjecture}, we first investigate arbitrary (not necessarily self-orthogonal) generator over a self-injective algebra, and construct two pairs of triangle endofunctors for stable module categories. This enables us to establish a recollement of the relative stable categories for a self-orthogonal generator. By employing the heart of a torsion pair in the stable module category, we then characterize the compact objects of the right term of this recollement.

Let $A$ be a self-injective algebra and $M$ a generator for $A\modcat$. For a full subcategory $\mathscr{X}$ of $A\Modcat$, we denote by $\mathscr{X}/[M]$ the quotient category of $\mathscr{X}$ modulo $\Add(M)$. In particular, $A\Modcat/[A]$ is the same as the stable module category $\Stmc{A}$. For simplicity, we denote by $\underline{\Hom}(X,Y)$ the Hom-set in $\Stmc{A}$ for $A$-modules $X$ and $Y$, and define $\Gamma:=\underline{\End}_A(M)$, called the \emph {stable endomorphism algebra} of ${_A}M$. Let $$\underline{M}^{\bot}:=\{X\in\Stmc{A}\mid \StHom_A(M, X[n])=0\;\;\mbox{for all}\;\; n\in\mathbb{Z}\}.$$

Given the pair $(A,M)$, we construct explicitly two pairs of triangle endofunctors of $\Stmc{A}$ (see Section \ref{ECTF} for details):
$$(\Phi, \Psi)\;\;\mbox{and}\;\; (\Phi', \Psi'):\;\; \Stmc{A}\lra \Stmc{A},$$
and define $\mathscr{S}:=\{X\in A\Modcat\mid \Psi(\underline{X})=0\}.$
If ${_A}M$ is additionally self-orthogonal or $\Omega$-periodic (that is, $\Omega_A^n(X)\simeq X$ in $\Stmc{A}$ for a positive integer $n$), then $\mathscr{S}$ is the smallest thick subcategory of $A\Modcat$ containing $M$ and being closed under direct sums (see Corollary \ref{localizing}), and  the above endofunctors contribute to building the recollement of $\Stmc{A}$ in Theorem \ref{Main result}.

Now, suppose that $M$ is a self-orthogonal and Nakayama-stable generator for $A\modcat$. Further, we consider the following two categories associated with $M$:
$$\mathscr{H}:=\{X\in\underline{\mathscr{S}}\mid \StHom_A(M, X[n])=0\;\;\mbox{for}\;\; n\neq 0\} \mbox{ and }$$
$$\mathscr{E}:=\{X\in\mathscr{G}\mid \StHom_A(M, X),\, \StHom_A(M[1], X)\in \Gamma\modcat\},$$
where $\mathscr{H}$ is the heart of a torsion pair in $\Stmc{A}$ determined by $M$, and thus an abelian category (see the beginning of Section \ref{RRR}).

The main result on constructing recollements of relative stable categories reads as follows.

\begin{Theo}\label{Gorenstein-recollement}
Let $A$ be a self-injective Artin algebra and $M$ a self-orthogonal and Nakayama-stable generator for $A\modcat$. Then the following hold.

$(1)$ There exists a recollement of triangulated categories:
$$
\xymatrix@C=1.2cm{\underline{M}^{\bot}\ar[r]
&\mathscr{C}\ar@/^1.2pc/[l]\ar[r]\ar@/^1.2pc/[l]\ar@/_1.2pc/[l]
&(\mathscr{G}\cap\mathscr{S})/[M]\ar@/^1.2pc/[l]\ar@/^1.2pc/[l]\ar@/_1.2pc/[l]}
$$

$(2)$ The recollement in $(1)$ restricts to a recollement of triangulated categories:

$$
\xymatrix@C=1.2cm{\underline{M}^{\bot}\ar[r]
&\mathscr{E}/[M]\ar@/^1.2pc/[l]\ar[r]\ar@/^1.2pc/[l]\ar@/_1.2pc/[l]
&(\mathscr{E}\cap\mathscr{S})/[M].\ar@/^1.2pc/[l]\ar@/^1.2pc/[l]\ar@/_1.2pc/[l]}
\bigskip
$$

$(3)$ $\dim\big((\mathscr{E}\cap\mathscr{S})/[M]\big)\leq min\{2\; LL(\Gamma)-1, 2\;\gd(\Gamma)+1\}$, where $LL(\Gamma)$ and $\gd(\Gamma)$ denote the Loewy length and global dimension of the algebra $\Gamma$, respectively.
\end{Theo}

Of importance are the compact objects for triangulated categories. Since $\mathscr{C}$ is compactly generated, it follows from  basic properties of recollements  that $\underline{M}^{\bot}$ has compact objects. It seems, however, to be unclear that the category $(\mathscr{G}\cap\mathscr{S})/[M]$ has compact objects. In the following, we gives a complete description of its compact objects (see also Corollary \ref{Compact recollement} for details).

Let $\mathcal{T}$ be a triangulated category and  $\mathcal{U}$ a set  of objects in $\mathcal{T}$. For integers $i\leq j$ and $n \ge 0$, we denote by $\langle\mathcal{U}\rangle_{n+1}^{[i,j]}$ the full subcategory of $\mathcal{T}$ consisting of all objects obtained by taking $(n+1)$-fold extensions of finite direct sums of objects in the set $\{U[-s]\mid U\in\mathcal{U},\, s\in\mathbb{Z},\, i\leq s\leq j\}$.

\begin{Prop}\label{Compact objects in recollement} Let $A$ be a self-injective Artin algebra and $M$ a self-orthogonal, Nakayama-stable generator. Then the following hold.

$(1)$ Each object $X$ of $\mathscr{E}\cap\mathscr{S}$ is $M$-compact and isomorphic in $\Stmc{A}$ to an $M$-filtered module.
Moreover, $X$ is finitely generated if and only if $X\in\add({_A}M)$.

$(2)$ The category $(\mathscr{G}\cap\mathscr{S})/[M]$ is a compactly generated triangulated category and has $(\mathscr{E}\cap\mathscr{S})/[M]$ as its full subcategory consisting of all compact objects.

$(3)$ Let $\mathcal{S}$ be the set of isomorphism classes of simple objects of $\mathscr{H}$, and let $n$ be the Loewy length of the algebra $\Gamma$. Then
$\big(\mathscr{G}\cap\mathscr{S}\big)/[M]=\langle\Add(\mathcal{S})\rangle_{2n}^{[-1,0]}$ and
$(\mathscr{E}\cap\mathscr{S})/[M]=\langle\mathcal{S}\rangle_{2n}^{[-1,0]}$, where $\mathcal{S}$ is a finite set.
\end{Prop}

Under the assumption of Theorem \ref{Gorenstein-recollement}, the module
${_A}M$ is projective if and only if $(\mathscr{G}\cap\mathscr{S})/[M]=0$ if and only if $(\mathscr{E}\cap\mathscr{S})/[M]=0$ (see Corollary \ref{Vanishing}). Thus  (TC2) is true for the pair $(A, M)$ exactly when the recollements in Theorem \ref{Gorenstein-recollement} are trivial. Hence, to construct a counterexample to (TC2),  our results, Theorem \ref{Gorenstein-recollement} and Proposition \ref{Compact objects in recollement}, provide necessary homological information on self-orthogonal modules.

\medskip
\subsection{Overview of the contexts}
The contents of this article are sketched as follows. In Section \ref{sect2} we briefly recall definitions of quotient categories, recollements and Gorenstein-projective modules over algebras. In Section \ref{RSMC} we construct two pairs of triangle endofunctors of the stable module category $\Stmc{A}$ for a self-injective algebra $A$ with a generator ${_A}M$. With these endofunctors, we establish the recollement in Theorem \ref{Main result} of $\Stmc{A}$ determined by $M$. Moreover, we show that the subcategory $\mathscr{G}$ of $A\Modcat$ relative to $M$ is equivalent to the category of Gorenstein-projective modules over the endomorphism algebra of $M$ (see Lemma \ref{Gorenstein}). In Section \ref{RRR} we prove Theorem \ref{Gorenstein-recollement} and establish a representability theorem for a series of homological functors (see Theorem \ref{Compact}). In Section \ref{Characterization} we show Theorem \ref{Conjecture} and Corollary \ref{Gendo-symmetric}.

\section{Preliminaries\label{sect2}}
In this section we briefly recall definitions, basic facts and notation used in this paper.

\subsection{Quotient categories and recollements}

Let $\mathcal C$ be an additive category.

A full subcategory $\mathcal B$ of $\mathcal
C$ is always assumed to be closed under isomorphisms, that is, if
$X\in {\mathcal B}$ and $Y\in\cal C$ with $Y\simeq X$, then
$Y\in{\mathcal B}$.

Let $X$ be an object in $\mathcal{C}$. The full subcategory of $\mathcal{C}$ consisting of all direct summands of
finite coproducts of copies of $X$ is denoted by $\add(X)$. If $\mathcal{C}$ admits
coproducts (that is, coproducts indexed over sets exist in
${\mathcal C}$), then $\Add(X)$ denotes the full subcategory of
$\mathcal{C}$ consisting of all direct summands of  coproducts
of copies of $X$. Dually, if $\mathcal{C}$ admits products,
then $\Prod(X)$ denotes the full subcategory of $\mathcal{C}$
consisting of all direct summands of products of copies of
$X$.

For morphisms $f: X\to Y$ and $g: Y\to Z$ in $\mathcal C$, the composition of $f$ and $g$ is written by $fg$, a morphism from
$X$ to $Z$. The induced morphisms $\Hom_{\mathcal
C}(Z,f):\Hom_{\mathcal C}(Z,X)\ra \Hom_{\mathcal C}(Z,Y)$ and
$\Hom_{\mathcal C}(f,Z): \Hom_{\mathcal C}(Y, Z)\ra \Hom_{\mathcal
C}(X, Z)$ are denoted by $f^*$ and $f_*$, respectively.

For functors $F:\mathcal {C}\to
\mathcal{D}$ and $G: \mathcal{D}\to \mathcal{E}$, the composition of $F$ and $G$ is denoted by $G\circ F$ which is a functor from $\mathcal C$ to
$\mathcal E$. Let $\Ker(F)$ and $\Img(F)$ be the kernel and image of
the functor $F$, respectively. In particular, $\Ker(F)$ is closed
under isomorphisms in $\mathcal{C}$. In this paper, we require that
$\Img(F)$ is closed under isomorphisms in $\mathcal{D}$.

Suppose that $\mathcal{B}$ is a full subcategory of $\mathcal{C}$. A morphism $f: X\to Y$ in $\mathcal{C}$ is called a
\emph{right $\mathcal{B}$-approximation} of $Y$ if $X\in \mathcal{B}$ and
$\Hom_\mathcal{C}(B, f): \Hom_\mathcal{C}(B, X)\to \Hom_\mathcal{C}(B, Y)$ is surjective for any $B\in\mathcal{B}$; and
\emph{right minimal} if $\alpha\in\End_{\mathcal{C}}(X)$ is an isomorphism whenever $f=\alpha f$.
If $f$ is both a right minimal morphism and a right $\mathcal{B}$-approximation of $Y$, then $f$ is called a \emph{minimal right $\mathcal{B}$-approximation} of $Y$. In this case, the object $X$ is unique up to isomorphism and is called the minimal right $\mathcal{B}$-approximation of $Y$ (without mentioning $f$).
If each object of $\mathcal{C}$ admits a right $\mathcal{B}$-approximation, then $\mathcal{B}$ is said to be \emph{contravariantly finite} in $\mathcal{C}$.
Dually, there are the notions of (minimal) left approximations and \emph{covariantly finite} subcategories in $\mathcal{C}$.
If $\mathcal{B}$ is both contravariantly and covariantly finite in $\mathcal{C}$, then it is said to be \emph{functorially finite} in  $\mathcal{C}$.

We recall Wakamatsu's Lemma (see \cite[Proposition 1.3]{ar1990}):
Let $\mathcal{S}$ be a class of $R$-modules over a ring $R$ closed under extensions. If $f: C\ra X$ is a minimal right $\mathcal{S}$-approximation of an $R$-module $X$, then $\Ext^1_R(L,\Ker(f))=0$ for $L\in \mathcal{S}$. Dually, if $g: X\ra C'$ is a minimal left $\mathcal{S}$-approximation of an $R$-module $X$, then $\Ext^1_R(\Coker(g),M)=0$ for $M\in \mathcal{S}$.

Next, we recall the definition of quotient categories of additive categories.

Let $\mathcal{D}$ be a full subcategory of $\mathcal{C}$. Denote by $\mathcal{C}/\mathcal{D}$ the \emph {quotient category} of $\mathcal{C}$ modulo  $\mathcal{D}$. It has the same objects as $\mathcal{C}$, but its morphism set for any two objects
$X$ and $Y$ is given by $\Hom_{\mathcal{C}/\mathcal{D}}(X, Y):=\Hom_\mathcal{C}(X, Y)/\mathcal{D}(X, Y)$
where $\mathcal{D}(X, Y)$ is the subgroup of $\Hom_{\mathcal{C}}(X,Y)$ consisting of all morphisms
factorizing through objects in $\mathcal{D}$. The canonical quotient functor  $q: \mathcal{C}\to \mathcal{C}/\mathcal{D}$
sends a morphism $f:X\to Y$ in $\mathcal{C}$ to $f+\mathcal{D}(X, Y)$ in $\mathcal{C}/\mathcal{D}$.
Clearly, $\Ker(q)$ consists of all direct summands (in $\mathcal{C}$) of objects of $\mathcal{D}$.

Suppose that $\mathcal{C}$ admits coproducts. An object $X$ is said to be \emph {compact} in $\mathcal{C}$ if
the functor $\Hom_\mathcal{C}(X,-)$ from $\mathcal{C}$ to the category of abelian groups commutes with coproducts.
The full subcategory of $\mathcal{C}$ consisting of compact objects is denoted by $\mathcal{C}^{\rm c}$.
A set $\mathcal{U}$ of objects of $\mathcal{C}$ is called a \emph{compact generating set} of $\mathcal{C}$ if each object of $\mathcal{U}$ is compact in $\mathcal{C}$ and an object $X\in\mathcal{C}$ is zero whenever $\Hom_\mathcal{C}(U, X)=0$ for all $U\in\mathcal{U}$.
When $\mathcal{C}$ is a triangulated category, it is said to be \emph{compactly generated} if it has a compact generating set.
If $\mathcal{U}$ is a set of compact objects of $\mathcal{C}$ closed under shifts, then $\mathcal{U}$ is a compact generating set of $\mathcal{C}$ if and only if $\mathcal{C}$ itself is the smallest full triangulated subcategory of $\mathcal{C}$ containing $\mathcal{U}$ and being closed under coproducts.

The following result is elementary.

\begin{Lem}\label{AJCP}
$(1)$ Let $\mathcal{C}$ and $\mathcal{D}$ be additive categories, and let $\mathcal{X}\subseteq \mathcal{C}$ and $\mathcal{Y}\subseteq\mathcal{D}$ be full subcategories. Suppose that $F:\mathcal{C}\to\mathcal{D}$ and $G:\mathcal{D}\to \mathcal{C}$ form an adjoint pair $(F, G)$ of additive functors.
If $F(\mathcal{X})\subseteq\mathcal{Y}$ and $G(\mathcal{Y})\subseteq\mathcal{X}$, then the adjoint pair $(F,G)$ induces an adjoint pair $(F_0,G_0)$ of  additive functors $F_0:\mathcal{C}/\mathcal{X}\to\mathcal{D}/\mathcal{Y}$ and $G_0:\mathcal{D}/\mathcal{Y}\to\mathcal{C}/\mathcal{X}$.

$(2)$ Suppose that an additive category $\mathcal{C}$ admits coproducts and $\mathcal{X}$ is a full subcategory of $\mathcal{C}$ closed under coproducts. Then $\mathcal{C}/\mathcal{X}$ admits coproducts, and the quotient functor $\mathcal{C}\to \mathcal{C}/\mathcal{X}$ preserves coproducts and compact objects.
\end{Lem}

We denote by $\C{\mathcal{C}}$ the category of all complexes over $\mathcal{C}$ with chain maps, and $\K{\mathcal{C}}$ the homotopy category of $\C{\mathcal{C}}$. For a chain map $\cpx{f}: \cpx{X}\to\cpx{Y}$ in $\C{\mathcal{C}}$, we denote by $\Cone(\cpx{f})$ the mapping cone of $\cpx{f}$. There is a distinguished triangle $\cpx{X}\to\cpx{Y}\to\Cone(\cpx{f})\to \cpx{X}[1]$ in $\K{\mathcal{C}}$. When $\mathcal{C}$ is a full subcategory of an abelian category $\mathcal{A}$, we denote by $\Kac{\mathcal{C}}$ the full subcategory of $\K{\mathcal{A}}$ consisting of acyclic complexes of $\mathcal{C}$. Clearly, if $\cpx{f}$ is a quasi-isomorphism, then $\Cone(\cpx{f})$ is acyclic. If $\mathcal{C}$ is an abelian category, we denote by $\D{\mathcal{C}}$ the \emph{unbounded derived category} of $\mathcal{C}$, which is the localization of $\K{\mathcal{C}}$ by inverting all quasi-isomorphisms. Clearly, $\K{\mathcal{A}}$ and $\D{\mathcal{C}}$ are triangulated categories.

Next, we recall the notion of recollements of triangulated
categories, introduced in \cite{BBD} for studying derived categories of
perverse sheaves over singular spaces.

\begin{Def}\label{def01} \rm
Let  $\mathcal{T}$, $\mathcal{T'}$ and $\mathcal{T''}$ be
triangulated categories. $\mathcal{T}$ is called a
\emph{recollement} of $\mathcal{T'}$ and $\mathcal{T''}$ (or there is a recollement among $\mathcal{T'}, \mathcal{T}$ and $\mathcal{T}''$) if there
are six triangle functors
$$\xymatrix@C=1.5cm{\mathcal{T''}\ar^-{i_*=i_!}[r]&\mathcal{T}\ar^-{j^!=j^*}[r]
\ar^-{i^!}@/^1.2pc/[l]\ar_-{i^*}@/_1.6pc/[l]
&\mathcal{T'}\ar^-{j_*}@/^1.2pc/[l]\ar_-{j_!}@/_1.6pc/[l]}$$ among the three categories such
that

$(1)$ $(i^*,i_*),(i_!,i^!),(j_!,j^!)$ and $(j^*,j_*)$ are adjoint
pairs,

$(2)$ $i_*,j_*$ and $j_!$ are fully faithful functors,

$(3)$ $j^! i_!=0$ (and thus also $i^!j_*=0$ and $i^*j_!=0$), and

$(4)$ for an object $X\in\mathcal{T}$, there are two triangles
$i_!i^!(X)\to  X\to j_*j^*(X)\to i_!i^!(X)[1]$ and $ j_!j^!(X)\to
X\to i_*i^*(X)\to j_!j^!(X)[1]$ in $\mathcal{T}$ induced by the counits and units of the adjunctions,
where [$1$] denotes the shift functor of $\mathcal{T}$.
\end{Def}

By a \emph{half recollement} among $\mathcal{T'}$, $\mathcal{T}$
and $\mathcal{T''}$, we mean that $i^*,i_*, j_!$ and $j^!$ satisfy the corresponding properties $(1)$-$(4)$ involved in Definition \ref{def01}. Note that there is a one-to-one correspondence between equivalence classes of half recollements (respectively, recollements) of triangulated categories and hereditary torsion pairs (respectively, TTF triples) of triangulated categories.
Recall that a \emph{torsion pair} in $\mathcal{T}$ is a pair  $(\mathcal{X}, \mathcal{Y})$ of full subcategories $\mathcal{X},\mathcal{Y}$ of $\mathcal{T}$ satisfying the three conditions:

$(a)$ $\Hom_{\mathcal{T}}(X, Y)=0$ for $X\in\mathcal{X}$ and $Y\in\mathcal{Y}$;

$(b)$ $\mathcal{X}[1]\subseteq\mathcal{X}$ and $\mathcal{Y}[-1]\subseteq\mathcal{Y}$; and

$(c)$ for any $M\in \mathcal{T}$, there is a triangle
$X\to M\to Y\to X[1]$ in $\mathcal{T}$ with $X\in\mathcal{X}$ and $Y\in\mathcal{Y}$.

\medskip
A torsion pair $(\mathcal{X}, \mathcal{Y})$ of $\mathcal{T}$ is said to be \emph{hereditary} if $\mathcal{X}$ (or equivalently, $\mathcal{Y})$ is a full triangulated subcategory of $\mathcal{T}$.
In this case, the inclusion $\mathcal{X}\to\mathcal{T}$ has a right adjoint, the inclusion $\mathcal{Y}\to\mathcal{T}$ has a left adjoint, and there is a half recollement among $\mathcal{X}, \mathcal{T}$ and $\mathcal{Y}$. If $(\mathcal{X}, \mathcal{Y})$ and $(\mathcal{Y}, \mathcal{Z})$ are hereditary torsion pairs in $\mathcal{T}$, then $(\mathcal{X}, \mathcal{Y}, \mathcal{Z})$ is called a TTF \emph{(torsion-torsionfree) triple} in $\mathcal{T}$. In this case, there is a recollement among $\mathcal{X}$, $\mathcal{T}$ and $\mathcal{Y}$. Conversely, the recollement in Definition \ref{def01} gives a TTF triple $(\Img(j_!), \Img(i_*),  \Img(j_*))$ in $\mathcal{T}$.
For more details, see \cite[Chap. 9]{Neemanbook}, \cite[Chap. I. 2]{BI} or \cite[Section 2.3]{xc1}.

Hereditary torsion pairs can be constructed in compactly generated triangulated categories as follows.

Let $\mathcal{T}$ be a compactly generated triangulated category. Then coproducts and products indexed by sets exist in $\mathcal{T}$.
Let $\mathcal{S}$ be a set of objects in $\mathcal{T}$. Denote by $\Loc_\mathcal{T}(\mathcal{S})$, $\Coloc_\mathcal{T}(\mathcal{S})$ and $\thick_\mathcal{T}(\mathcal{S})$ the smallest full triangulated subcategories of $\mathcal{T}$ containing $\mathcal{S}$ and being closed under coproducts, products and direct summands, respectively. If $\mathcal{T}$ is clearly understood in the context, we shall write $\Loc(\mathcal{S})$, $\Coloc(\mathcal{S})$ and $\thick(\mathcal{S})$ for $\Loc_\mathcal{T}(\mathcal{S})$, $\Coloc_\mathcal{T}(\mathcal{S})$ and $\thick_\mathcal{T}(\mathcal{S})$, respectively.
By $\mathcal{S}^{\bot}$ we denote the right
orthogonal full subcategory of $\mathcal{T}$ with respect to $\mathcal{S}$, that is, $\mathcal{S}^{\bot}=\{C\in \mathcal{T}\mid \Hom_{\mathcal{T}}(S[n],C)=0, \forall\; S\in \mathcal{S}, n\in\mathbb{Z} \}.$ Then $\mathcal{S}^{\bot}$ is a triangulated subcategory of
$\mathcal{T}$ closed under products. Similarly, ${^{\bot}}\mathcal{S}$
stands for the left orthogonal full subcategory of $\mathcal{T}$ with respect
to $\mathcal{S}$.

In general, the opposite category of a compactly generated triangulated category is not compactly generated, but is perfectly generated in the sense of Krause (see \cite{Kr0}).
It is worth mentioning that perfectly generated triangulated categories not only generalize compactly generated triangulated categories, but also satisfy the Brown representability theorem in \cite[Theorem A]{Kr0}. This implies the following result:

\emph{If a triangle functor from a perfectly generated triangulated category to another triangulated category preserves coproducts, then it has a right adjoint.}

The next result seems to be known. For the convenience of the reader, we include here a proof.

\begin{Prop}\label{CPG}
If $\mathcal{T}$ is a compactly generated triangulated category and  $\mathcal{S}$ is a set of objects in $\mathcal{T}$, then $(\Loc(\mathcal{S}), \mathcal{S}^{\bot})$ and $({^{\bot}}\mathcal{S}, \Coloc(\mathcal{S}))$ are hereditary torsion pairs in $\mathcal{T}$.
\end{Prop}
{\it Proof.} The Verdier localization $Q_\mathcal{S}: \mathcal{T}\to \mathcal{T}/\Loc(\mathcal{S})$ preserves coproducts
and $\Ker(Q_\mathcal{S})=\Loc(\mathcal{S})$. Since $\mathcal{T}$ is compactly generated, it is perfectly generated. It follows that $Q_\mathcal{S}$ has a right adjoint. Then $(\Loc(\mathcal{S}), \Loc(\mathcal{S})^{\bot})$ is a torsion pair in $\mathcal{T}$ by \cite[Theorem 9.1.13]{Neemanbook}. Since $\Loc(\mathcal{S})^{\bot}=\mathcal{S}^{\bot}$,
$(\Loc(\mathcal{S}), \mathcal{S}^{\bot})$ is a torsion pair in $\mathcal{T}$. Moreover, it is hereditary because $\Loc(\mathcal{S})$ is a full triangulated subcategory of $\mathcal{T}$.
This proof also implies that $(\Loc(\mathcal{S}), \mathcal{S}^{\bot})$ is a hereditary torsion pair in $\mathcal{T}$ whenever $\mathcal{T}$ is perfectly generated.

Clearly, $\Coloc(\mathcal{S})$ can be regarded as a localizing subcategory of $\mathcal{T}\opp$, that is, $\Coloc_\mathcal{T}(\mathcal{S})=\Loc_{\mathcal{T}\opp}(\mathcal{S}\opp)$, where $\mathcal{S}\opp:=\{X\opp\in\mathcal{T}\opp\mid X\in\mathcal{S}\}$. Since $\mathcal{T}\opp$ is perfectly generated,  $(\Loc(\mathcal{S}\opp), {(\mathcal{S}\opp)}^{\bot})$ is a hereditary torsion pair in $\mathcal{T}\opp$, that is, $({^{\bot}}\mathcal{S}, \Coloc(\mathcal{S}))$ is a hereditary torsion pair in $\mathcal{T}$. $\square$

\medskip
Let $F$, $G$ and $H$ be triangle endofunctors of $\mathcal{T}$. We say that a sequence of natural transformations
$$F\lraf{\tau} G\lraf{\eta} H\lraf{\sigma} F[1]:\;\;\mathcal{T}\lra \mathcal{T}$$ is \emph{exact} if for each $X\in \mathcal{T}$
the sequence $F(X)\lraf{\tau_X} G(X)\lraf{\eta_X} H(X)\lraf{\sigma_X} F(X)[1]$ is a triangle in $\mathcal{T}$.

\subsection{Cotorsion pairs and Gorestein-projective modules}

In the subsection we recall the definitions of cotorsion pairs in abelian categories and Gorenstein-projective modules over algebras.

Let $\mathcal{A}$ be an abelian category and $n\geq 1$ a natural number. Given a class $\mathcal{S}$ of objects in $\mathcal{A}$,  we define
$${^{\bot n}}\mathcal{S}:=\{X\in\mathcal{A}\mid \Ext_\mathcal{A}^n(X, S)=0\;\; \mbox{for}\;\; S\in\mathcal{S}\},
\;\;{^{\bot >0}}\mathcal{S}:=\bigcap_{n\geq 1}{^{\bot n}}\mathcal{S},
$$
$$
\mathcal{S}^{\bot n}:=\{X\in\mathcal{A}\mid \Ext_\mathcal{A}^n(S, X)=0\;\; \mbox{for}\; \;S\in\mathcal{S}\}\;\;\mbox{and}\;\;\mathcal{S}^{\bot >0}:=\bigcap_{n\geq 1}\mathcal{S}{^{\bot n}}.$$

\begin{Def}\label{def-cotor}
$(1)$ A pair $(\mathcal{U}, \mathcal{V})$ of full subcategories of $\mathcal{A}$  is called a \emph{cotorsion pair} in $\mathcal{A}$ if

$(i)$ $\mathcal{U}={^{\bot  1}}\mathcal{V}$ and $\mathcal{V}=\mathcal{U}^{\bot 1}$;

$(ii)$ For each object $X\in\mathcal{A}$, there are exact sequences
$ 0\to V_X\to U_X\lraf{\pi_X} X\to 0$ and $0\to X\lraf{\lambda_X}V^X\to U^X\to 0$
in $\mathcal{A}$ such that $U_X, U^X\in\mathcal{U}$ and $V_X, V^X\in\mathcal{V}$.

$(2)$ A cotorsion pair $(\mathcal{U}, \mathcal{V})$ in $\mathcal{A}$ is \emph{hereditary} if  $\mathcal{U}={^{\bot >0}}\mathcal{V}$ and $\mathcal{V}=\mathcal{U}^{\bot >0}$.
\end{Def}

Let $(\mathcal{U}, \mathcal{V})$ be a cotorsion pair in $\mathcal{A}$. Then $\pi_X$ is a right $\mathcal{U}$-approximation of $X$ and $\lambda_X$ is a left $\mathcal{V}$-approximation of $X$. Let $\mathcal{E}=\mathcal{U}\cap\mathcal{V}$. Then $U_X$ and $V^X$ are unique up to isomorphism in the quotient category $\mathcal{A}/\mathcal{E}$.
Further, the inclusion $\mathcal{U}/\mathcal{E}\to \mathcal{A}/\mathcal{E}$ has a right adjoint sending $X$ to $U_X$, while  the inclusion $\mathcal{V}/\mathcal{E}\to \mathcal{A}/\mathcal{E}$ has a left adjoint sending $X$ to $V^X$. Moreover, $X$ has a $\mathcal{U}$-resolution and a $\mathcal{V}$-coresolution in the following sense.

\begin{Def}\label{Co-Resolution}
Let $\mathcal{C}$ be a full subcategory of $\mathcal{A}$. A \emph{$\mathcal{C}$-resolution} of an object $U\in\mathcal{A}$ is a complex $\cpx{X}:=(X^i)_{i\in\mathbb{Z}}\in\C{\mathcal{C}}$ with $X^i=0$ for all $i\geq 1$, together with a chain map $\cpx{\pi}: \cpx{X}\to U$ (regarded as a stalk complex), such that $\cpx{\Hom}_\mathcal{A}(C, \cpx{\pi}):\cpx{\Hom}_\mathcal{A}(C, \cpx{X})\to \cpx{\Hom}_\mathcal{A}(C, U)$ is a quasi-isomorphism
for any $C\in\mathcal{C}$. In this case, $\cpx{\pi}$ is called a $\mathcal{C}$-resolution of $U$ for simplicity. Dually, a \emph {$\mathcal{C}$-coresolution} of $U$ can be defined.
\end{Def}

Let $A$ be an Artin algebra. We denote by $\Prj{A}$ and $\Inj{A}$ the full subcategories of $A\Modcat$ consisting of projective and injective $A$-modules, respectively.
Moreover, we write $\K{A}$, $\Kac{A}$ and $\D{A}$ for $\K{A\Modcat}$, $\Kac{A\Modcat}$ and $\D{A\Modcat}$, respectively.

An $A$-module $X$ is said to be \emph{Gorenstein-projective} if there is an exact complex
$\cpx{P}:\cdots\to P^{-2}\to P^{-1}\to P^{0}\to P^{1}\to P^{2}\to\cdots$ of projective $A$-modules
such that $X\simeq \Ker(P^0\to P^1)$ and the complex $\cpx{\Hom}_A(\cpx{P}, A)$ is exact.
The complex $\cpx{P}$ is called a \emph{complete projective resolution} of $X$. Dually, one can define Gorenstein-injective modules and complete injective coresolutions.  Let $\GP{A}$ and $\GI{A}$ be the full subcategories of all Gorenstein-projective and Gorenstein-injective $A$-modules, respectively. By \cite[Theorem X. 2.4]{BI}, $(\GP{A},\GP{A}^{\bot >0})$ and $({^{\bot >0}}\GI{A}, \GI{A})$ are hereditary cotorsion pairs in $A\Modcat$. Moreover, each $A$-module admits a minimal right $\GP{A}$-approximation and also a minimal left $\GI{A}$-approximation by \cite[Proposition 3.8(iv) and Corollary 6.8]{Bel}.

An Artin algebra $A$ is said to be \emph{virtually Gorenstein} (see \cite[Definition 8.1]{Bel}) if $\GP{A}^{\bot >0}={^{\bot >0}}\GI{A}$.
Virtually Gorenstein algebras include Gorenstein algebras and algebras of finite representation type, and are closed under derived equivalences and stable equivalences of Morita type.

By \emph{filtered colimits} of $A$-modules we mean colimits of filtered diagrams $I\to A\Modcat$ with $I$ an essentially small, filtered category.

\section{Recollements of stable module categories}\label{RSMC}

In this section, we investigate an arbitrary (not necessarily self-orthogonal) generator over a self-injective algebra, and construct two pairs of triangle endofunctors for the stable module category by means of the endomorphism algebra of the generator. When the generator is additionally self-orthogonal or $\Omega$-periodic, we show that these endofunctors coincidentally appear in a recollement of the stable module category determined by the generator (see Theorem \ref{Main result} and Corollary \ref{Self-orthogonal recollement}). This recollement will be restricted to the one of a relative Gorenstein stable category in the next section. Finally, we describe the category of Gorenstein-projective modules over the endomorphism algebra of a self-orthogonal generator in term of a relative Gorenstein category (see Lemma \ref{Gorenstein}), and provide equivalent characterizations for the endomorphism algebra to be virtually Gorenstein (see Proposition \ref{Limits}).

Throughout this section, let $A$ denote a \textbf{self-injective Artin algebra}. Then $\Stmc{A}$
is a triangulated category with a shift functor $[1]: \Stmc{A}\to \Stmc{A}$, given by the cosyzygy functor $\Omega^{-}_A$.
Clearly, $\Ext_A^n(X_1, X_2)\simeq\StHom_A(X_1, X_2[n])$ for all $n\geq 1$ and $X_1,X_2\in A\Modcat$. Let $q: A\Modcat\to\Stmc{A}$ be the canonical functor. By Lemma \ref{AJCP}(2) and its dual, $q$ preserves direct sums and direct products.
To emphasize objects in $\Stmc{A}$, the image of $X\in A\Modcat$ under $q$ is denoted by $\underline{X}$.

A full subcategory $\mathcal{U}$ of $A\Modcat$ is called a \emph {thick subcategory} if it is closed under direct summands in $A\Modcat$ and has
the \emph {two out of three property}: if two terms of an exact sequence $0\to X\to Y\to Z\to 0$ in $A\Modcat$ belong to $\mathcal{U}$, then so does the third. If $\mathcal{U}$ contains all projective $A$-modules, then $\mathcal{U}$ is a thick subcategory of $A\Modcat$ if and only if $\underline{\mathcal{U}}$ is a full triangulated subcategory of $\Stmc{A}$.

From now on, let $M = A\oplus M_0$ be a generator in $A\modcat$ such that $M_0$ contains no nonzero projective direct summands. As $M$ is finitely generated, we know $\Add(M)=\Prod(M)$. 

\subsection{Endomorphism algebras of generators over self-injective algebras}
In this subsection, we establish additive equivalences between the subcategories of the module category of a self-injective algebra and the ones of the module category of the endomorphism algebra of a generator over the algebra.

Let $\Lambda=\End_A(M)$ and $e$ be the idempotent element of $\Lambda$ such that $\Lambda e=\Hom_A(M,A)$. Then $A\simeq e\Lambda e$ as algebras. Let $S_e: \Lambda\Modcat\to A\Modcat, Y\mapsto eY$ be the \emph{Schur functor} determined by $e$.
This functor has a fully faithful left and right adjoint functors $F:=\Lambda e\otimes_A-$ and  $G:=\Hom_A(e\Lambda,-)$ from $A\Modcat$ to $\Lambda\Modcat$, respectively. Note that both $F$ and $G$ commute with direct products and direct sums. Due to $e\Lambda\simeq M$ as $A$-$\Lambda$-bimodules, the functor $G$ can be identified with $\Hom_A(M,-)$. Moreover, there is a canonical natural transformation $\delta:F\to G$
defined by
$$
\delta_X:\;F(X)\lra G(X),  \quad ae\otimes x\mapsto [eb\mapsto ebaex]
$$
where $X\in A\Modcat$, $a, b\in\Lambda$ and $x\in X$. For a full subcategory $\mathcal{X}\subseteq A\Modcat$
and $\mathcal{Y}\subseteq\Lambda\Modcat$
we define
$$\Kac{\mathcal{X}}:=\{\cpx{X}\in \K{\mathcal{X}}\mid \cpx{X}\in\Kac{A}\}\quad\mbox{and}\quad
\Keac{\mathcal{Y}}:=\{\cpx{Y}\in \K{\mathcal{Y}}\mid S_e(\cpx{Y})\in\Kac{A}\}.$$
Since $S_e$ is exact, $\Keac{\mathcal{Y}}$ consists of those $\cpx{Y}\in\K{\mathcal{Y}}$ with $H^i(\cpx{Y})\in (\Lambda/\Lambda e\Lambda)\Modcat$ for all $i\in\mathbb{Z}$.

Given an Artin algebra $B$, we denote by $\nu_B$ and $\nu_B^{-}$ the Nakayama functor $D(B)\otimes_B-:  B\Modcat\to B\Modcat$ and its right adjoint functor $\Hom_B(D(B),-)$, respectively. By restriction to $B\modcat$, we have natural isomorphisms $\nu_B\simeq D\circ\Hom_B(-,B)$ and $\nu_B^{-}\simeq \Hom_{B\opp}(-, B)\circ D$.
As $D(B)$ is finitely presented, $\Hom_B(D(B),-)$ commutes with filtered colimits, while $D(B)\otimes_B-$ commutes always with filtered colimits. Since a projective $B$-module is always a direct sum of finitely generated projective $B$-modules (see \cite [Theorem 27.11, p.306]{AF}, we can show by projective presentations of modules that, for a self-injective Artin algebra $B$, both $\nu_B$ and $\nu_B^{-}$ are exact and auto-equivalent on $B\Modcat$.

\begin{Lem}\label{Equivalence}
$(1)$ The restriction of $\delta$ to $\Prj{A}$ is a natural isomorphism, that is, $\delta_X$ is an isomorphism for any $X\in \Prj{A}$.

$(2)$ There are natural isomorphisms of additive functors:
$$D\circ F\circ \nu_A\simeq \Hom_A(-, M): A\Modcat\to\Lambda\opp\Modcat,$$
$$\nu_\Lambda\circ G\simeq F\circ \nu_A\;\;\mbox{and}\;\; \nu_\Lambda^-\circ F\simeq G\circ \nu_A^-:\;\; A\Modcat\to \Lambda\Modcat.$$
In particular, $F(\nu_A(M))\simeq D(\Lambda)$ as $\Lambda$-$\Lambda$-bimodules.

$(3)$ The functor $S_e$ restricts to equivalences of additive categories:
$$
\Add(\Lambda e)\lraf{\simeq} \Prj{A}, \;\; \Prj{\Lambda}\lraf{\simeq}\Add({_A}M)\;\;\mbox{and}\;\;\Inj{\Lambda}\lraf{\simeq} \Add(\nu_A(M)).
$$

$(4)$ Let $\mathcal{X}\subseteq A\Modcat$ and $\mathcal{Y}\subseteq\Lambda\Modcat$ be additive, full subcategories. If the functor $S_e$ induces an equivalence $\mathcal{Y}\lraf{\simeq}\mathcal{X}$ of additive categories, then there is an equivalence $\Keac{\mathcal{Y}}\lraf{\simeq}\Kac{\mathcal{X}}$
of triangulated categories.
\end{Lem}

{\it Proof.} $(1)$ holds because $\delta_A$ is an isomorphism and both $F$ and $G$ commute with direct sums.

$(2)$ We need the following result, its proof is left to the reader:

\smallskip
$(\ast)$ Let $A_1, A_2$ and $A_3$ be Artin algebras and let $F_1, F_2: A_1\Modcat\to A_2\Modcat$ be (covariant) additive functors. If $X$ is an $A_1$-$A_3$-bimodule, then $F_1(X)$ is an $A_2$-$A_3$-bimodule, where the right $A_3$-module structure is given by the composition of associated ring homomorphisms $A_3\to \End_{A_1}(X)$ and $\End_{A_1}(X)\to\End_{A_2}(F_1(X))$. Additionally, if $\eta: F_1\to F_2$ is a natural transformation, then $\eta_X: F_1(X)\to F_2(X)$ is a homomorphism of $A_2$-$A_3$-bimodules.

\smallskip
Note that $D\circ(F\circ \nu_A)=D\big(\Lambda e\otimes_A\nu_A(-)\big)\simeq \Hom_A(\nu_A(-), D(\Lambda e))$ and $D(\Lambda e)=D\Hom_A(M, A)\simeq D(A)\otimes_AM=\nu_A(M)$ as $A$-$\Lambda$-bimodules by $(\ast)$. Since $\nu_A: A\Modcat\to A\Modcat$ is an auto-equivalence, there are isomorphisms
$D\circ(F\circ \nu_A)\simeq\Hom_A\big(\nu_A(-), \nu_A(M)\big)\simeq\Hom_A(-, M): A\Modcat\to \Lambda\opp\Modcat$. This implies $D\circ(F\circ \nu_A)(M)\simeq \Lambda$ as $\Lambda\opp$-modules. Since $M$ is an $A$-$\Lambda$-bimodule, it follows from $(\ast)$ that $D\circ(F\circ \nu_A)(M)\simeq \Lambda$ as $\Lambda$-$\Lambda$-bimodules. This leads to $D(\Lambda)\simeq F(\nu_A(M))$ as $\Lambda$-$\Lambda$-bimodules.  Consequently, there are natural isomorphisms:
$$\nu_\Lambda\circ G=D(\Lambda)\otimes_\Lambda G(-)\simeq \Lambda e\otimes_AD(A)\otimes_AM\otimes_\Lambda\Hom_A(M,-)\simeq \Lambda e\otimes_AD(A)\otimes_A-= F\circ \nu_A,$$
$$\nu_\Lambda^-\circ F\simeq \Hom_\Lambda(\Lambda e\otimes_A\nu_A(M), \Lambda e\otimes_A-)\simeq\Hom_A(\nu_A(M), -)\simeq \Hom_A(M, \nu_A^-(-))=G\circ \nu_A^-.$$

$(3)$ The restrictions of $S_e$ and $G$ to the corresponding subcategories give the first and second equivalences since $S_e$ and $G$ commutes with direct sums. Now, we claim that $(F,S_e)$ restricts to equivalences between $\Inj{\Lambda}$ and $\Add(\nu_A(M))$.

In fact, $S_e(D(\Lambda_{\Lambda}))=eD(\Lambda)\simeq D(\Lambda e)=\nu_A(M)$. As $S_e$ commutes with direct sums and $\Inj{\Lambda}=\Add(D(\Lambda_\Lambda))=\Prod(D(\Lambda_\Lambda))$, the Schur functor
$S_e$ restricts to a functor from $\Inj{\Lambda}$ to $\Add(\nu_A(M))$. Moreover, by $(2)$, $F(\nu_A(M))\simeq D(\Lambda)$ as $\Lambda$-modules.
As $F$ commutes with direct sums, it also restricts to a functor from $\Add(\nu_A(M))$ onto $\Inj{\Lambda}$.
Recall that $(F, S_e)$ is an adjoint pair and $F$ is fully faithful. Thus $F:\Add(\nu_A(M))\to\Inj{\Lambda}$ is an equivalence with the quasi-inverse $S_e$.

$(4)$ Since $S_e:\mathcal{Y}\to\mathcal{X}$ is an equivalence, it induces an equivalence $\K{\mathcal{Y}}\simeq\K{\mathcal{X}}$.
Then $(4)$ follows from the definition of $\Keac{\mathcal{Y}}$. $\square$

\begin{Rem}\label{DDE}
$(1)$ Let ${\bf P}_1(\Lambda e)$ (respectively, ${\bf I}_1(\Lambda e)$) be the full subcategory of $\Lambda\Modcat$ consisting of all modules $Y$ such that there is an exact sequence $E_1\to E_0\to Y\to 0$ (respectively, $ 0\to Y\to E_0\to E_1$) in $\Lambda\Modcat$ with $E_0, E_1\in \Add(\Lambda e)$. By \cite[Lemma 3.1]{APT}, there are additive equivalences
$F: A\Modcat\to {\bf P}_1(\Lambda e)$ and $G:A\Modcat\to{\bf I}_1(\Lambda e).$
Thus, by Lemma \ref{Equivalence}(2), the functor $\nu_\Lambda$ restricts to an additive equivalence
${\bf I}_1(\Lambda e)\to {\bf P}_1(\Lambda e)$, with the inverse $\nu_\Lambda^{-}$.

$(2)$ Each $A$-module admits a minimal right and left $\Add(M)$-approximations.
This property will be used in Section \ref{relative stable category}.

Indeed, since $\Lambda$ is an Artin algebra, each $\Lambda$-module admits a projective cover. For any $A$-module $X$, let $f_X: Q\to G(X)$ be a projective cover of ${_\Lambda}G(X)$. Since $G$ is fully faithful and restricts to an equivalence from $\Add(M)$ to $\Prj{\Lambda}$ by Lemma \ref{Equivalence}(3), there is a homomorphism $r_X:M_X\to X$ of $A$-modules such that $M_X\in \Add(M)$ and $G(r_X)=f_X$. Then $r_X$ is a  minimal right $\Add(M)$-approximation of $X$.

Let $\overline{F}= F\circ\nu_A: A\Modcat\to \Lambda\Modcat$ be the composition of $\nu_A$ with $F$. Then $\overline{F}$ is fully faithful. By Lemma \ref{Equivalence}(3), $\overline{F}$ restricts to an equivalence from $\Add(M)$ to $\Inj{\Lambda}$. Thus, if  $g_X: \overline{F}(X)\to I$ is an injective envelop of $\overline{F}(X)$, then there is  a homomorphism $\ell_X: X\to M^X$  in $A\Modcat$ with $M^X\in\Add(M)$ such that $g_X=\overline{F}(\ell_X)$. One can check that $\ell_X$ is a minimal left $\Add(M)$-approximation of $X$.
\end{Rem}

Finally, we recall some notation and facts on homotopy and derived categories of rings.

Let $\K{\Lambda}_P$ (respectively, $\K{\Lambda}_I$) be the smallest full
triangulated subcategory of $\K{\Lambda}$ which

$(i)$ contains all bounded above (respectively, bounded below)
complexes of projective (respectively, injective) $\Lambda$-modules, and

$(ii)$ is closed under arbitrary direct sums (respectively, direct
products).

Note that $\K{\Lambda}_P\subseteq \K{\Prj{\Lambda}}$ and $\K{\Lambda}_I\subseteq \K{\Inj{\Lambda}}$. Moreover, the compositions
$\K{\Lambda}_P\hookrightarrow\K{\Lambda}\to\D{\Lambda}$ and $\K{\Lambda}_I\hookrightarrow\K{\Lambda}\to\D{\Lambda}$ are equivalences.
This means that, for any $\cpx{Y}\in\D{\Lambda}$, there exists a complex ${_p}\cpx{Y}\in \K{\Lambda}_P$ together
with a quasi-isomorphism ${_p}\cpx{Y}\to\cpx{Y}$. Dually, there is a complex ${_i}\cpx{Y}\in \K{\Lambda}_I$ together with a quasi-isomorphism
$\cpx{Y}\to{_i}\cpx{Y}$. Moreover, ${_p}\cpx{Y}$ and ${_i}\cpx{Y}$ are unique up to isomorphism in $\K{\Lambda}$. As usual,
${_p}\cpx{Y}$ and  ${_i}\cpx{Y}$ are called the \emph{projective resolution} and \emph{injective coresolution} of $\cpx{Y}$ in $\K{\Lambda}$, respectively.
For example, if $Y$ is a $\Lambda$-module, then ${_p}Y$ is a deleted projective resolution of $Y$.

Let $Q:\K{\Lambda}\to \D{\Lambda}$ be the localization functor.
Then $Q$ has a left adjoint $Q_\lambda$ and a right adjoint $Q_\rho$ defined by
$$Q_\lambda={_p}(-): \;\;\D{\Lambda}\to\K{\Lambda}\quad \mbox{and}\quad Q_\rho={_i}(-):\;\; \D{\Lambda}\to\K{\Lambda}.$$
It is known that $\Img(Q_\lambda)=\K{\Lambda}_P$ and $\Img(Q_\rho)=\K{\Lambda}_I$.

\begin{Lem} {\rm \cite{Kr1}}\label{Pro-Inj-Rec}
$(\K{\Lambda}_P, \Kac{\Prj{\Lambda}})$  and $(\Kac{\Inj{\Lambda}}, \K{\Lambda}_I)$ are hereditary torsion pairs in $\K{\Prj{\Lambda}}$ and $\K{\Inj{\Lambda}}$, respectively.  In other words, there are half recollements of triangulated categories:
$$
\xymatrix@C=1.2cm{\Kac{\Prj{\Lambda}}\ar[r]_-{I}
&\K{\Prj{\Lambda}}\ar[r]_-{Q}\ar_-{I_\lambda}@/_1.2pc/[l]
&\D{\Lambda}\ar@/_1.2pc/[l]_{Q_\lambda}}\;\; \mbox{and}\;\;
\xymatrix@C=1.2cm{\D{\Lambda}\ar[r]_-{Q_\rho}
&\K{\Inj{\Lambda}}\ar[r]_-{J_\rho}\ar_-{Q}@/_1.2pc/[l]
&\Kac{\Inj{\Lambda}}\ar@/_1.2pc/[l]_{J}}$$
where $I$ and $J$ are the inclusion functors, and $Q$ denotes restrictions of the localization functor.
\end{Lem}

\subsection{Construction of triangle endofunctors of stable module categories} \label{ECTF}

In the subsection we introduce two pairs of triangle endofunctors of stable module categories for self-injective algebras,
and then discuss two relevant thick subcategories of module categories.

Let $\epsilon: Q_\lambda \circ  Q\to {\rm Id}$ and $\eta: {\rm Id} \to I\circ I_\lambda$ be the counit and unit adjunctions with respect to the adjoint pairs $(Q_\lambda, Q)$ and $(I_\lambda, I)$ in Lemma \ref{Pro-Inj-Rec}, respectively, where ${\rm Id}$ denotes the identity functor.

\begin{Lem}\label{EDF}
$(1)$ There is a half recollement of triangulated categories:
$$
\xymatrix@C=1.2cm{\Kac{\Prj{\Lambda}}\ar[r]_-{I}
&\Keac{\Prj{\Lambda}}\ar[r]_-{Q_\lambda \circ  Q}\ar_-{I_\lambda}@/_1.2pc/[l]
&\Keac{\Prj{\Lambda}}\cap\K{\Lambda}_P\ar@/_1.3pc/[l]_{{\rm inc}}},$$
where $\rm{inc}$ is the inclusion functor.

$(2)$ There is an exact sequence of triangle endofunctors of $\Keac{\Prj{\Lambda}}$:
$$
Q_\lambda \circ  Q \lraf{\epsilon} {\rm Id} \lraf{\eta} I\circ I_\lambda\lra Q_\lambda \circ  Q[1]
$$
\end{Lem}

{\it Proof.} Since $\Kac{\Prj{\Lambda}}\subseteq\Keac{\Prj{\Lambda}}$ and $(\K{\Lambda}_P, \Kac{\Prj{\Lambda}})$ is a hereditary torsion pair in $\K{\Prj{\Lambda}}$, the pair $(\Keac{\Prj{\Lambda}}\cap\K{\Lambda}_P, \Kac{\Prj{\Lambda}})$ is a hereditary torsion pair in $\Keac{\Prj{\Lambda}}$. By Lemma \ref{Pro-Inj-Rec}, each object $\cpx{Y}\in\K{\Prj{\Lambda}}$ is endowed with a canonical triangle
$$
Q_\lambda \circ  Q(\cpx{Y})\lraf{\epsilon_{\cpx Y}} \cpx{Y} \lraf{\eta_{\cpx Y}} I\circ I_\lambda(\cpx{Y})\lra Q_\lambda \circ  Q(\cpx{Y})[1].
$$
If $\cpx{Y}\in\Keac{\Prj{\Lambda}}$, then it follows from  $I\circ I_\lambda(\cpx{Y})\in\Kac{\Prj{\Lambda}}$ and the exactness of $S_e$ that
$Q_\lambda\circ Q(\cpx{Y})\in\Keac{\Prj{\Lambda}}\cap\K{\Lambda}_P$. Thus $(1)$ and $(2)$ hold. $\square$

\medskip
Now let $\ell: \K{A}\to\K{\Inj{A}}$ be a \emph{left adjoint} of the inclusion $\K{\Inj{A}}\hookrightarrow\K{A}$.
By \cite[Corollary 1.3]{c}, $\ell$ is given by taking the total complexes of Cartan-Eilenberg injective  coresolutions of complexes over $A\Modcat$. Moreover, it has the following property.

\begin{Lem}\label{ell}
The functor $\ell$ restricts to a triangle functor $\ell_{ac}: \Kac{A}\to \Kac{\Inj{A}}$ and the composition
$$
\ell_M:\;\; \Kac{\Add(M)}\hookrightarrow \Kac{A}\lraf{\ell_{ac}} \Kac{\Inj{A}}=\Kac{\Prj{A}}.
$$
is a left adjoint of the inclusion $\Kac{\Prj{A}}\hookrightarrow \Kac{\Add(M)}$.
\end{Lem}

{\it Proof.} Clearly,  $(\Ker(\ell), \K{\Inj{A}})$ and $(\Kac{A}, \K{A}_I)$ are hereditary torsion pairs in $\K{A}$. Since $\K{A}_I\subseteq \K{\Inj{A}}$, we have $\Ker(\ell)\subseteq \Kac{A}$. This implies that $(\Ker(\ell), \Kac{\Inj{A}})$ is a hereditary torsion pair in $\Kac{A}$. Consequently, $\ell$ restricts to a functor $\Kac{A}\to \Kac{\Inj{A}}$.  Since ${_A}M$ is a generator, $\Kac{\Prj{A}}\subseteq\Kac{\Add(M)}$.
Now the second part of Lemma \ref{ell} holds because $\ell$ is a left adjoint of the inclusion $\K{\Inj{A}}\to\K{A}$. $\square$
\medskip

Let $X$ be an $A$-module. We denote by $\cpx{\pi}_X: \cpx{P}_X\to X$ and $\cpx{\lambda}_X: X\to \cpx{I}_X$
a minimal projective resolution and injective coresolution of $X$, respectively. Then there exists a triangle equivalence
$$S:\;\; \Stmc{A}\lraf{\simeq} \Kac{\Prj{A}},\quad X\mapsto S(X):=\Con(\cpx{\pi}_X\cpx{\lambda}_X).$$
This functor $S$ is called  the \emph{stabilization functor} of $A$ (for example, see \cite{Kr1}), while $S(X)$ is a complete projective resolution of $X$.
A quasi-inverse of $S$ is given by taking the $0$-th cocycle of complexes:
$$
Z^0: \Kac{\Prj{A}}\lra \Stmc{A},\quad \cpx{I}\mapsto Z^0(\cpx{I}) :=\Ker(I^0\to I^1).
$$
Further, let $\ell_M: \Kac{\Add(M)}\to \Kac{\Prj{A}}$ be the triangle functor defined in Lemma \ref{ell}, and let $\mu: \Keac{\Add(\Lambda e)}\to \Keac{\Prj{\Lambda}}$ be the inclusion induced from $\Add(\Lambda e)\subseteq \Prj{\Lambda}$.

By Lemmas \ref{EDF} and \ref{ell}, we can define a pair of triangle endofunctors of $\Stmc{A}$ by
$$
\Phi=Z^0\circ \ell_M\circ S_e\circ (Q_\lambda\circ Q)\circ \mu\circ G\circ S:\;\; \Stmc{A}\lra \Stmc{A},
$$
$$
\Psi=Z^0\circ \ell_M\circ S_e\circ (I\circ I_\lambda)\circ \mu\circ G\circ S:\;\;\;  \Stmc{A}\lra \Stmc{A}.
$$
They are illustrated by the following diagram
$$
\xymatrix{
\Stmc{A}\ar[d]_-{\rm Id}\ar_-{\Phi}@/_1.3pc/[d]\ar^-{\Psi}@/^1.4pc/[d]\ar[r]^-{S}_-{\simeq}& \Kac{\Prj{A}}\ar[r]^-{G}_-{\simeq}& \Keac{\Add(\Lambda e)}\ar^-{\mu}@{^{(}->}[r]&\Keac{\Prj{\Lambda}}
\ar[d]_-{\rm Id}\ar_-{Q_\lambda\circ Q}@/_1.3pc/[d]\ar^-{I\circ I_\lambda}@/^1.4pc/[d]\\
\Stmc{A}&\ar[l]_-{Z^0}^-{\simeq}\Kac{\Prj{A}}&\ar[l]_-{\ell_M} \Kac{\Add(M)} &\ar[l]_-{S_e}^-{\simeq} \Keac{\Prj{\Lambda}}
}
$$
where the equivalences of $G$ and $S_e$ follow from Lemma \ref{Equivalence}(3)-(4).
There is the natural isomorphism
$$Z^0\circ \ell_M\circ S_e\circ {\rm Id}\circ \mu\circ G\circ S\simeq {\rm Id}:\;\; \Stmc{A}\lra \Stmc{A}.$$
This follows from the equivalence
$ S_e:\;\Keac{\Add(\Lambda e)}\lraf{\simeq} \Kac{\Prj{A}}$
given by Lemma \ref{Equivalence}(3)-(4), together with the fact that the restriction of $\ell_M$ to $\Kac{\Prj{A}}$ is isomorphic to Id.

Dually, we can construct another pair $(\Phi',\Psi')$ of endofunctors of $\Stmc{A}$. Here, we only list some key points of this construction, and omit the details.

By Lemma \ref{Pro-Inj-Rec}, there is a half recollement of triangulated categories:
$$
\xymatrix@C=1.2cm{\Keac{\Inj{\Lambda}}\cap\K{\Lambda}_I\ar[r]^-{\rm inc}
&\Keac{\Inj{\Lambda}}\ar[r]^-{J_\rho }\ar_-{Q_\rho\circ Q}@/_1.2pc/[l]
&\Kac{\Inj{\Lambda}}\ar@/_1.3pc/[l]_{J}}.$$
This implies an exact sequence of endofunctors of $\Keac{\Inj{\Lambda}}$:
$$
J \circ J_\rho \lraf{\epsilon'} {\rm Id} \lraf{\eta'} Q_\rho\circ Q\lra J\circ J_\rho[1]:\;\; \Keac{\Inj{\Lambda}}\lra \Keac{\Inj{\Lambda}}
$$
where $\epsilon'$ and $\eta'$ are counit and unit adjunctions of the adjoint pairs $(J, J_\rho)$ and $(Q, Q_\rho)$, respectively (see Lemma \ref{Pro-Inj-Rec}). Let $r:\K{A}\to \K{\Prj{A}}$ be a right adjoint of the inclusion $\K{\Prj{A}}\to\K{A}$. It is given by taking the total direct product complexes of Cartan-Eilenberg projective resolutions of complexes over $A\Modcat$, due to \cite[Corollary 1.3]{c}.
Moreover, it restricts to a triangle functor $r_{ac}: \Kac{A}\to \Kac{\Prj{A}}$. Let
$$r_M: \;\;\Kac{\Add(\nu_A(M))}\hookrightarrow\Kac{A}\lraf{r_{ac}} \Kac{\Prj{A}}$$
be the composition of the inclusion with $r_{ac}$, and let $\mu':\Keac{\Add(\Lambda e)}\hookrightarrow\Keac{\Inj{\Lambda}}$
be the inclusion induced from $\Add(\Lambda e)\subseteq\Inj{\Lambda}$.
Now we define the two endofunctors $\Phi'$ and $\Psi' $ of $\Stmc{A}$ by
$$
\Phi'=Z^0\circ r_M\circ S_e\circ (Q_\rho\circ Q)\circ \mu'\circ F\circ S:\;\; \Stmc{A}\lra \Stmc{A},
$$
$$
\Psi'=Z^0\circ r_M\circ S_e\circ (J\circ J_\rho)\circ \mu'\circ F\circ S:\;\;\;  \Stmc{A}\lra \Stmc{A},
$$
which are illustrated by the diagram
$$
\xymatrix{
\Stmc{A}\ar[d]_-{\rm Id}\ar_-{\Psi'}@/_1.3pc/[d]\ar^-{\Phi'}@/^1.4pc/[d]\ar[r]^-{S}_-{\simeq}& \Kac{\Prj{A}}\ar[r]^-{F}_-{\simeq}& \Keac{\Add(\Lambda e)}\ar^-{\mu'}@{^{(}->}[r]&\Keac{\Inj{\Lambda}}
\ar[d]_-{\rm Id}\ar_-{J\circ J_\rho}@/_1.3pc/[d]\ar^-{Q_\rho\circ Q}@/^1.4pc/[d]\\
\Stmc{A}&\ar[l]_-{Z^0}^-{\simeq}\Kac{\Prj{A}}&\ar[l]_-{r_M} \Kac{\Add(\nu_A(M))} &\ar[l]_-{S_e}^-{\simeq} \Keac{\Inj{\Lambda}}.
}
$$
Similarly, there is a natural isomorphism of functors:
$$
Z^0\circ r_M\circ S_e\circ {\rm Id}\circ \mu'\circ F\circ S\simeq {\rm Id}: \;\; \Stmc{A}\lra \Stmc{A}.
$$
By construction, $\Phi$ and $\Psi$ commute with direct sums, while $\Phi'$ and $\Psi'$ commute with direct products. Clearly, if $M=A$, then $\Phi=\Phi'=0$
and $\Psi={\rm Id}=\Psi'$. In general, we have the result.

\begin{Prop}\label{Exact sequence}
There exist exact sequences of triangle endofunctors of $\Stmc{A}$:
$$
\Phi\lraf{\widetilde{\epsilon}} {\rm Id}\lraf{\widetilde{\eta}} \Psi\lra \Phi[1]\quad\mbox{and}\quad
\Psi'\lraf{\widetilde{\epsilon'}} {\rm Id}\lraf{\widetilde{\eta'}} \Phi'\lra \Psi'[1],
$$
where $\widetilde{\epsilon}$, $\widetilde{\eta}$, $\widetilde{\epsilon'}$ and $\widetilde{\eta'}$
are induced from $\epsilon$, $\eta$, $\epsilon'$ and $\eta'$, respectively.
\end{Prop}

Next, we investigate two full subcategories of $A\Modcat$ associated with the functors $\Psi$ and $\Psi'$.
$$
\mathscr{S}:=\{X\in A\Modcat\mid \Psi(\underline{X})=0\}\quad\mbox{and}\quad \mathscr{T}:=\{X\in A\Modcat\mid \Psi'(\underline{X})=0\}.
$$
Since $\Psi$ is a triangle functor commutating with direct sums, $\mathscr{S}$ is a thick subcategory of $A\Modcat$ containing projective modules and being closed under direct sums. Dually, $\mathscr{T}$ is a thick subcategory of $A\Modcat$ containing projective modules and being closed under direct products. Also, $\underline{\mathscr{S}}$ and $\underline{\mathscr{T}}$are full triangulated subcategories of $\Stmc{A}$.

Recall that an $A$-module $X$ is said to be \emph {$\Omega$-periodic} if  $\underline{\Omega_A^n(X)}\simeq \underline{X}$ in $\Stmc{A}$ for a positive integer $n$.

\begin{Lem}\label{Periodic}
$(1)$ If  $X\in\Add({_A}M)$ is $\Omega$-periodic, then $X\in \mathscr{S}$ and $\nu_A(X)\in\mathscr{T}$.

$(2)$ If $\Lambda$ has finite global dimension, then both $\Phi$ and $\Phi'$ are isomorphic to the identity functor, and therefore $\mathscr{S}=\mathscr{T}=A\Modcat$.
\end{Lem}
{\it Proof.} For an $A$-module $X$, we write $G\circ S(X)=\cpx{Y}:=(Y^n, d^n)_{n\in\mathbb{Z}}$ and $F\circ S(\nu_A(X))=\cpx{Z}:=(Z^n, h^n)_{n\in\mathbb{Z}}.$
Then $\cpx{Y}, \cpx{Z}\in \Keac{\Add(\Lambda e)}$.  Since the $\Lambda$-module
$\Lambda e$ is projective-injective,  both $Y^n$ and $Z^n$ are projective-injective.  As $G$ is left exact, $\Ker(d^n)=G(\Omega_A^{-n}(X))$. Dually,
$\Coker(h^n)=F\circ \Omega_A^{-(n+2)}(\nu_A(X))$ since $F$ is right exact.

$(1)$ Suppose $X\in\Add(M)$ and $X\simeq\Omega_A^s(X)$ in $\Stmc{A}$ for some $s\geq 1$. Then $X\simeq\Omega_A^{-sm}(X)$ in $\Stmc{A}$ for any $m\in\mathbb{N}$ and $G(X)\in\Prj{\Lambda}$ by Lemma \ref{Equivalence}(3). It follows that $\Ker(d^{sm})=G(\Omega_A^{-sm}(X))\simeq G(X)=\Ker(d^0)\in \Prj{\Lambda}$. Here, the isomorphism is regarded in the stable category.
Let  $\tau_{\leq sm}(\cpx{Y})$ be the subcomplex of $\cpx{Y}$:
$$
\cdots \lra Y^{-1}\lra \cdots\lra Y^{sm-2}\lra Y^{sm-1}\lra \Ker(d^{sm})\lra 0.
$$
Then $\tau_{\leq sm}(\cpx{Y})\in\Kf{\Prj{\Lambda}}$. Since
$\cpx{Y}$ is isomorphic in $\K{\Prj{\Lambda}}$ to the homotopy colimit of the sequence of inclusions:
$
\tau_{\leq 0}(\cpx{Y})\hookrightarrow\tau_{\leq s}(\cpx{Y})\hookrightarrow \tau_{\leq 2s}(\cpx{Y})\hookrightarrow \cdots\hookrightarrow\tau_{\leq sm}(\cpx{Y})\hookrightarrow \cdots,
$
there is a canonical triangle in $\K{\Prj{\Lambda}}$:
$$\bigoplus_{m=0}^{\infty}\tau_{\leq sm}(\cpx{Y})\lra \bigoplus_{m=0}^{\infty}\tau_{\leq sm}(\cpx{Y})\lra \cpx{Y}\lra\bigoplus_{m=0}^{\infty}\tau_{\leq sm}(\cpx{Y}).  $$
As $\Kf{\Prj{\Lambda}}\subseteq \K{\Lambda}_P$ and
$\K{\Lambda}_P$ is closed under direct sums in
$ \K{\Prj{\Lambda}}$, we get $\cpx{Y}\in\K{\Lambda}_P$.
Since $(\K{\Lambda}_P, \Kac{\Prj{\Lambda}})$ is a hereditary torsion pair in $\K{\Prj{\Lambda}}$, we obtain $Q_\lambda\circ Q(\cpx{Y})\simeq\cpx{Y}$ and
$I\circ I_\lambda(\cpx{Y})=0$. Thus  $ \Psi(\underline{X})=0$ and $X\in \mathscr{S}.$

Set $U:=\nu_A(X)$. It follows from $X\in\Add(M)$ and Lemma \ref{Equivalence}(3) that $F(U)\in \Inj{\Lambda}$.  Since $X\simeq\Omega_A^{-s}(X)$ for some $s\ge 1$ and $\nu_A$ is an auto-equivalence on $A\Modcat$, we have $U\simeq\Omega_A^{sm}(U)$ for all $m\in\mathbb{N}$. Thus $\Coker(h^{-2-sm})=F\circ \Omega_A^{sm}(U)\simeq F(U)=\Coker(h^{-2})$ and the truncated quotient complex of $\cpx{Z}$
$$
\tau_{\geq -1-sm}(\cpx{Z}):\quad \cdots\lra 0\lra \Coker(h^{-2-sm})\lra Z^{-sm}\lra Z^{-sm+1}\lra\cdots
$$
is in $\Kz{\Inj{\Lambda}}$. Dually, $\cpx{Z}$ is isomorphic in $\K{\Inj{\Lambda}}$ to a homotopy limit of
the canonical surjections:
$$
\cdots \twoheadrightarrow\tau_{\geq -1-sm}(\cpx{Z})\twoheadrightarrow\cdots\twoheadrightarrow\tau_{\geq -1-2s}(\cpx{Z})\twoheadrightarrow
\tau_{\leq -1-s}(\cpx{Y})\twoheadrightarrow\tau_{\leq -1}(\cpx{Y}).
$$ We can show $\cpx{Z}\in\K{\Lambda}_I$. Hence $Q_\rho\circ Q(\cpx{Z})\simeq\cpx{Z}$ and $J\circ J_\rho(\cpx{Z})=0$.
Thus $\Psi'(\underline{U})=0$ and $U\in \mathscr{T}.$

$(2)$ Suppose $\Lambda$ has finite global dimension. Then $\Kac{\Prj{\Lambda}}=0=\Kac{\Inj{\Lambda}}$, $\K{\Prj{\Lambda}}=\K{\Lambda}_P$ and $\K{\Inj{\Lambda}}=\K{\Lambda}_I$. Hence $Q_\lambda\circ Q$ and $Q_\rho\circ Q$ are naturally isomorphic to Id. This implies that $\Phi$ and $\Phi'$ are naturally isomorphic to Id, while $\Psi$ and $\Psi'$ are zero functor. $\square$

\begin{Prop}\label{S-cat} Let $X$ be an $A$-module. Then the following hold.

$(1)$ If $X\in M^{\bot >0}$, then $X\in \mathscr{S}$ if and only if for any $A$-module $Y$,
$\Hom_{\K{A}}(\Con(\cpx{\mu}_X), \Con(\cpx{\pi}_Y))=0$, where
$\cpx{\mu}_X: \cpx{M}_X\to X$ is an $\Add(M)$-resolution of $X$.

$(2)$ If $X\in {^{\bot >0}}M$, then $\nu_A(X)\in \mathscr{T}$ if and only if for any $A$-module $Y$, $\Hom_{\K{A}}(\Con(\cpx{\lambda}_Y), \Con(\cpx{\sigma}_X))=0$ where $\cpx{\sigma}_X: X\to \cpx{M}_X$ is an $\Add(M)$-coresolution of $X$.
\end{Prop}

{\it Proof.} $(1)$  Let $\cpx{\pi}:=\cpx{\pi}_X$, $\cpx{\lambda}:=\cpx{\lambda}_X$ and $\cpx{\mu}:=\cpx{\mu}_X$. Then $S(X)=\Con(\cpx{\pi}\cpx{\lambda})$. Since $\cpx{\mu}$ is an $\Add(M)$-resolution of $X$ (see Definition \ref{Co-Resolution}),  $G(\cpx{\mu}): G(\cpx{M}_X)\to G(X)$ is a quasi-isomorphism; equivalently,
 $\Con(G(\cpx{\mu}))$ is exact. Moreover, $\Con(G(\cpx{\mu}))=G(\Con(\cpx{\mu}))$. As ${_A}M$ is a generator, $\Con(\cpx{\mu})$ is exact.
By $\cpx{P}_X\in\Kf{\Prj{A}}$, the identity map of $X$ can be lifted to a unique morphism $\cpx{h}: \cpx{P}_X\to\cpx{M}_X$ in $\K{A}$ such that $\cpx{h}\cpx{\mu}=\cpx{\pi}$. This implies $G(\cpx{h})G(\cpx{\mu}\cpx{\lambda})=G(\cpx{\pi}\cpx{\lambda})$. By the octahedral axiom of triangulated category, there exists a triangle in $\Keac{\Prj{\Lambda}}$:
$$(\ast):\;\; \Con(G(\cpx{h}))\lra \mu\circ G\circ S(X)\lra\Con(G(\cpx{\mu}\cpx{\lambda}))\lra \Con(G(\cpx{h}))[1]$$
where $\mu\circ G\circ S(X)=\Con(G(\cpx{\pi}\cpx{\lambda}))$. Since $X\in M^{\bot >0}$, the morphism $G(\cpx{\lambda}):G(X)\to G(\cpx{I}_X)$ is a quasi-isomorphism.
Thus $G(\cpx{\mu}\cpx{\lambda})=G(\cpx{\mu})G(\cpx{\lambda})$ is a composition of two quasi-isomorphisms. This means
$\Con(G(\cpx{\mu}\cpx{\lambda}))\in \Kac{\Prj{\Lambda}}$.
It follows from $\Con(G(\cpx{h}))\in\mathscr{K}^{-}(\Prj{\Lambda})\subseteq \K{\Lambda}_P$ that
$$Q_\lambda \circ  Q(\mu\circ G\circ S(X))\simeq \Con(G(\cpx{h}))\;\;\mbox{and}\;\; I\circ I_\lambda(\mu\circ G\circ S(X))\simeq \Con(G(\cpx{\mu}\cpx{\lambda})).$$
Since the composition of $G$ with $S_e$ is isomorphic to the identity functor, we apply $S_e$ to the triangle $(\ast)$ and get another triangle $\Con(\cpx{h})\to  S(X)\to \Con(\cpx{\mu}\cpx{\lambda})\to \Con(\cpx{h})[1]$ in $\Kac{\Add(M)}$. Further, by applying the functor $\ell_M$ (see Lemma \ref{ell}) to this triangle, we are led to a  triangle in $\Kac{\Prj{A}}$:
$$\ell_M(\Con(\cpx{h}))\lra  \ell_M\circ S(X)\lra \ell_M(\Con(\cpx{\mu}\cpx{\lambda}))\lra \ell_M(\Con(\cpx{h}))[1].$$
Clearly, $\ell_M\circ S(X)=S(X)$ since $S(X)\in\Kac{\Prj{A}}$. From $Z^0\circ S(X)\simeq X$, we obtain a triangle in
$\Stmc{A}$:
$$Z^0\circ \ell_M(\Con(\cpx{h}))\lra X\lra Z^0\circ \ell_M(\Con(\cpx{\mu}\cpx{\lambda}))\lra Z^0\circ \ell_M(\Con(\cpx{h}))[1].$$ Thus
$\Phi(X)=Z^0\circ \ell_M(\Con(\cpx{h}))$ and $\Psi(X)=Z^0\circ \ell_M(\Con(\cpx{\mu}\cpx{\lambda})).$
 Note that $Z^0: \Kac{\Prj{A}}\to \Stmc{A}$ is an equivalence. Hence $\Psi(X)=0$ (equivalently, $X\in\mathscr{S}$) if and only if $\ell_M(\Con(\cpx{\mu}\cpx{\lambda}))=0$ if and only if $\Hom_{\K{A}}(\Con(\cpx{\mu}\cpx{\lambda}), \cpx{Q})=0$ for any $\cpx{Q}\in\Kac{\Prj{A}}$.

Since $\Con(\cpx{\lambda})\in \mathscr{K}^{+}_{\rm ac}(A)$ and $\mathscr{K}^{+}_{\rm ac}(A)\subseteq {^\bot}\K{\Inj{A}}$,
there holds $\Hom_{\K{A}}(\Con(\cpx{\lambda}), \cpx{Q}[n])=0$ for any $n\in\mathbb{Z}$. Applying $\Hom_{\K{A}}(-, \cpx{Q}[n])$ to the triangle $$\Con(\cpx{\mu})\to \Con(\cpx{\mu}\cpx{\lambda})\to \Con(\cpx{\lambda})\to \Con(\cpx{\mu})[1]$$ in $\Kac{A}$, we get $\Hom_{\K{A}}(\Con(\cpx{\mu}\cpx{\lambda}), \cpx{Q})\simeq \Hom_{\K{A}}(\Con(\cpx{\mu}), \cpx{Q})$. Now, let $Y$ be the kernel of the $0$-th differential of $\cpx{Q}$. Taking the canonical truncation on $\cpx{Q}$ at degree $0$, we obtain a subcomplex $\tau_{\leq 0}\cpx{Q}$ of $\cpx{Q}$, which is acyclic and isomorphic to $\Con(\cpx{\pi}_Y)$ in $\K{A}$. Since the inclusion $\tau_{\leq 0}\cpx{Q}\subseteq \cpx{Q}$ induces $\Hom_{\K{A}}(\Con(\cpx{\mu}), \tau_{\leq 0}\cpx{Q})\simeq \Hom_{\K{A}}(\Con(\cpx{\mu}), \cpx{Q})$, it follows that $\Hom_{\K{A}}(\Con(\cpx{\mu}), \Con(\cpx{\pi}_Y))\simeq \Hom_{\K{A}}(\Con(\cpx{\mu}\cpx{\lambda}), \cpx{Q})$. So $X\in\mathscr{S}$ if and only if
$\Hom_{\K{A}}(\Con(\cpx{\mu}), \Con(\cpx{\pi}_Y))=0$. This shows $(1)$.

$(2)$ Let $\nu:=\nu_A$, $Z:=\nu(X)$ and $\cpx{\sigma}:=\cpx{\sigma}_X$. Then
$F\circ S(Z)=F\big(\Con(\nu(\cpx{\pi}\cpx{\lambda}))\big)=\Con\big(F(\nu\cpx{\pi})F(\nu\cpx{\lambda})\big).$
Since $X\in {^{\bot >0}}M$ and $\nu$ is an auto-equivalence of $A\Modcat$, there hold $D\Tor_i^A(\Lambda e, Z)\simeq \Ext_A^i(Z, D(\Lambda e))\simeq \Ext_A^i(Z, \nu M)\simeq\Ext_A^i(X,M)=0$ for any $i\geq 1$. Hence $F(\nu\cpx{\pi} )$ is a quasi-isomorphism, that is, $\Con(F(\nu\cpx{\pi}))$ is exact.  Thus $F\circ S(Z)\simeq\Con(F(\nu\cpx{\lambda}))$ in $\D{\Lambda}$, where $\nu\cpx{\lambda}: Z\to \nu(\cpx{I}_X)$ is an injective coresolution of $Z$.  To calculate  $Q_\rho\circ Q\circ \mu'(\Con(F(\nu\cpx{\lambda})))$ in $\Keac{\Inj{\Lambda}}$, we show that $F(\nu \cpx{\sigma}): F(Z)\to F\circ \nu(\cpx{M}_X)$ is an injective coresolution of $F(Z)$, and then replace $F(Z)$ by its deleted injective coresolution.

In fact, by the proof of Lemma \ref{Equivalence}(3), the adjoint pair $(F, S_e)$ induces an equivalence $\Add(\nu M)\lraf{\simeq} \Inj{\Lambda}$. This implies $F(\nu(\cpx{M}_X))\in\mathscr{K}^{+}(\Inj{\Lambda})$. It remains to show that
$F(\nu\cpx{\sigma})$ is a quasi-isomorphism. Since $D:\Lambda\Modcat\to \Lambda\opp\Modcat$ is exact and detects zero objects, we only need to show that $DF(\nu\cpx{\sigma})$ is a quasi-isomorphism. However, by Lemma \ref{Equivalence}(2), $DF(\nu \cpx{\sigma})\simeq \Hom_A(\cpx{\sigma}, M)$ which is a quasi-isomorphism by the construction of $\sigma$. Thus $F(\nu\cpx{\sigma})$ is an injective coresolution of $F(Z)$.

Now, let $\cpx{f}: \cpx{M}_X\to\cpx{I}_X$ be a chain map which lifts the identity map of $X$. Then there is a canonical triangle
$$
\Con(F\circ \nu(\cpx{\pi}\cpx{\sigma}))\lra \mu'\circ F\circ S(Z)\lra \Con(F\circ \nu(\cpx{f}))
\lra \Con(F\circ \nu(\cpx{\pi}\cpx{\sigma}))[1]
$$
 in $\Keac{\Inj{\Lambda}}$, where the first term lies in $\Kac{\Inj{\Lambda}}$ and the third one lies in $\mathscr{K}^+(\Inj{\Lambda})$.
Thus $$J\circ J_\rho \circ \mu'\circ F\circ S(Z)\simeq \Con(F\circ \nu(\cpx{f}))\quad\mbox{and}\quad Q_\rho\circ Q\circ \mu'\circ F\circ S(Z)\simeq \Con(F\circ \nu(\cpx{\pi}\cpx{\sigma})).$$
Since the composition of $F$ with $S_e$ is isomorphic to the identity functor, it follows that
$$
\Psi'(Z)=Z^0\circ r_M(\Con(\nu(\cpx{\pi}\cpx{\sigma})))\quad\mbox{and}\quad \Phi'(Z)=Z^0\circ r_M(\Con(\nu\cpx{f})).
$$
Dually,  by the equivalence of $Z^0$ and the inclusion $\mathscr{K}^{+}_{\rm ac}(A)\subseteq \K{\Prj{A}}^{\bot}$,
we can show that $Z\in \mathscr{T}$ if and only if $\Hom_{\K{A}}(\Con(\cpx{\lambda}_U), \Con(\nu\cpx{\sigma}))=0$ for any $A$-module $U$.
Since $\nu$ is an auto-equivalence of $A\Modcat$, we have $\Hom_{\K{A}}(\Con(\cpx{\lambda}_U), \Con(\nu\cpx{\sigma}))\simeq
\Hom_{\K{A}}(\Con(\cpx{\lambda}_U), \nu(\Con(\cpx{\sigma})))\simeq\Hom_{\K{A}}(\Con(\cpx{\lambda}_Y), \Con(\cpx{\sigma}))$
with $Y:=\nu^{-}(U)$. Thus $(2)$ holds. $\square$

\medskip
The following result is a consequence of Lemma \ref{Periodic}(1) and Proposition \ref{S-cat}.

\begin{Koro}\label{SO}
If ${_A}M$ is  self-orthogonal or $\Omega$-periodic, then $M\in \mathscr{S}$ and $\nu_A(M)\in\mathscr{T}$.
\end{Koro}

\subsection{Recollements of stable module categories induced by self-orthogonal modules} \label{ECTF2}
In this subsection, we apply the triangle endofunctors in Section \ref{ECTF} to construct a recollement of the stable module category of a self-injective algebra from a generator with conditions that are  satisfied for a self-orthogonal or an $\Omega$-periodic generator.

As is known, $\Stmc{A}$ is compactly generated and the inclusion of $A\modcat\to A\Modcat$ induces a triangle equivalence from $\stmc{A}$ to the full subcategory of $\Stmc{A}$ consisting of all compact objects.
Since ${_A}M$ is finitely generated, $M$ is compact in $\Stmc{A}$.
For a set $\Delta$ of integers, let
$$\underline{M}^{\bot\Delta}:=\{X\in\Stmc{A}\mid \StHom_A(M, X[n])=0\;\;\mbox{for any}\;\;n\in\Delta\},$$
$${^{\bot\Delta}}\underline{M}:=\{X\in\Stmc{A}\mid \StHom_A(X, M[n])=0\;\;\mbox{for any}\;\;n\in\Delta\}.$$
For simplicity, we write $\underline{M}^{\bot}$ and ${^\bot}\underline{M}$ for $\underline{M}^{\bot\mathbb{Z}}$ and ${^{\bot\mathbb{Z}}}\underline{M}$, respectively.
Then $\underline{M}^{\bot}$ is a full triangulated subcategory of $\Stmc{A}$ closed under direct sums and direct products.

\begin{Lem} \label{Vanishing-F}
Let $X$ be an $A$-module.

$(1)$ $X\in\underline{M}^{\bot}$ if and only if $G\circ S(X)$ (equivalently, $F\circ S(X)$) is an exact complex.

$(2)$ If $X\in\underline{M}^{\bot}$, then $\Phi(X)=0$ and $\Phi'(X)=0$.
\end{Lem}
{\it Proof.} By Lemma \ref{Equivalence}(1), $F$ and $G$ are naturally isomorphic on $\Prj{A}$. Thus $G\circ S(X)\simeq F\circ S(X)$ as complexes. So we show (1) for $G\circ S(X)$. Let $0\to X_1\to P\to X_0\to 0$ be an exact sequence of $A$-modules with $P\in\Prj{A}$.
Then $0\to G(X_1)\to G(P)\to G(X_0)\to 0$ is exact if and only if $\StHom_A(M, X_0)=0$.
This implies $(1)$.   Moreover, $(2)$ follows from $(1)$ and the definitions of $\Phi$ and $\Phi'$. $\square$

\medskip
We need the following result which is concluded from the Auslander-Reiten formula (see \cite{AR1}).

\begin{Lem}\label{AR}
If $X\in A\modcat$, then there is a natural isomorphism
$$D\StHom_A(X,-)\simeq\StHom_A(-, \nu_A(X)[-1]):\;\; \Stmc{A}\lra \underline{{\rm End}}_A(X)\opp\Modcat. $$
\end{Lem}

{\it Proof.} Since $X$ is finitely generated, it follows from  \cite[Proposition 2.2]{AR1} that
$$D\StHom_A(X,-)\simeq\Ext_A^1(-, D{\rm Tr}(X)):\;\;\Stmc{A}\lra \underline{{\rm End}}_A(X)\opp\Modcat $$
where $D{\rm Tr}$ is the Auslander-Reiten translation on $\stmc{A}$. Since $A$ is self-injective,
$D{\rm Tr} \simeq \Omega^2_A \circ \nu_A$ as functors on $\stmc{A}$. This implies that
$$\begin{array}{ll} \Ext_A^1(-, D{\rm Tr}(X)) & \simeq \Ext_A^1(-, \Omega^2_A\circ \nu_A(X)) \simeq \StHom_A(\Omega_A(-), \Omega^2_A\circ \nu_A(X))\\ & \simeq \StHom_A(-, \Omega_A\circ \nu_A(X))=\StHom_A(-, \nu_A(X)[-1])\end{array}$$on $\Stmc{A}$.  Thus Lemma \ref{AR} holds.
$\square$

\medskip
To construct recollements of $\Stmc{A}$ from compact objects, we establish a result on torsion pairs.

\begin{Lem}\label{Endofunctor}
$(1)$ $\big(\Loc(\underline{M}),\, \underline{M}^{\bot}\big)$ and $\big(\underline{M}^{\bot}, \Coloc(\underline{\nu_A(M)})\big)$ are hereditary torsion pairs in $\Stmc{A}$.

$(2)$ If ${_A}M\in \mathscr{S}$, then $\underline{\mathscr{S}}=\Img(\Phi)=\Loc(\underline{M})$ and $\Img(\Psi)=\underline{M}^{\bot}$.

$(3)$ If $\nu_A(M)\in \mathscr{T}$, then $\underline{\mathscr{T}}=\Img(\Phi')=\Coloc(\underline{\nu_A(M)})$ and $\Img(\Psi')=\underline{M}^{\bot}$.

\end{Lem}

{\it Proof.} (1) Since ${_A}M$ is finitely generated, we have $\underline{M}^{\bot}={^\bot}\underline{\nu_A(M)}\subseteq\Stmc{A}$ by Lemma \ref{AR}.
Clearly, $\Stmc{A}$ is compactly generated by simple modules because each $A$-module has a radical series of length less than or equal the Loewy length of $A$. Thus $(1)$ holds by Proposition \ref{CPG}.

$(2)$ Let $\mathscr{X}:=\Loc(\underline{M})$ and $\mathscr{Y}:=\underline{M}^{\bot}$. Suppose $M\in \mathscr{S}$. Then $\mathscr{X}\subseteq\underline{\mathscr{S}}$ because $\underline{\mathscr{S}}$ is a full triangulated subcategory of $\Stmc{A}$ closed under direct sums. By Lemma \ref{Vanishing-F}(2), $\mathscr{Y}\subseteq\Ker(\Phi)$.

Let $N\in A\Modcat$. By $(1)$, up to isomorphism, there is a unique triangle
$X_N\to N\to Y^N\to X_N[1]$ in $\Stmc{A}$ such that $X_N\in\mathscr{X}$ and $Y^N\in\mathscr{Y}$.
This yields $\Psi(X_N)=0=\Phi(Y^N)$. Now, we apply Proposition \ref{Exact sequence} to the triangle and obtain the commutative diagram in $\Stmc{A}$:
$$
\xymatrix{
\Phi(X_N)\ar[r]^-{\simeq}\ar[d]^-{\simeq}_-{\widetilde{\epsilon}_{X_N}}& \Phi(N)\ar[r]\ar[d]_-{\widetilde{\epsilon}_N}&
0\ar[r]\ar[d] & \Phi(X_N)[1]\ar[d]^-{\simeq}\\
X_N\ar[r]\ar[d] & N\ar[r]\ar[d]_-{\widetilde{\eta}_N} & Y^N\ar[r]\ar[d]^-{\simeq}_-{\widetilde{\eta}_{Y^N}}& X_N[1]\ar[d]\\
0\ar[r]\ar[d]& \Psi(N)\ar[r]^-{\simeq}\ar[d]& \Psi(Y^N)\ar[r]\ar[d]& \Psi(X_N)[1]\ar[d]\\
\Phi(X_N)[1]\ar[r]^-{\simeq}& \Phi(N)[1]\ar[r]& 0\ar[r] &\Phi(X_N)[2]
}
$$
where all rows and columns are triangles. Thus $\Phi(N)\simeq X_N\in\mathscr{X}$ and $\Psi(N)\simeq Y^N\in\mathscr{Y}$.
This implies $\Img(\Phi)\subseteq \mathscr{X}$ and $\Img(\Psi)\subseteq\mathscr{Y}$. The inclusions in other direction follows from Proposition \ref{Exact sequence}. Thus $\Img(\Phi)=\mathscr{X}$ and $\Img(\Psi)=\mathscr{Y}$. Since $\mathscr{X}\subseteq\underline{\mathscr{S}}\subseteq\Img(\Phi)$ and $\mathscr{Y}\subseteq \underline{\mathscr{T}}\subseteq \Img(\Phi)$, we have $\underline{\mathscr{S}}=\Img(\Phi)$ and $\Img(\Psi)=\underline{M}^{\bot}$.

(3) Similar, we can show $(3)$ by Lemma \ref{Vanishing-F}(2) and the pair $\big(\underline{M}^{\bot}, \Coloc(\underline{\nu_A(M)})\big)$ in $(1)$. $\square$

\medskip
A consequence of Lemma \ref{Endofunctor} and Corollary \ref{SO} is the following.

\begin{Koro}\label{localizing}
Suppose that ${_A}M$ is self-orthogonal or $\Omega$-periodic. Then $\mathscr{S}$ is the smallest thick subcategory of $A\Modcat$ containing $M$ and being closed under direct sums, while $\mathscr{T}$ is the smallest thick subcategory of $A\Modcat$ containing $\nu_A(M)$ and being closed under direct products.
\end{Koro}

Now, we are in position to prove the following main result of this section.

\begin{Theo}\label{Main result}
Suppose $A$ is a self-injective Artin algebra and $_AM$ is generator in $A\modcat$. If ${_A}M\in \mathscr{S}$ and  $\nu_A(M)\in \mathscr{T}$, then there exists a recollement of triangulated categories:
$$\xymatrix@C=1.2cm{\underline{M}^{\bot}\ar[r]^-{{\rm inc}}
&\Stmc{A}\ar_-{\widetilde{\Psi'}}@/^1.2pc/[l]\ar[r]^-{\widetilde{\Phi}}\ar@/^1.2pc/[l]\ar_-{\widetilde{\Psi}}@/_1.2pc/[l]
&\Loc(\underline{M})\ar_-{\Phi''}@/^1.2pc/[l]\ar@/^1.2pc/[l]\ar@/_1.2pc/[l]_{{\rm inc}}}
$$
such that
$$\Phi={\rm inc}\circ\widetilde{\Phi},\quad \Psi={\rm inc}\circ\widetilde{\Psi},\quad\Psi'={\rm inc}\circ\widetilde{\Psi'}
\quad{\rm and}\quad\Phi''=\Phi'\circ{\rm inc}.
$$
Moreover, the functor $\Phi''$ restricts to a triangle equivalence $\Loc(\underline{M})\lraf{\simeq}\Coloc(\underline{\nu_A(M)}).$
\end{Theo}

{\it Proof.} Suppose ${_A}M\in \mathscr{S}$. By Lemma \ref{Endofunctor}(2), we have the factorisation of $\Phi$ and $\Psi$:
$$\Phi: \;\;\Stmc{A}\lraf{\widetilde{\Phi}} \Loc(\underline{M})\hookrightarrow\Stmc{A}\quad\mbox{and}\quad\Psi:\;\;
\Stmc{A}\lraf{\widetilde{\Psi}} \underline{M}^{\bot}\hookrightarrow\Stmc{A}.$$
By Lemma \ref{Endofunctor}(1), $\big(\Loc(\underline{M}),\, \underline{M}^{\bot}\big)$ is a hereditary torsion pair in $\Stmc{A}$. Then the proof of Lemma \ref{Endofunctor}(2) together with \cite[Chapter I, Prop. 2.3]{BI} implies that $\widetilde{\Phi}$ is a right adjoint of the inclusion $\Loc(\underline{M})\to\Stmc{A}$ and that $\widetilde{\Psi}$ is a left adjoint of the inclusion $\underline{M}^{\bot}\to\Stmc{A}$.

Suppose $\nu_A(M)\in \mathscr{T}$. Dually, from the torsion pair $\big(\underline{M}^{\bot}, \Coloc(\underline{\nu_A(M)})\big)$ in Lemma  \ref{Endofunctor}(1) and  from Lemma \ref{Endofunctor}(3),  we obtain the factorisations of $\Phi'$ and $\Psi'$:
$$\Phi': \;\;\Stmc{A}\lraf{\widetilde{\Phi'}} \Coloc(\underline{\nu_A(M)})\hookrightarrow\Stmc{A}\quad\mbox{and}\quad\Psi':\;\;
\Stmc{A}\lraf{\widetilde{\Psi'}} \underline{M}^{\bot}\hookrightarrow\Stmc{A}$$
such that $\widetilde{\Phi'}$ is a left adjoint of the inclusion $\Coloc(\underline{\nu_A(M)})\to\Stmc{A}$ and
$\widetilde{\Psi'}$ is a right adjoint of the inclusion $\underline{M}^{\bot}\to\Stmc{A}$. Recall that there is a correspondence between TTF (torsion-torsionfree) triples and recollements of triangulated categories (see, for example, \cite[Chapter I. 2]{BI} or \cite[Section 2.3]{xc1}). Thus
Theorem \ref{Main result} follows from Lemma \ref{Endofunctor}(1) and \cite[Lemma 2.6]{xc1}. $\square$
\medskip

Combining Theorem \ref{Main result} with Corollary \ref{SO}, we obtain the corollary.

\begin{Koro}\label{Self-orthogonal recollement} Let $A$ be a self-injective algebra.
If ${_A}M$ is self-orthogonal or $\Omega$-periodic, then there exists a recollement of triangulated categories:
$$
\xymatrix@C=1.2cm{\underline{M}^{\bot}\ar[r]^-{{\rm inc}}
&\Stmc{A}\ar_-{\widetilde{\Psi'}}@/^1.2pc/[l]\ar[r]^-{\widetilde{\Phi}}\ar@/^1.2pc/[l]\ar_-{\widetilde{\Psi}}@/_1.2pc/[l]
&\Loc(\underline{M})\ar_-{\Phi''}@/^1.2pc/[l]\ar@/^1.2pc/[l]\ar@/_1.2pc/[l]_{{\rm inc}},}
$$
\smallskip
in which the functors are the same as the ones in  Theorem \ref{Main result}.
\end{Koro}

Later we will see that the above recollement restricts to the one of relative stable categories.

\subsection{Categories of Gorenstein-projective modules}\label{relative stable category}
In the subsection we describe the category of Gorenstein-projective modules over the endomorphism
algebra of a self-orthogonal generator.

Recall that the \emph{$M$-stable category} $A\Modcat/[M]$ of $A\Modcat$ is defined to be the quotient category of $A\Modcat$ modulo $\Add(M)$. From now on, we set $$\mathscr{D}:=A\Modcat/[M]\quad\mbox{and}\quad\StHom_M(X_1, X_2):=\Hom_\mathscr{D}(X_1, X_2)$$ for $X_1, X_2\in A\Modcat$. We say that $X_1$ and $X_2$ are \emph{$M$-stably isomorphic} if  they are isomorphic in $\mathscr{D}$. For a full subcategory $\mathscr{U}$ of $A\Modcat$, we denote by $\mathscr{U}/[M]$ the full subcategory of $\mathscr{D}$ consisting of all objects $X$ which are $M$-stably isomorphic to objects of $\mathscr{U}$.

As $M$ is finitely generated and $\Add(M)=\Prod(M)$, we know that $\Add(M)$ is a functorially finite subcategory in $A\Modcat$.
According to \cite[Chap. II. 1]{BI}, $\mathscr{D}$ is a \emph{pretriangulated category} mainly consisting of the following data:

$(1)$ An adjoint pair $(\Omega_M^{-}, \Omega_M)$ of additive endofunctors $\Omega_M^-, \Omega_M: \mathscr{D}\to \mathscr{D}$.

For an $A$-module $X$, $\Omega_M^-(X)$ is defined to be the cokernel of a minimal left $\Add(M)$-approximation $\ell_X: X\to M^X$ of $X$, while
$\Omega_M(X)$ is defined to be the kernel of a minimal right $\Add(M)$-approximation $r_X: M_X\to X$ of $X$. The existence of minimal approximations
follows from Remark \ref{DDE}. Moreover, $\Omega_M^-(X)$ and $\Omega_M(X)$ are unique up to isomorphism.

$(2)$ A collection of \emph{right triangles} (up to isomorphism) of the form
$
X\lraf{f}Y\lraf{g}Z\lraf{h} \Omega_M^-(X)
$
arising from an exact commutative diagram in $A\Modcat$:
$$\xymatrix{
0\ar[r]&X\ar[r]^-{f}\ar@{=}[d]&Y\ar[d]\ar[r]^-{g}& Z\ar[r]
\ar[d]^-{h}&0\\
0\ar[r]&X\ar[r]^-{\ell_X} & M^X\ar[r]  & \Omega_M^-(X)\ar[r]  &0
}
$$
where $\Hom_A(f, M): \Hom_A(Y, M)\to\Hom_A(X, M)$ is surjective.

$(3)$  A collection of \emph{left triangles} (up to isomorphism) of the form $\Omega_M(Z)\to X\to Y\to Z$ which is defined in a dual way as in $(2)$.

$(4)$ Right triangles (respectively, left triangles) satisfy all the axioms of a triangulated category, except that $\Omega_M^-$ (respectively,
$\Omega_M$) is not necessarily an equivalence.

\smallskip
In the following,  $\Omega_M$ is called the \emph{$M$-syzygy} functor on $\mathscr{D}$. For an $A$-module $X$, we put $\Omega_M^0(X):=X$, and $\Omega_M^n(X):=\Omega_M(\Omega_M^{n-1}(X))$ for $n\geq 1$. The functor $\Omega_M^n$ is called the \emph{$n$-th $M$-syzygy} functor on $\mathscr{D}$. Dually, $\Omega_M^{-}$ is called the \emph{$M$-cosyzygy} functor and the \emph{$n$-th $M$-cosyzygy} functor $\Omega_M^{-n}$ is defined dually. If $M$ is a self-orthogonal generator for $A\modcat$, then $(\Omega_M^{-n}, \Omega_M^n):\mathscr{D}\to\mathscr{D}$ is an adjoint pair for $n\geq 2$, see \cite[Lemmas 3.2 and 3.3]{cs}.

The following simple observation will be used in the later discussions.

\begin{Lem}\label{M-stable category}
$(1)$ Let $X$ and $Y$ be $A$-modules. Then $X$ and $Y$ are $M$-stably isomorphic if and only if there are $M_1, M_2\in\Add(M)$ such that $X\oplus M_1\simeq Y\oplus M_2$ in $A\Modcat$.

$(2)$ Let $X\stackrel{f}{\ra} Y\stackrel{g}{\ra} Z\ra X[1]$ be a triangle in $\Stmc{A}$. If $Z\in {^{\bot 1}}M$, then there is a right
triangle $X\stackrel{f}{\ra} Y\stackrel{g}{\ra} Z\ra \Omega_M^-(X)$ in $\mathscr{D}$.
\end{Lem}
{\it Proof.} $(1)$  This can be proved similarly as done in $\Stmc{A}$.

$(2)$ Up to isomorphism of triangles in $\Stmc{A}$, we can assume that the sequence $0\to X\stackrel{f}{\ra} Y\stackrel{g}{\ra}Z\to 0$ is exact. If $Z\in {^{\bot 1}}M$, then $\Hom_A(f, M)$ is surjective. This implies $(2)$.   $\square$

\medskip
By Lemma \ref{M-stable category}, if $\mathscr{U}$ is a full subcategory of $A\Modcat$ containing $\Add(M)$ and being closed under both direct summands and finite direct sums, then  $\mathscr{U}/[M]$ is closed under direct summands in  $\mathscr{D}$.

Let
$$\mathscr{G}:=M^{\bot >0}\cap {^{\bot >0}}M,\quad  \mathscr{C}:=\mathscr{G}/[M],\quad\mathscr{G}_0:=\mathscr{G}\cap A\modcat.$$
Then $\mathscr{G}$ is always closed under filtered colimits in $A\Modcat$, that is, colimits of filtered diagrams $I\to A\Modcat$ with $I$ an essentially small, filtered category.

In fact, since $M$ is finitely generated, the functor $\Ext_A^n(M,-)$ commutes with filtered colimits for each $n\geq 0$.
Then $M^{\bot >0}\subseteq A\Modcat$ is closed under filtered colimits. Thanks to $M\simeq D(DM)$ as $A$-modules,
$M$ is pure-injective. It follows that $\Ext_A^n(-, M)$ sends filtered colimits to filtered limits, and therefore
${^{\bot >0}}M\subseteq A\Modcat$ is also closed under filtered colimits. Thus $\mathscr{G}\subseteq A\Modcat$ is closed under filtered colimits.

Let $\lim\limits_{\longrightarrow}\mathscr{G}_0$ denote the full subcategory of $A\Modcat$ consisting of all filtered colimits of modules in $\mathscr{G}_0$. Then $\lim\limits_{\longrightarrow}\mathscr{G}_0\subseteq\mathscr{G}$.

\medskip
\textbf{In the rest of this subsection, we assume that ${_A}M$ is \emph {self-orthogonal}.}
The following result is a unbounded version of \cite[Lemma 3.5]{cs}.

\begin{Lem}\label{Gorenstein}
$(1)$ The category $\mathscr{G}$ (respectively, $\mathscr{G}_0$) is a Frobenius category with the shift functor given by $\Omega_M^{-}$. The full subcategory of projective-injective objects of $\mathscr{G}$ (respectively, $\mathscr{G}_0$) equals $\Add(M)$ (respectively, $\add(M)$). In particular, $\mathscr{G}/[M]$ is a triangulated category.

$(2)$ The functor $G: A\Modcat\to \Lambda\Modcat$ restricts to equivalences of Frobenius categories:
$$\mathscr{G}\lraf{\simeq}\GP{\Lambda}\quad\mbox{\rm and}\quad  \mathscr{G}_0\lraf{\simeq}\gp{\Lambda}.$$
In particular, there are equivalences of triangulated categories:
$$\mathscr{C}\lraf{\simeq}\Lambda\mbox{-}\underline{\rm GProj}\quad\mbox{\rm and}\quad
\mathscr{G}_0/[M]\lraf{\simeq}\Lambda\mbox{-}\underline{\rm Gproj}.$$
\end{Lem}

{\it Proof.} For $n\ge 1$, we define
$$
\mathscr{G}_n:=\{X\in A\modcat\mid \Ext_A^i(M, \Omega_M^{-j}(X))=0=\Ext_A^i(\Omega_M^j(X), M)\;\;\mbox{for any}\;\; j\geq 0\;\;\mbox{and}\;\; 1\leq i\leq n\}$$
(see also \cite[Definition 3.4]{cs}). Then $\mathscr{G}_n$ is a Frobenius category and there is a chain of full subcategories of $A\modcat$:
$\mathscr{G}_1\supseteq \mathscr{G}_2\supseteq\cdots\supseteq\mathscr{G}_n\supseteq
\mathscr{G}_{n+1}\supseteq\cdots.$
Moreover, $G$  restricts to an equivalence of Frobenius categories: $\mathscr{G}_n\lraf{\simeq}\gp{\Lambda}$ for $n\ge 1$.
It follows that the inclusions $\mathscr{G}_{n+1}\to \mathscr{G}_n$ are equivalences of additive categories. Since $\mathscr{G}_n$ are closed under isomorphisms in $A\modcat$, we have $\mathscr{G}_n=\mathscr{G}_{n+1}$ for all $n\geq 1$. As $M$ is self-orthogonal, $\mathscr{G}_0$ is closed under taking $\Omega_M$ and $\Omega^{-}_M$ in $A\modcat$. Consequently, $\mathscr{G}_0\subseteq \bigcap_{n\geq 1}^{\infty}\mathscr{G}_n\subseteq\mathscr{G}_0$. This implies $\mathscr{G}_0=\mathscr{G}_n$, and thus Lemma \ref{Gorenstein} holds for $\mathscr{G}_0$.

Note that $G$ commutes with direct sums and restricts to an equivalence $\Add(M)\lraf{\simeq}\Prj{\Lambda}$.
As done in the proof of \cite[Lemma 3.5]{cs}, we can show that Lemma \ref{Gorenstein} holds first for
$$
\mathscr{\bf{G}}_n:=\{X\in A\Modcat\mid \Ext_A^i(M, \Omega_M^{-j}(X))=0=\Ext_A^i(\Omega_M^j(X), M)\;\;\mbox{for }\;\; j\geq 0\;\;\mbox{and}\;\; 1\leq i\leq n\}$$ and then for $\mathscr{G}$. $\square$

\medskip
Consequently, we get a characterization of virtually Gorenstein algebras in terms of compact objects.

\begin{Prop}\label{Limits}
Let $A$ be a self-injective algebra and $M$ a self-orthogonal generator for $A\modcat$. Then $\End_A(M)$ is virtually Gorenstein if and only if $\mathscr{G}=\lim\limits_{\longrightarrow}\mathscr{G}_0$ if and only if each compact object of $\mathscr{C}$ is $M$-stably isomorphic to a finitely generated module.
\end{Prop}

{\it Proof.} By \cite[Theorem 5]{BK} and \cite[Theorem 8.2]{Bel}, $\Lambda:=\End_A(M)$ is virtually Gorenstein if and only if $\GP{\Lambda}=\lim\limits_{\longrightarrow}(\gp{\Lambda})$ if and only if each compact object of $\Lambda\mbox{-}\underline{\rm GProj}$ is isomorphic to an object of $\Lambda\mbox{-}\underline{\rm Gproj}$. Note that $G: A\Modcat\to\Lambda\Modcat$ is fully faithful and commutes with filtered colimits. Moreover, since ${_A}M$ is a generator, two $A$-modules $X$ and $Y$ are isomorphic if and only if $G(X)$ and $G(Y)$ are isomorphic. Now the equivalences in Proposition \ref{Limits} follow from Lemmas \ref{Gorenstein} and \ref{M-stable category}. $\square$

\medskip
A complex $\cpx{P}\in\K{\Prj{\Lambda}}$ is called \emph {totally acyclic} if both $\cpx{P}$ and $\cpx{\Hom}_\Lambda(\cpx{P}, \Lambda)$ are acyclic. Let $\mathscr{K}_{\rm tac}(\Prj{\Lambda})$ be the full subcategory of $\K{\Prj{\Lambda}}$ consisting of totally acyclic complexes. It is known that there is a triangle equivalence $\Lambda\mbox{-}\underline{\rm GProj}\lraf{\simeq}\mathscr{K}_{\rm tac}(\Prj{\Lambda})$ which sends a Gorenstein-projective $\Lambda$-module to its complete projective resolution. Composing this equivalence with the equivalence
$\mathscr{C}\lraf{\simeq}\Lambda\mbox{-}\underline{\rm GProj}$ in Lemma \ref{Gorenstein}(2), we obtain a triangle equivalence
$$
\Hom_A(M, \cpx{M}_{-}):\;\;\mathscr{C}\lraf{\simeq} \mathscr{K}_{\rm tac}(\Prj{\Lambda}),\;\;  X\mapsto \Hom_A(M, \cpx{M}_X)
$$
where $\cpx{M}_X:=(M_X^n)_{n\in\mathbb{Z}}\in\Kac{\Add(M)}$ is defined by concatenating an $\Add(M)$-resolution $\cdots \to M_X^{-2}\to M_X^{-1}\to M_X^{0}\to X\to 0$ of $X$ with an $\Add(M)$-coresolution $0\to X\to M_X^{1}\to M_X^{2}\to M_X^{3}\to\cdots$ of $X$ at the position $X$. The complex $\cpx{M}_X$ is called a \emph {complete $\Add(M)$-resolution} of $X$. For each $n\in\mathbb{Z}$, there is an additive functor
$$
\Psi_n:=H^n\big(\StHom_A(M, \cpx{M}_{-})\big):\;\;\mathscr{C}\lra \Gamma\Modcat,\;\;  X\mapsto H^n\big(\StHom_A(M, \cpx{M}_X)\big)
$$
which is \emph{homological} in the sense that applying $\Psi_n$ to every triangle $X_1\to X_2\to X_3\to \Omega_M^{-}(X_1)$ in $\mathscr{C}$ yields an exact sequence
$\cdots\to\Psi_n(X_1)\to \Psi_n(X_2)\to \Psi_n(X_3)\to \Psi_{n+1}(X_1)\to \cdots.$

\begin{Rem}\label{Decomposition}
By Lemma \ref{Gorenstein}(1) and Remark \ref{DDE}(2), we have the following observation.

For each $X\in\mathscr{G}$, there are isomorphisms $X\simeq \Omega_M^-\Omega_M(X)\oplus M_0\simeq \Omega_M\Omega_M^-(X)\oplus M_1$ in $A\Modcat$
with $M_0, M_1\in\Add(M)$ such that $\Omega_M^-\Omega_M(X)$ and $\Omega_M\Omega_M^-(X)$ have no direct summands in $\Add(M)$.
Thus, if $X\in\mathscr{G}$ has no direct summands in $\Add(M)$, then we can choose
the $n$-th differential $d_X^n: M_X^n\to M_X^{n+1}$ such that the induced map $M_X^n\to \Img(d_X^n)$ is a minimal right $\Add(M)$-approximation of $\Img(d_X^n)$.
\end{Rem}

\begin{Lem} \label{Cotorsion pair}
$(1)$ The functor $F\circ \nu_A: A\Modcat\to\Lambda\Modcat$ restricts to an equivalence:
$\mathscr{G}\lraf{\simeq} \GI{\Lambda}$ of Frobenious categories.

$(2)$  $(\mathscr{G}, \mathscr{G}^{\bot 1})$ and $({^{\bot 1}}\mathscr{G}, \mathscr{G})$ are cotorsion pairs in $A\Modcat$ such that $\mathscr{G}\cap\mathscr{G}^{\bot 1}=\Add(M)={^{\bot 1}}\mathscr{G}\cap\mathscr{G}$. In particular, $\mathscr{G}$ is functorially finite in $A\Modcat$.

$(3)$ The inclusion $\mathscr{C}\hookrightarrow \mathscr{D}$ admits a left adjoint $T: \mathscr{D}\to\mathscr{C}$ which is induced from minimal left $\mathscr{G}$-approximations of modules. Moreover, $T$ preserves compact objects, sends right triangles to triangles and
commutes with the functor $\Omega_M^-$.

$(4)$ The set $\{T(X)\mid X\in A\modcat\}$ is a compact generating set of $\mathscr{C}$.
\end{Lem}

{\it Proof.} $(1)$ Since the adjoint pair $(\nu_\Lambda, \nu_\Lambda^{-})$ induces quasi-inverse equivalences
$\GP{\Lambda}\lraf{\simeq}\GI{\Lambda}$, $(1)$ follows from Lemmas \ref{Gorenstein}(2) and \ref{Equivalence}(2) together with Remark \ref{DDE}(1).

$(2)$ Since $\Omega^-_M(\mathscr{G})\subseteq \mathscr{G}$ by Lemma \ref{Gorenstein}(1),  the sequence $0\to X\lraf{\ell_X} M^X\to \Omega^-_M(X)\to 0$ splits whenever $X\in\mathscr{G}\cap\mathscr{G}^{\bot 1}$. This implies
$\mathscr{G}\cap\mathscr{G}^{\bot 1}\subseteq \Add(M)$. Clearly, $\Add(M)\subseteq\mathscr{G}\cap\mathscr{G}^{\bot 1}$. Thus $\mathscr{G}\cap\mathscr{G}^{\bot 1}=\Add(M)$. Dually, ${^{\bot 1}}\mathscr{G}\cap\mathscr{G}=\Add(M)$.

Next, we show that $(\mathscr{G}, \mathscr{G}^{\bot 1})$ is a cotorsion pair in $A\Modcat$.
Since $A\Modcat$ is an abelian category with enough projectives and injectives, it suffices to show that, for any $A$-module $X$, there is an exact sequence $0\to X_2\to X_1\to X\to 0$ in $A\Modcat$ such that $X_1\in\mathscr{G}$ and $X_2\in\mathscr{G}^{\bot 1}$ (for example, see \cite[Lemma V. 3.3]{BI}).

Since each $\Lambda$-module admits a minimal right $\GP{\Lambda}$-approximation, we take a minimal right $\GP{\Lambda}$-approximation of $G(X)$, say $g: Y\to G(X)$. By Lemma \ref{Gorenstein}(2), we can assume $Y=G(X_1)$ for some $X_1\in\mathscr{G}$. As $G$ is fully faithful and $M$ is a generator, there is a surjective map
$f: X_1\to X$ of $A$-modules such that $g=G(f)$ and $f$ is a minimal right $\mathscr{G}$-approximation of $X$. Since $\mathscr{G}$ is closed under extensions in $A\Modcat$, it follows from Wakamatsu's Lemma that $\Ker(f)\in\mathscr{G}^{\bot 1}$. Hence the sequence $0\to \Ker(f)\to X_1\to X\to 0$ is a desired one.

Similarly, to show that  $({^{\bot 1}}\mathscr{G}, \mathscr{G})$ is a cotorsion pair in $A\Modcat$, it is enough to prove that there is an exact sequence $0\to X\to T^X\to C^X\to  0$ of $A$-modules such that $T^X\in\mathscr{G}$ and $C^X\in{^{\bot 1}}\mathscr{G}$.

Let $\overline{F}=F\circ \nu_A:  A\Modcat\to \Lambda\Modcat$. Then $\overline{F}$ is fully faithful and restricts to an equivalence $\mathscr{G}\lraf{\simeq}\GI{\Lambda}$ by $(1)$. Since each $\Lambda$-module admits a minimal left $\GI{\Lambda}$-approximation,
there is a map $h^X: X\to T^X$ in $A\Modcat$ with $T^X\in\mathscr{G}$ such that $\overline{F}(h^X)$ is a minimal left $\GI{\Lambda}$-approximation of $\overline{F}(X)$. This implies that $h^X$ is a minimal left $\mathscr{G}$-approximation of $X$. Moreover, $h^X$ is injective. This is due to $D(A)\in\mathscr{G}$ and $\Coker(h^X)\in{^{\bot 1}}\mathscr{G}$ by the dual of  Wakamatsu's Lemma. Thus the sequence
$0\to X\to T^X\to \Coker(h^X)\to 0$ is the one as desired.

$(3)$ Since $({^{\bot 1}}\mathscr{G}, \mathscr{G})$ is a cotorsion pair with ${^{\bot 1}}\mathscr{G}\cap\mathscr{G}=\Add(M)$ by (2), the functor $T: \mathscr{D}\to \mathscr{C}$ exists (see the comments after Definition \ref{def-cotor}).
As the inclusion functor $\mathscr{C}\to \mathscr{D}$ preserves
direct sums, we know that $T$ preserves compact objects. It is known that $\mathscr{C}$ is a triangulated category, $\mathscr{D}$ is a pretriangulated category and
$\mathscr{C}$ is a pretriangulated subcategory of $\mathscr{D}$. Thus the last assertion follows from the dual version of \cite[Proposition II. 2.6]{BI}.

$(4)$ We show that $A\modcat/[M]$ is a compact generating set of $\mathscr{D}$.
Clearly, if $X\in A\modcat$, then $X$ is compact in $A\Modcat$, and also compact in $\mathscr{D}$ by Lemma \ref{AJCP}(2). Let $Y\in A\Modcat$ such that $\StHom_M(X, Y)=0$ for all $X\in A\modcat$. Then each map from $X$ to $Y$ factorizes through an object of $\Add(M)$, and particularly, through $M_Y$ via the minimal right $\Add(M)$-approximation $r_Y: M_Y\to Y$ of $Y$.
Recall that, for an Artin algebra $B$, an exact sequence $0\to X_1\ra X_2\lraf{g} X_3\to 0$ of $B$-modules is called \emph{pure-exact} if $\Hom_B(Z, g)$ is surjective for any $Z\in B\modcat$; equivalently, $0\to L\otimes_B X_1\to L\otimes_BX_2\to L\otimes_BX_3\to 0$ is exact for any $L\in B\opp\Modcat$. This implies that  $0\to \Omega_M(Y)\to M_Y\to Y\to 0$ is pure-exact in $A\Modcat$. Note that there is a natural isomorphism $\Hom_\Lambda(U, G(-))\simeq \Hom_A(M\otimes_\Lambda U, -)$ for any $\Lambda$-module $U$ and that ${_A}M\otimes_\Lambda U\in A\modcat$ if $U\in\Lambda\modcat$. Hence $0\to G(\Omega_M(Y))\to G(M_Y)\to G(Y)\to 0$ is pure-exact in $\Lambda\Modcat$. Since $G(M_Y)$ is projective, there holds $\Tor^\Lambda_1(V, G(Y))=0$ for all $V\in\Lambda\opp\Modcat$. Thus $G(Y)$ is flat, and therefore projective since $\Lambda$ is an Artin algebra. It follows that $G(Y)$ is a direct summand of $G(M_Y)$. Then $Y\in\Add(M)$ due to $A\in\add({_A}M)$. Hence $Y=0$ in $\mathscr{D}$ and $A\modcat/[M]$ is a compact generating set of $\mathscr{D}$.

By $(3)$, if $X\in A\modcat$, then $T(X)$ is compact in $\mathscr{C}$. Moreover, since $T$ is a left adjoint of the inclusion $\mathscr{C}\to \mathscr{D}$,
one can check that $T$ always preserves generating sets. This shows $(4)$.  $\square$

\section{Restrictions of recollements to relative stable categories}\label{RRR}

In this section we prove Theorem \ref{Gorenstein-recollement}. As a preparation of the proof, we first show that the recollement in Corollary \ref{Self-orthogonal recollement} restricts to a recollement of $\mathscr{C}$.
Throughout this section we set up the following.

{\bf Assumption:} \emph{Let $A$ be a self-injective algebra and ${_A}M$ a self-orthogonal and Nakayama-stable generator for $A\modcat$.}

We set $\Gamma:=\StEnd_A(M)$.
By Lemma \ref{AR}, there exists a natural isomorphism of additive functors:
$$(\lozenge)\quad D\StHom_A(M,-)\simeq\StHom_A(-, \nu_A(M)[-1]):\;\; \Stmc{A}\lra \Gamma\opp\Modcat.$$

We define the following categories related to $M$:
$$
 \mathscr{Y}:=\{Y\in\Stmc{A}\mid \StHom_A(M, Y[n]=0\;\;\mbox{for any}\;\;n\leq  0\},
$$
$$
\mathscr{X}:=\{X\in\Stmc{A}\mid \StHom_A(X,Y)=0\;\;\mbox{for any}\;\; Y\in\mathscr{Y}\},
$$
$$\mathscr{H}:=\mathscr{X}\cap \mathscr{Y}[1].$$
Then $\mathscr{X}$ is the smallest full subcategory of $\Stmc{A}$ containing $M$ and being closed under $[1]$, extensions and direct sums, $(\mathscr{X}, \mathscr{Y})$ is a torsion pair in $\Stmc{A}$, and $\mathscr{H}$ is an abelian category and called the \emph{heart} of $(\mathscr{X}, \mathscr{Y})$ (see \cite{BBD}). Clearly, $\underline{M}^{\bot}\subseteq\mathscr{Y}$ and $\mathscr{X}\subseteq \Loc(\underline{M})$. In general, $\mathscr{X}$ has not to be a triangulated subcategory of $\Stmc{A}$ since it is not necessarily closed under $[-1]$.

\begin{Prop}\label{Self-orthogonal}
$\mathscr{X}=\underline{M}^{\bot >0}\cap \underline{\mathscr{S}}\,$, $\mathscr{Y}={^{\bot \geq -1}}\underline{M}$ and $\mathscr{H}=\underline{M}^{\bot\neq 0}\cap\underline{\mathscr{S}}$.
\end{Prop}
{\it Proof.}  Since ${_A}M$ is self-orthogonal, it follows from Lemma \ref{Endofunctor}(2) that $\underline{\mathscr{S}}=\Loc(\underline{M})$.  As $\underline{M}^{\bot >0}$
contains $M$ and is closed under $[1]$, extensions and direct sums in $\Stmc{A}$, we have $\mathscr{X}\subseteq \underline{M}^{\bot >0}$.
Thanks to $\mathscr{X}\subseteq\Loc(\underline{M})$,  we obtain $\mathscr{X}\subseteq\underline{M}^{\bot >0}\cap \Loc(\underline{M})$. To show the converse inclusion, we pick up $X\in \underline{M}^{\bot >0}\cap \Loc(\underline{M})$ and show $\StHom_A(X, Y)=0$ for any $Y\in\mathscr{Y}$.

Actually, it follows from $\underline{M}^{\bot >0}=\underline{M^{\bot >0}}$ that $X\in M^{\bot >0}$. Note that $X$ lies in $ \mathscr{S}$. Let $\cpx{\mu}_X:\cpx{M}_X\to X$ be an $\Add(M)$-resolution of $X$. It follows from Proposition \ref{S-cat}(1) that $\Hom_{\K{A}}(\Con(\cpx{\mu}_X), \Con(\cpx{\pi}_Y))=0$ for any $Y\in\mathscr{Y}$. Since $\StHom_A(M, \Omega_A^n(Y))=\StHom_A(M, Y[-n])=0$ for any $n\geq 0$, the chain map $\cpx{\pi}_Y: \cpx{P}_Y\to Y$ is an $\Add(M)$-resolution of  $Y$. Thus each homomorphism $f:X\to Y$ in $A\Modcat$ can be lifted to a chain map $\cpx{f}: \Con(\cpx{\mu}_X)\to \Con(\cpx{\pi}_Y)$ in $\K{A}$. From $\Hom_{\K{A}}(\Con(\cpx{\mu}_X), \Con(\cpx{\pi}_Y))=0$, we see that $f$ factorizes through the projective module $P_Y^0$. Thus $f=0$ in $\Stmc{A}$ and $\StHom_A(X, Y)=0$.

By definition, $\mathscr{Y}=\underline{M}^{\bot\leq 0}$. It follows from $\add(\nu_A(M))=\add({_A}M)$ and $(\lozenge)$ that $\mathscr{Y}={^{\bot \geq -1}}\underline{M}$. Since $\mathscr{H}=\mathscr{X}\cap\mathscr{Y}[1]$ and $\mathscr{Y}[1]=\underline{M}^{\bot<0}$,  one gets
$\mathscr{H}=\underline{M}^{\bot\neq 0}\cap\Loc(\underline{M})=\underline{M}^{\bot\neq 0}\cap\underline{\mathscr{S}}$. $\square$

\subsection{Restrictions of recollements\label{sect4.1}}
In this subsection we consider the restriction of the recollement in Corollary \ref{Self-orthogonal recollement}
to relative stable categories. This leads to a part of the proof of Theorem \ref{Gorenstein-recollement}.

We begin with the following preparation.

\begin{Lem}\label{Equality}
$(1)$ $\underline{M}^{\bot}={^{\bot}}\underline{M},\quad \underline{\mathscr{G}}=\underline{M}^{\bot\neq 0,-1}={^{\bot \neq 0, -1}}\underline{M}\quad \mbox{and}\quad \mathscr{H}\subseteq \underline{\mathscr{G}\cap\mathscr{S}}\subseteq\mathscr{X}.$

$(2)$ Let $\pi: \underline{\mathscr{G}}\to \mathscr{C}$ be the canonical quotient functor. Then the following are true.

\quad $(a)$ The composition $\underline{M}^{\bot}\hookrightarrow\underline{\mathscr{G}}\lraf{\pi} \mathscr{C}$ is a fully faithful triangle functor.
In particular, the image of this composition is a full triangulated subcategory of $\mathscr{C}$.

\quad $(b)$ The composition $\mathscr{H}\hookrightarrow\underline{\mathscr{G}\cap\mathscr{S}}\lraf{\pi} (\mathscr{G}\cap\mathscr{S})/[M]$ is a fully faithful additive functor.
\end{Lem}

{\it Proof.} (1) Since $\add(\nu_A(M))=\add({_A}M)$, $(1)$ follows from $(\lozenge)$ and Proposition \ref{Self-orthogonal}.

(2) If $X\in\underline{M}^{\bot}$, then $X\in{^{\bot}}\underline{M}$ by $(1)$. In this case, $\StHom_A(M, X)=0=\StHom_A(X, M)$, and therefore $\StHom_A(M_0, X)=0=\StHom_A(X, M_0)$ for $M_0\in\Add(M)$. This implies that $\Omega_A(X)\simeq \Omega_M(X)$, $\Omega_A^-(X)\simeq \Omega_M^-(X)$ in $\mathscr{D}$ and $\StHom_A(X, X')=\underline{\Hom}_M(X, X')$ for any $A$-module $X'$. Thus the composition in $(a)$ is fully faithful. It is a triangle functor since $\underline{M}^{\bot}$ is a full triangulated subcategory of $\Stmc{A}$. This shows $(a)$.

Now, let $X\in\mathscr{H}$. Then $\StHom_A(M, X[-1])=0$. Since $D\StHom_A(M, X[-1])\simeq \StHom_A(X, \nu_A(M))$ by $(\lozenge)$ and $\add(\nu_A(M))=\add(M)$, we have $\StHom_A(X, M)=0$. It then follows from $\Add(\underline{M})=\Prod(\underline{M})$ that $\StHom_A(X, M')=0$ for $M'\in\Add(M)$. Thus $\StHom_A(X, Y)=\StHom_M(X, Y)$ for any $A$-module $Y$. This implies $(b)$. $\square$

\begin{Prop}\label{Restriction}
The recollement in Corollary \ref{Self-orthogonal recollement} induces a recollement of triangulated categories:
$$
\xymatrix@C=1.2cm{\underline{M}^{\bot}\ar[r]^-{ \pi\circ{\rm inc}}
&\mathscr{C}\ar_-{\widetilde{\Psi'}}@/^1.2pc/[l]\ar[r]^-{\widetilde{\Phi}}\ar@/^1.2pc/[l]\ar_-{\widetilde{\Psi}}@/_1.2pc/[l]
&(\mathscr{G}\cap\mathscr{S})/[M]\ar_-{\Phi''}@/^1.2pc/[l]\ar@/^1.2pc/[l]\ar@/_1.2pc/[l]_{{\rm inc}}.}
$$
\end{Prop}

{\it Proof.} We first show that the recollement in Corollary \ref{Self-orthogonal recollement} can be restricted to a ``recollement'' of additive categories with six additive functors:
$$
(\sharp)\quad \xymatrix@C=1.2cm{\underline{M}^{\bot}\ar[r]^-{{\rm inc}}
&\underline{\mathscr{G}}\ar_-{\widetilde{\Psi'}}@/^1.2pc/[l]\ar[r]^-{\widetilde{\Phi}}\ar@/^1.2pc/[l]\ar_-{\widetilde{\Psi}}@/_1.2pc/[l]
&\underline{\mathscr{G}\cap\mathscr{S}}\ar_-{\Phi''}@/^1.2pc/[l]\ar@/^1.2pc/[l]\ar@/_1.2pc/[l]_{{\rm inc}}}
$$
which satisfy the conditions $(1)$-$(3)$ in Definition \ref{def01}. Obviously, ${\rm inc}: \underline{M}^{\bot}\to\underline{\mathscr{G}}$ has left and right adjoints which are the restriction of the functors $\widetilde{\Psi}$ and $\widetilde{\Psi'}$ in Corollary \ref{Self-orthogonal recollement} to $\underline{\mathscr{G}}$, respectively. Now, we claim that $\widetilde{\Phi}(\underline{\mathscr{G}})\subseteq\underline{\mathscr{G}\cap\mathscr{S}}$ and $\Phi''(\underline{\mathscr{G}\cap\mathscr{S}})\subseteq\underline{\mathscr{G}}$. Then $({\rm inc}, \widetilde{\Phi})$ and $(\widetilde{\Phi}, \Phi'')$ in $(\sharp)$ are adjoint pairs.

By Corollary \ref{localizing}, $\underline{\mathscr{S}}=\Loc(\underline{M})$. By Theorem \ref{Main result}, each $A$-module $X$ is endowed with a triangle $\Psi(X)[-1]\to\Phi(X)\to X\to \Psi(X)$ in $\Stmc{A}$ such that $\Phi(X)\in\underline{\mathscr{S}}$ and $\Psi(X)\in \underline{M}^{\bot}$. Note that $\underline{\mathscr{G}}$ contains $\underline{M}^{\bot}$ and is closed under extensions of triangles in $\Stmc{A}$, due to Lemma \ref{Equality}. Since $\Psi(X)$ and $\Psi(X)[-1]$ lie in $\underline{M}^{\bot}$, we see that $X\in\underline{\mathscr{G}}$ if and only if $\Phi(X)\in\underline{\mathscr{G}}$.
Thus $\widetilde{\Phi}(\underline{\mathscr{G}})=\Phi(\underline{\mathscr{G}})\subseteq\underline{\mathscr{G}\cap\mathscr{S}}$. Similarly, by the triangle
$\Psi'(X)\to X\to \Phi'(X)\to \Psi'(X)[1]$ in Proposition \ref{Exact sequence}, $X\in\underline{\mathscr{G}}$ if and only if $\Phi'(X)\in\underline{\mathscr{G}}$. This implies $\Phi''(\underline{\mathscr{G}\cap\mathscr{S}})=\Phi'(\underline{\mathscr{G}\cap\mathscr{S}})\subseteq\underline{\mathscr{G}}$.

Next, we show that the functors in $(\sharp)$ induce triangle functors among quotient categories.

By Corollary \ref{localizing}, $\mathscr{S}$ contains $\Add(M)$ and is closed under taking $\Omega_M$ and $\Omega_M^{-}$ in $A\Modcat$. Since $\Add(M)\subseteq\mathscr{G}$ and $\mathscr{C}$ is a triangulated category by Lemma \ref{Gorenstein}, $(\mathscr{G}\cap\mathscr{S})/[M]$ is a full triangulated subcategory of $\mathscr{C}$. Note that $\Phi(\Add(\underline{M}))=\Add(\underline{M})$ and $\Psi(\Add(\underline{M}))=0$ due to $\underline{M}\in\underline{\mathscr{S}}$. By Lemma \ref{AJCP}(1), the adjoint pairs $({\rm inc}, \widetilde{\Phi})$ and $(\widetilde{\Psi}, {\rm inc})$ in $(\sharp)$ induce adjoint pairs $({\rm inc}, \widetilde{\Phi}_0)$ and $(\widetilde{\Psi}_0, \pi\circ{\rm inc})$ of additive functors among triangulated categories:
$$
(\natural)\quad \xymatrix@C=1.2cm{\underline{M}^{\bot}\ar[r]^-{\pi\circ {\rm inc}}
&\mathscr{C}\ar[r]^-{\;\widetilde{\Phi}_0}\ar_-{\widetilde{\Psi}_0}@/_1.3pc/[l]
&(\mathscr{G}\cap\mathscr{S})/[M]\ar@/_1.3pc/[l]_{{\rm inc}}}.
$$
In this diagram, both ${\rm inc}$ and $\pi\circ {\rm inc}$ (see Lemma \ref{Equality}(2)) are fully faithful triangle functors.
It is known that any left or right adjoint of a triangle functor between triangulated categories is again a triangle functor.
Thus $\widetilde{\Phi}_0$ and $\widetilde{\Psi}_0$ are triangle functors.
In the following, we show that both $\pi\circ {\rm inc}$ and $\widetilde{\Phi}_0$ have right adjoints.

According to Lemma \ref{Endofunctor}(3), $\underline{\mathscr{T}}=\Coloc(\underline{\nu_A(M)})=\Img(\Phi')$. By the definition of $\mathscr{T}$,
if  $X\in\underline{\mathscr{T}}$, then $\Psi'(X)=0$ and $\Phi'(X)\simeq X $. This means that $$\Psi'(\Prod(\underline{\nu_A(M)}))=0\quad\mbox{and}\quad\Phi'(\Prod(\underline{\nu_A(M)}))=\Prod(\underline{\nu_A(M)}).$$
Since $\add(\nu_A(M))=\add(M)$ and ${_A}M\in A\modcat$, there holds $\Prod(\underline{\nu_A(M)})=\Prod(\underline{M})=\Add(\underline{M})$.
It follows that $\Psi'(\Add(\underline{M}))=0$ and $\Phi'(\Add(\underline{M}))=\Add(\underline{M})$. Since $\Add(\underline{M})\subseteq\underline{\mathscr{G}\cap\mathscr{S}}$ and $\Phi''=\Phi'\circ{\rm inc}$, we have $\Phi''(\Add(\underline{M}))=\Add(\underline{M})$.
Due to Lemma \ref{AJCP}(1), $\Phi'': \underline{\mathscr{G}\cap\mathscr{S}}\to\underline{\mathscr{G}}$ induces a functor  $\Phi''_0: (\mathscr{G}\cap\mathscr{S})/[M]\to\mathscr{C}$ which is a right adjoint of $\widetilde{\Phi}_0$, while
$\widetilde{\Psi'}:\underline{\mathscr{G}}\to \underline{M}^{\bot}$ induces a functor $\widetilde{\Psi'}_0: \mathscr{C}\to \underline{M}^{\bot}$ which is a right adjoint of $\pi\circ {\rm inc}$.

When acting on objects, $\widetilde{\Phi}_0$ and $\widetilde{\Phi}$ are the same. So we denote $\widetilde{\Phi}_0$ by $\widetilde{\Phi}$ for simplicity.
Similarly, we denote $\widetilde{\Psi}_0$, $\widetilde{\Psi'}_0$ and $\Phi''_0$ by $\widetilde{\Psi}$, $\widetilde{\Psi'}$ and $\Phi''$, respectively.

To show the existence of the recollement of $\mathscr{C}$ in Proposition \ref{Restriction},
it remains to show the existence of two canonical triangles in Definition \ref{def01}(4).
Let $X\in\mathscr{G}$. We have shown that there is a canonical triangle ${\rm inc}\circ\widetilde{\Phi}(X)\to X\to  {\rm inc}\circ\widetilde{\Psi}(X)\to
{\rm inc}\circ\widetilde{\Phi}(X)[1] $ in $\Stmc{A}$ such that $\widetilde{\Phi}(X)\in\underline{\mathscr{G}\cap\mathscr{S}}$ and $\widetilde{\Psi}(X)\in \underline{M}^{\bot}$. Since $\underline{M}^{\bot}={^{\bot}}\underline{M}\subseteq \underline{{^{\bot 1}}M}$ by Lemma \ref{Equality}(1), it follows from Lemma $\ref{M-stable category}(2)$ that this triangle induces a triangle ${\rm inc}\circ\widetilde{\Phi}(X)\to X \to {\rm inc}\circ\widetilde{\Psi}(X)\to\Omega_M^{-}({\rm inc}\circ\widetilde{\Phi}(X))$ in $\mathscr{C}$, which is the required second triangle in Definition \ref{def01}(4). Similarly, by the triangle $\Psi'(X)\to X\to \Phi'(X)\to \Psi'(X)[1]$ with $\Psi'(X)\in \underline{M}^{\bot}$ and $\Phi'(X)\in\underline{\mathscr{G}}$ and by the inclusion
$\underline{M}^{\bot}\subseteq \underline{M^{\bot 1}}$, we can obtain the first triangle in Definition \ref{def01}(4).
Thus the proof of Proposition \ref{Restriction} is completed.
$\square$

\medskip
Now, we define a full subcategory $\mathscr{E}$ of $\mathscr{G}$, which contains both $\mathscr{G}_0$ and $\Prj{A}$. $$\mathscr{E}:=\{X\in\mathscr{G}\mid \StHom_A(M, X),\, \StHom_A(M[1], X)\in \Gamma\modcat\}.$$

\begin{Lem}\label{E-property}
$(1)$ $\mathscr{E}=\{X\in\mathscr{G}\mid \StHom_A(M, X)\in\Gamma\modcat,\, \StHom_A(X, M)\in \Gamma\opp\modcat\}$.

$(2)$ $\mathscr{E}$ is a thick subcategory of $\mathscr{G}$ and $\mathscr{E}/[M]$ is a full triangulated subcategory of $\mathscr{C}$.
\end{Lem}

{\it Proof.} $(1)$ As $A$ is an Artin algebra over an Artin commutative ring $R$ and ${_A}M$ is finitely generated, $\Gamma$ is an Artin algebra
over $R$. Moreover, a $\Gamma$-module $N$ is finitely generated if and only if $_RN$ is finitely generated as an $R$-module. Thus $(1)$ follows
from  $D\StHom_A(M[1], X)\simeq \StHom_A(X, \nu_A(M))$ and $\add(\nu_A(M))=\add(M)$.

$(2)$ Clearly, $\mathscr{E}$ is closed under direct summands in $\mathscr{G}$. Now, let $0\to X_1\to X_2\to X_3\to 0$ be an exact sequence in $A\Modcat$ with $X_i\in\mathscr{G}$ for $1\leq i\leq 3$. Since $\underline{\mathscr{G}}=\underline{M}^{\bot\neq 0,-1}$ by Lemma \ref{Equality}, there is a long exact sequence in $R\modcat$:
{\small $$
0\to \StHom_A(M[1], X_1)\to \StHom_A(M[1], X_2)\to \StHom_A(M[1], X_3)\to \StHom_A(M, X_1)\to \StHom_A(M, X_2)\to \StHom_A(M, X_3)\to 0.
$$}
This implies that $\mathscr{E}$ has the two out of three property in $\mathscr{G}$. Thus $\mathscr{E}$ is a thick subcategory of $\mathscr{G}$.

By definition, $\mathscr{E}/[M]$ consists of all objects $X\in\mathscr{C}$ which is $M$-stably isomorphic to an object of $\mathscr{E}$.
To show that $\mathscr{E}/[M]$ is a full triangulated subcategory of $\mathscr{C}$, it suffices to show that $\mathscr{E}/[M]$ is closed under taking $\Omega_M$ and $\Omega_M^-$ in $\mathscr{C}$.

Let $X\in\mathscr{E}$. Due to $\StHom_A(M, X)\in\Gamma\modcat$, we see that $X$ has a right $\Add(M)$-approximation $f: M_0\oplus P_0\to X$ such that $M_0\in\add(M)$ and $P_0\in\Prj{A}$. Note that $\mathscr{E}$ contains $\add(M)$ and $\Prj{A}$. This forces $M_0\oplus P_0\in\mathscr{E}$, and further
$\Ker(f)\in\mathscr{E}$. It follows from $\Omega_M(X)\simeq \Ker(f)$ in $\mathscr{C}$ that $\Omega_M(X)\in \mathscr{E}/[M]$. Since $\StHom_A(X, M)$ is a $\Gamma\opp$-module, we can show  $\Omega_M^-(X)\in \mathscr{E}/[M]$.

The proof of $(2)$ also implies that if $X\in\mathscr{G}$, then $X\in \mathscr{E}$ if and only if it has a complete $\Add(M)$-resolution $\cpx{M}_X:=(M_X^i)_{i\in\mathbb{Z}}$ satisfying that $M_X^i=N^i\oplus P^i$ with $N^i\in\add(M)$ and $P^i\in\Prj{A}$. $\square$

\begin{Koro}\label{Compact-reco}
The recollement in Proposition \ref{Restriction} can be restricted a recollement of $\mathscr{E}/[M]$:
$$
\xymatrix@C=1.2cm{\underline{M}^{\bot}\ar[r]^-{\pi\circ {\rm inc}}
&\mathscr{E}/[M]\ar_-{\widetilde{\Psi'}}@/^1.2pc/[l]\ar[r]^-{\widetilde{\Phi}}\ar@/^1.2pc/[l]\ar_-{\widetilde{\Psi}}@/_1.2pc/[l]
&(\mathscr{E}\cap\mathscr{S})/[M].\ar_-{\Phi''}@/^1.2pc/[l]\ar@/^1.2pc/[l]\ar@/_1.2pc/[l]_{{\rm inc}}}
$$
\end{Koro}

\smallskip
{\it Proof.} Note that $\big((\mathscr{G}\cap\mathscr{S})/[M]\big)\cap (\mathscr{E}/[M])=(\mathscr{E}\cap\mathscr{S})/[M]$ and the image  of the functor $\pi\circ {\rm inc}: \underline{M}^{\bot}\to \mathscr{C}$ belongs to $\mathscr{E}/[M]$. Thus Corollary \ref{Compact-reco} is a consequence of Proposition \ref{Restriction}. $\square$

\subsection{Compacts objects and representability of homological functors}\label{Torsion pair}
In this subsection, we find out a special compact object  in $\mathscr{C}$ (or even in a bigger relative stable category) and establish a series of homological functors from $\mathscr{C}$ to $\Gamma\Modcat$ (see Theorem \ref{Compact}).

From now on, let ${\bf L}: \Stmc{A}\to \mathscr{Y}$ be the left adjoint of the inclusion $\mathscr{Y}\subseteq \Stmc{A}$ and let ${\bf R}: \Stmc{A}\to \mathscr{X}$ be the right adjoint of the inclusion  $\mathscr{X}\subseteq\Stmc{A}$.
Define
$$H^0_*:=[1]\circ {\bf L}\circ[-1]\circ{\bf R}:\;\; \Stmc{A}\lra\mathscr{H}, \;\; N:=H^0_*(M)\in\mathscr{H}\;\;\mbox{and} $$
$$
\mathscr{H}^{\rm fg}:=\{X\in\mathscr{H}\mid \StHom_A(M, X)\in\Gamma\modcat\}.$$
Then $H^0_*$ is a homological functor, that is, for a triangle $X_1\to X_2\to X_3\to X_1[1]$ in $\Stmc{A}$,
the sequence $H^0_*(X_1)\to H^0_*(X_2)\to H^0_*(X_3)$ is exact in $\mathscr{H}$.
Moreover, $\mathscr{H}^{\rm fg}=\mathscr{H}\cap\underline{\mathscr{E}}$ by Proposition \ref{Self-orthogonal} and Lemma \ref{Equality}(1). Note that $D\StHom_A(M, N[-1])\simeq \StHom_A(N, \nu_A(M))$ by $(\lozenge)$. Since $\add(\nu_A(M))=\add(M)$ and $N\in\mathscr{H}$, we get $\StHom_A(N, M)=0$.
This implies $\StEnd_A(N)=\End_\mathscr{C}(N)$.

\begin{Lem}\label{Basic propertices} Let $X\in\mathscr{X}$.

$(1)$ There is a canonical triangle
$[1]\circ {\bf R}\circ[-1](X)\to X\lraf{\tau_X} H^0_*(X)\to [2]\circ {\bf R}\circ[-1](X)$
in $\Stmc{A}$.

$(2)$ If $Y\in\mathscr{Y}[1]$, then $\StHom_A(\tau_X, Y): \StHom_A(H^0_*(X), Y)\to \StHom_A(X, Y)$ is an isomorphism.
In particular, there is a natural isomorphism
$$\StHom_A(\tau_M, -):\;\; \StHom_A(N, -)\lraf{\simeq} \StHom_A(M, -):\;\; \mathscr{H}\lra \Gamma\Modcat.$$
\end{Lem}

{\it Proof.}   Since $(\mathscr{X}, \mathscr{Y})$ is a torsion pair in $\Stmc{A}$, each $A$-module $W$ is endowed with a triangle ${\bf R}(W)\to W\to {\bf L}(W)\to {\bf R}(W)[1]$ in $\Stmc{A}$.  We take $W=X[-1]$ and get a triangle ${\bf R}(X[-1])\to X[-1]\to {\bf L}(X[-1])\to {\bf R}(X[-1])[1]$ in $\Stmc{A}$.  Shifting this triangle by $[1]$ yields another triangle $$(\dag)\quad {\bf R}(X[-1])[1]\to X\to {\bf L}(X[-1])[1]\to  {\bf R}(X[-1])[2].$$ Due to $X\in\mathscr{X}$, we have ${\bf R}(X)\simeq X$. Thus $H^0_*(X) = [1]\circ {\bf L}\circ[-1]\circ{\bf R}(X)\simeq {\bf L}(X[-1])[1].$ Then the triangle $(\dag)$ can be rewritten as
$[1]\circ {\bf R}\circ[-1](X)\to X\lraf{\tau_X} H^0_*(X)\to [2]\circ {\bf R}\circ[-1](X).$ This shows $(1)$.

Since $\mathscr{X}\subseteq \Stmc{A}$ is closed under $[1]$, it follows from ${\bf R}(X[-1])\in \mathscr{X}$ that  both $[1]\circ {\bf R}\circ[-1](X) $ and  $[2]\circ {\bf R}\circ[-1](X)$ belong to $\mathscr{X}$.  As $\underline{\Hom}_{A}(X,Y)=0$ for $X\in \mathscr{X}$ and $Y\in \mathscr{Y}$,  the first part of $(2)$ holds by $(1)$, while the second part of $(2)$ follows from  $M\in\mathscr{X}$ and $\mathscr{H}\subseteq \mathscr{Y}[1]$. $\square$
\medskip

The next result follows from \cite[Chap. III, Lemma 3.3 and Theorem 3.4]{BI}; see also \cite[Theorem 1.3(3)]{HKM} for the assertion $(3)$.

\begin{Lem}\label{Isomorphism}
$(1)$ For any $X\in\mathscr{X}$, the morphism $\tau_X$ in Lemma \ref{Basic propertices}(1) induces an isomorphism of $\Gamma$-modules:
$$\StHom_A(M, \tau_X): \StHom_A(M, X)\lraf{\simeq}\StHom_A(M, H^0_*(X)).$$

$(2)$ The functor $H^0_*$ induces an isomorphism $\Gamma\lraf{\simeq}\StEnd_A(N)$ of algebras such that
$\gamma\,\tau_M=\tau_MH^0_*(\gamma)$ for any $\gamma\in\Gamma$. In this sense, $\Gamma$ can be identified with $\StEnd_A(N)$.

$(3)$ $N$ is a small projective generator of $\mathscr{H}$ and the functor
$\StHom_A(N, -):\mathscr{H}\to \Gamma\Modcat$ is an equivalence of abelian categories.
\end{Lem}

An easy observation is the following result, of which $(2)$ conveys that $N$ has only finitely many indecomposable direct summands in $\mathscr{C}$.

\begin{Koro}\label{Gamma-mod}
$(1)$ There is a fully faithful functor $\Theta: \Gamma\modcat\to (\mathscr{E}\cap\mathscr{S})/[M]$ which sends $\Gamma$ to $N$.

$(2)$ Let $M=A\oplus \bigoplus_{i=1}^{m}M_i$, where $m\in\mathbb{N}$ and $M_i$ are indecomposable and non-projective. Then $N\simeq \bigoplus_{i=1}^{m}H_*^0(M_i)$
in $\mathscr{C}$ and $H_*^0(M_i)$ are indecomposable in $\mathscr{C}$.

$(3)$ Let $X\in \mathscr{G}\cap\mathscr{S}$. Then there is a triangle
$X_0\to X\to H^0_*(X)\to \Omega_M^{-}(X_0)$ in $\mathscr{C}$ such that $X_0\in \mathscr{G}\cap\mathscr{S}$ and $\Omega_M(X_0)\in\mathscr{H}$.
Further, if $X\in \mathscr{E}\cap\mathscr{S}$, then $\Omega_M(X_0), H^0_*(X)\in\mathscr{H}^{\rm fg}$.\end{Koro}

{\it Proof.} $(1)$ By Lemma \ref{Equality}(2)(b), the composition $U: \mathscr{H}^{\rm fg}=\mathscr{H}\cap\underline{\mathscr{E}}\hookrightarrow\underline{\mathscr{E}\cap\mathscr{S}}\to (\mathscr{E}\cap\mathscr{S})/[M]$ is fully faithful. Moreover, from Lemmas \ref{Basic propertices}(2) and \ref{Isomorphism}(3) it follows that $\StHom_A(M, -)\simeq \StHom_A(N, -):\mathscr{H}^{\rm fg}\to\Gamma\modcat$, which is an equivalence of abelian categories sending $N$ to $\Gamma$. In particular, $N\in\mathscr{H}^{\rm fg}$. Now, let $\Theta:\Gamma\modcat\to (\mathscr{E}\cap\mathscr{S})/[M]$ be the composition of a quasi-inverse of $\StHom_A(M, -)$ with $U$. Then $\Theta$ is fully faithful and sends $\Gamma$ to $N$.

$(2)$ By Lemmas \ref{Basic propertices}(2) and \ref{Isomorphism}(1), $H_*^0$ induces an isomorphism $\StHom_A(M, M_i)\simeq \StHom_A(N, H_*^0(M_i))$.
Thus $\StEnd_A(M_i)\simeq \StEnd_A(H_*^0(M_i))$ as algebras. It then follows from $H_*^0(M_i)\in\mathscr{H}$ that $\StEnd_A(H_*^0(M_i))\simeq\End_\mathscr{C}(H_*^0(M_i))$ as algebras by Lemma \ref{Equality}(2).
As $M_i$ is indecomposable and non-projective, it is also indecomposable in $\Stmc{A}$. Thus $H_*^0(M_i)$ is indecomposable in $\mathscr{C}$.

$(3)$ Let $X_0:={\bf R}(X[-1])[1]$. Then $X_0[-1]\in\mathscr{X}$ and there is a triangle $X_0\to X\lraf{\tau_X} H^0_*(X)\to X_0[1]$ in $\Stmc{A}$ by Lemma \ref{Basic propertices}(1). Thus $\StHom_A(M, X_0[n])=0$ for all $n\geq 0$ by Proposition \ref{Self-orthogonal}. In particular, $\Ext_A^1(M, X_0)=0$.
Since both $X$ and $H^0_*(X)$ lie in $\mathscr{G}\cap\mathscr{S}$, we see that $X_0\in\mathscr{G}\cap\mathscr{S}$ and the above triangle induces a triangle
$X_0\to X\to H^0_*(X)\to \Omega_M^{-}(X_0)$ in $\mathscr{C}$ by Lemma \ref{M-stable category}(2). Further, $\StHom_A(M, X_0[-2])=0$ by Lemma \ref{Equality}(1),
and $X_0[-1]=\Omega_A(X_0)\simeq\Omega_M(X_0)\in\mathscr{G}$ by $\StHom_A(M, X_0)=0$. Thus $X_0[-1]\in\mathscr{H}$ by Proposition \ref{Self-orthogonal}.
Suppose that $X$ is in $\mathscr{E}\cap\mathscr{S}$. Then $\StHom_A(M[1], X_0)\simeq \StHom_A(M, X_0[-1])\simeq\StHom_A(M, {\bf R}(X[-1]))\simeq\StHom_A(M, X[-1])\in \Gamma\modcat$, where the last isomorphism follows from the fact that ${\bf R}: \Stmc{A}\to \mathscr{X}$ is the right adjoint of the inclusion
$\mathscr{X}\subseteq\Stmc{A}$. Now, it follows from $X_0\in\mathscr{G}$ and $\StHom_A(M, X_0)=0$ that $X_0\in\mathscr{E}$. Thus $X_0\in\mathscr{E}\cap\mathscr{S}$.
By Lemma \ref{E-property}(2), both $\Omega_M(X_0)$ and $H^0_*(X)$ belong to $\mathscr{H}^{\rm fg}$. $\square$

\begin{Koro}\label{Vanishing}
The module ${_A}M$ is projective if and only if $(\mathscr{E}\cap\mathscr{S})/[M]=0$ if and only if  $(\mathscr{G}\cap\mathscr{S})/[M]=0$.
\end{Koro}

{\it Proof.} If ${_A}M$ is projective, then $\mathscr{S}$ consists of projective $A$-modules by Corollary \ref{localizing}, which implies $(\mathscr{E}\cap\mathscr{S})/[M]=0=(\mathscr{G}\cap\mathscr{S})/[M]$.
If $(\mathscr{E}\cap\mathscr{S})/[M]=0$, then $\Gamma\modcat=0$ by Corollary \ref{Gamma-mod}(1), that is, $\Gamma=0$, and therefore ${_A}M$ is projective. $\square$

\begin{Lem}\label{CR}
The object $N$ has a complete $\Add(M)$-resolution:
$$ \cdots\lra M^{-3}\lra  M^{-2}\lra P^{-1}\lraf{\partial} M\oplus P^0\lra P^1\lra M^{2}\lra M^{3}\lra \cdots
$$
satisfying

 $(1)$ $N=\Coker(\partial)$, $M^i\in\Add(M)$ for all $|i|\geq 2$ and $P^i\in\Add(A)$ for all $|i|\leq 1$; and

 $(2)$ $\Ker(\partial)\in\mathscr{H}$ and $\StHom_A(M,\Ker(\partial))\simeq D\big(\StHom_A(M, \nu_A(M))\big)$ as $\Gamma$-modules.

 \end{Lem}
{\it Proof.} $(1)$ Let $\tau_M: M\to N$ be the morphism in $A\Stmc$ in Lemma \ref{Basic propertices}(1). Up to projective direct summand, we may assume that $\tau_M$ is a homomorphism in $A\Modcat$ and a preimage of the $\tau_M$ in $A\Stmc$.  Let $f:P^0\to\Coker(\tau_M) $ be a projective cover of $\Coker(\tau_M)$. Then there exists a homomorphism $f_0: P^0\to N$ such that $f$ is the composition of $f_0$ with the quotient map $N\to \Coker(\tau_M)$. Moreover, the map $\pi_0:=(\tau_M, f_0): M\oplus P^0\to  N$ is a right $\Add(M)$-approximation of $N$. Let $K:=\Ker(\pi_0)$ with a projective cover $\pi_1: P^{-1}\to K$.
Then $K\simeq \Omega_M(N)$ in $\mathscr{D}$. By the proof of Corollary \ref{Gamma-mod}(3), we see that $K\in\mathscr{G}\cap\mathscr{S}$, $\StHom_A(M, K)=0$ and
$K[-1]=\Omega_A(K)\simeq\Omega_M(K)\in\mathscr{H}$. Consequently, $\pi_1$ is a right $\Add(M)$-approximation of $K$. Now, let $\partial: P^{-1}\to M\oplus P^0$ be the composition of $\pi_1$ with the inclusion $\lambda_1: K\to  M\oplus P^0$. Then $\Coker(\partial)=N$ and $\Ker(\partial)\simeq K[-1]$ in $\Stmc{A}$. Since $\StHom_A(N, M)=0$ and $\Add(M)=\Prod(M)$, a minimal injective envelope $N\to P^1$ is a left $\Add(M)$-approximation of $N$. Thus $(1)$ holds.

$(2)$ Thanks to $N\in\mathscr{H}$, there holds $\StHom_A(M, N[-1])=0=\StHom_A(M, N[-2])$ by Proposition \ref{Self-orthogonal}. This implies that
$\StHom_A(M, \lambda_1[-1]): \StHom_A(M, K[-1])\to \StHom_A(M, M[-1])$ is an isomorphism. By $(\lozenge)$,
$$\StHom_A(M, M[-1])\simeq D\big(\StHom_A(M[-1], \nu_A(M)[-1])\big)\simeq D\big(\StHom_A(M, \nu_A(M))\big)$$ as $\Gamma$-modules.
Thus $\StHom_A(M, \Ker(\partial)\simeq \StHom_A(M, K[-1])\simeq D(\StHom_A(M, \nu_A(M)))$. $\square$

\begin{Koro}\label{2-periodic}
If $A$ is a symmetric algebra and $\Gamma$ is a Frobenius algebra, then $\Omega_M^2(N)\simeq N$ in $\mathscr{C}$.
\end{Koro}

{\it Proof.}
Since $A$ is symmetric, $\nu_A(M)\simeq M$ as $A$-$\Lambda$-bimodules. By Lemma \ref{CR}(2), $\StHom_A(M, \Ker(\partial))\simeq D\big(\StHom_A(M, M)\big)=D(\Gamma)$ as $\Gamma$-modules. Note that $\StHom_A(M, N)\simeq \Gamma$ by Lemmas \ref{Basic propertices}(2) and \ref{Isomorphism}(2). Since $\Gamma$ is a Frobenius algebra, $\Gamma\simeq D(\Gamma)$ as $\Gamma$-modules. Moreover, $\Ker(\partial)\in \mathscr{H}$ and $\Ker(\partial)\simeq\Omega_M^2(N)$ in $\mathscr{C}$
by Lemma \ref{CR}. It follows from Lemmas \ref{Basic propertices}(2) and  \ref{Isomorphism} that $\Ker(\partial)\simeq N$ in $\mathscr{H}$. Thus $\Omega_M^2(N)\simeq  N$ in $\mathscr{C}$. $\square$

\begin{Lem}\label{Left exact}
Let $X\in {^{\bot >0}}M$ and $l_X: X\to M^X$ be a left $\Add(M)$-approximation of $X$. Then
there is an exact sequence of $\Gamma$-modules:
$$0\lra \StHom_A(N, X)\lraf{(\tau_M)_*} \StHom_A(M,X)\lraf{(l_X)^*} \StHom_A(M, M^X).
$$
\end{Lem}

{\it Proof.} Keep all notations introduced in the proof of Lemma \ref{CR}.  Recall that $\Prod(M)=\Add(M)$ and $\StHom_A(N, M)=0$. If  $M_0\in\Add(M)$, then $\StHom_A(N, M_0)=0$. It follows from $M^X\in\Add(M)$ that $\StHom_A(N, M^X)=0$. This implies that the composition of $(\tau_M)_*$ with $(l_X)^*$ is $0$. Thus
$\Img((\tau_M)_*)\subseteq \Ker((l_X)^*)$.

Applying $\StHom_A(-, X)$ to the triangle $K\to M\lraf{\tau_M} N\to K[1]$ in $\Stmc{A}$ yields an exact sequence
$\StHom_A(K[1], X)\to \StHom_A(N, X)\lraf{(\tau_M)_*} \StHom_A(M,X)$. By $(\lozenge)$, we have $\underline{{^{\bot >0}}M}=\underline{M}^{\bot\leq -2}=\mathscr{Y}[2]$. Moreover, $\underline{K}[-1]\in\mathscr{H}\subseteq \mathscr{X}$ by the proof of Lemma \ref{CR}. Thus $\underline{K}[1]\in\mathscr{X}[2]$.
Since $X\in\underline{{^{\bot >0}}M}$ and $(\mathscr{X}, \mathscr{Y})$ is a torsion pair in $\Stmc{A}$, we have $\StHom_A(K[1], X)=0$.
Thus $(\tau_M)_*$ is injective. It remains to show $\Ker((l_X)^*)\subseteq\Img((\tau_M)_*)$.

Suppose that $g\in\StHom_A(M, X)=\StHom_A(M\oplus P_0, X)$ with  $g l_X=0$. Let
$0\to X\lraf{\lambda_0} Q^0\lraf{\mu_0}\Omega_A^-(X)\to 0$ be a short exact
sequence  of $A$-modules in which $\lambda_0$ is an injective envelope. We claim that the pair $(g\lambda_0, 0)$ with $g\lambda_0: M\oplus P^0\to Q^0$ and $0: N\to \Omega_A^{-}(X)$ can be completed into the  commutative diagram
$$
\xymatrix{
\cdots\ar[r]& M^{-3}\ar[r]\ar@{-->}[d]^{f^{-3}}&  M^{-2}\ar[r]\ar@{-->}[d]^{f^{-2}}& P^{-1}\ar[r]^-{\partial}\ar@{-->}[d]^-{f^{-1}}& M\oplus P^0\ar[r]^-{\pi_0}
\ar[d]^-{g\lambda_0}& N\ar[r]\ar[d]^-{0}& 0\\
\cdots\ar[r]& Q^{-3}\ar[r]&  Q^{-2}\ar[r]& Q^{-1}\ar[r]&  Q^0\ar[r]^-{\mu_0}& \Omega_A^{-}(X)\ar[r] & 0,
}
$$
in which the first arrow is an $\Add(M)$-resolution of $N$ (see Lemma \ref{CR}) and the second arrow is a minimal projective resolution of $\Omega_A^{-}(X)$. In other words, there is a chain map $\cpx{f}:=(\cdots, f^{-3}, f^{-2}, f^{-1}, g\lambda_0, 0).$
To show the existence of the chain map, one splits the long exact sequences into a series of short exact sequences, and then construct relevant homomorphisms together with commutative diagrams between these sequences. We carry out the details as follows.

For $K :=\Ker(\pi_0)$, let $\lambda_1: K\to M\oplus P^0$ be the canonical inclusion, and let $\alpha: K\to X$ be the composition of $\lambda_1$ with $g: M\oplus P^0\to X$. Since $P^{-1}$ is projective, there exists a map $f^{-1}: P^{-1}\to Q^{-1}$ making the diagram commute:
$$
\xymatrix{
0\ar[r]&\Omega_A(K)\ar[r]\ar@{-->}[d]^{\Omega_A(\alpha)}& P^{-1}\ar[r]\ar@{-->}[d]^{f^{-1}}& K\ar[r]\ar@{-->}[d]^-\alpha&0 \\
0\ar[r]& \Omega_A(X)\ar[r]&  Q^{-1}\ar[r]^-{\mu_1}& X\ar[r]& 0.
}
$$
To construct $f^{-2}$, we will show that $\StHom_A(M, \alpha[-1]): \StHom_A(M, K[-1])\to \StHom_A(M, X[-1])$
is zero, where $[-1]$ denotes $\Omega_A$ in $\Stmc{A}$.

By the proof of Lemma \ref{CR}(2),  $\StHom_A(M, \lambda_1[-1])$ is an isomorphism. It is enough to show that the map $\StHom_A(M, g[-1]): \StHom_A(M, M[-1])\to \StHom_A(M, X[-1])$ is zero. Applying $D:\Gamma\Modcat\to \Gamma\opp\Modcat$ to this map, we see that $(\lozenge)$
yields  $D\StHom_A(M, g[-1])\simeq \StHom_A(g, \nu_A(M)): \StHom_A(X, \nu_A(M))\to \StHom_A(M, \nu_A(M))$.
Although $D$ may not be an equivalence in general, it is always exact and reflects zero objects. Thus $\StHom_A(M, g[-1])=0$ if and only if $\StHom_A(g, \nu_A(M))=0$. Since $\add(\nu_A(M))=\add(M)=\Prod(M)$, $\StHom_A(g, \nu_A(M))=0$ is equivalent to saying that $\StHom_A(g, M'): \StHom_A(X, M')\to \StHom_A(M, M')$ is $0$ for any $M'\in\Prod(M)$. By assumption, $l_X$ is a left $\Add(M)$-approximation of $X$ and
$gl_X=0$ in $\Stmc{A}$. This leads to $\StHom_A(g, M')=0$. Thus $\StHom_A(M, \alpha[-1])=0$, and therefore
$\StHom_A(M^{-2}, \alpha[-1])=0$, due to $M^{-2}\in\Add(M)$. Hence there are two homomorphisms $f^{-2}$ and $\beta$ such that the diagram is commutative:
$$
\xymatrix{
0\ar[r]&\Omega_M(\Omega_A(K))\ar[r]\ar@{-->}[d]^{\beta}& M^{-2}\ar[r]\ar@{-->}[d]^{f^{-2}}& \Omega_A(K)\ar[r]\ar@{-->}[d]^-{\Omega_A(\alpha)}&0 \\
0\ar[r]& \Omega_A^2(X)\ar[r]&  Q^{-2}\ar[r]& \Omega_A(X)\ar[r]& 0.
}
$$

Since $X\in \underline{{^{\bot >0}}M}=\underline{M}^{\bot\leq -2}$, we have  $\StHom_A(M, \Omega_A^n(X))=\StHom_A(M, X[-n])=0$ for all $n\geq 2$. Then $\StHom_A(M^{-n-1}, \Omega_A^n(X))=0$, due to $M^{-n-1}\in\Add(M)$. Consequently, the components $f^{-n-1}: M^{-n-1}\to Q^{-n-1}$ for $n\geq 2$ in $\cpx{f}$ can be constructed.

Since $N\in\mathscr{H}\subseteq\underline{M}^{\bot >0}\cap\underline{\mathscr{S}}$ by Proposition \ref{Self-orthogonal}, it follows from Proposition \ref{S-cat}(1) that $\cpx{f}=0$ in $\K{A}$. Therefore there are homomorphisms $h: N\to Q^0$ and $s^0: M\oplus P^0\to Q^{-1}$ such that $g\lambda_0=s^0\mu_1\lambda_0+\pi_0h$ and $h\mu_0=0$. Since $\lambda_0$ is the kernel of $\mu_0$, there is a unique map $h_0:N\to X$ satisfying $h=h_0\lambda_0$. Since $\lambda_0$ is injective, $g=s^0\mu_1+\pi_0h_0$. This forces $g=\tau_Mh_0$ in $\Stmc{A}$ and shows $g\in\Img((\tau_M)_*)$.
Thus $\Ker((l_X)^*)\subseteq\Img((\tau_M)_*)$. Finally, via the isomorphism $\Gamma\simeq\StEnd_A(N)$ in Lemma \ref{Isomorphism}(2), $\StHom_A(N, X)$ becomes a $\Gamma$-module and $(\tau_M)_*$  is a homomorphism of $\Gamma$-modules.  $\square$

\begin{Koro}\label{Direct sums}
The object $N$ is compact in $\underline{{^{\bot >0}}M}$.
\end{Koro}
{\it Proof.}
It suffices to show that $\StHom_A(N, -): \underline{{^{\bot >0}}M}\to \Gamma\Modcat$ commutes with direct sums. Let $\{X_i\}_{i\in I}$ be a set of $A$-modules in ${^{\bot >0}}M$ with $I$ an index set. For each $i$, let $g_i: X_i\to M_i$ be a left $\Add(M)$-approximation of $X_i$. According to $\Add(M)=\Prod(M)$, the direct sum $g=(g_i)_{i\in I}: \bigoplus_{i\in I}X_i\to\bigoplus_{i\in I}M_i$ of all these $g_i$ is a left $\Add(M)$-approximation of $\bigoplus_{i\in I}X_i$. Note that $\StHom_A(M, -): \Stmc{A}\to\Gamma\Modcat$ commutes with direct sums since $M$ is compact in $\Stmc{A}$. By Lemma \ref{Left exact}, we can construct the following commutative diagram with exact arrows and canonical vertical maps:
$$
\xymatrix{
0\ar[r]&\bigoplus_{i\in I}\StHom_A(N, X_i)\ar[r]\ar[d]& \bigoplus_{i\in I}\StHom_A(M,X_i)\ar[r]\ar[d]^{\simeq}
& \bigoplus_{i\in I}\StHom_A(M, M_i)\ar[d]^-{\simeq} \\
0\ar[r]&\StHom_A(N, \bigoplus_{i\in I}X_i)\ar[r]&\StHom_A(M, \bigoplus_{i\in I}X_i)\ar[r]&\StHom_A(M, \bigoplus_{i\in I}M_i).
}
$$
Thus the first vertical map is an isomorphism. $\square$

\begin{Theo}\label{Compact}
$(1)$ The object $N$ belongs to $(\mathscr{E}\cap\mathscr{S})/[M]$ and is compact in $({^{\bot >0}}M)/[M]$. In particular, $N$ is compact in $\mathscr{C}$.

$(2)$ For each $n\in\mathbb{Z}$, there exists a natural isomorphism of homological functors:
$$
H^n\big(\StHom_A(M, \cpx{M}_{-})\big)\lraf{\simeq}\StHom_M\big(\Omega_M^n(N),-\big):\;\;\mathscr{C}\lra \Gamma\Modcat.
$$
\end{Theo}

{\it Proof}
$(1)$ By Lemma \ref{Equality}, $N\in(\mathscr{G}\cap\mathscr{S})/[M]$. Note that $\StHom_A(M[1], N)=0$ by Proposition \ref{Self-orthogonal}
and $\StHom_A(M, N)\simeq \StHom_A(M, M)$ by Lemma \ref{Isomorphism}(1). This shows $N\in\mathscr{E}$. Since ${^{\bot >0}}M$ contains $\Add(M)$ and is closed under direct sums in $A\Modcat$, we see that ${^{\bot >0}}M$, as an additive category, has coproducts. Moreover, by Lemma \ref{AJCP}(2), the quotient functor $\underline{{^{\bot >0}}M}\to({^{\bot >0}}M)/[M]$ preserves coproducts and compact objects. Now $(1)$ follows from Corollary \ref{Direct sums}.

$(2)$
Let $\cpx{M}_X:=(M_X^n, d_X^n)_{n\in\mathbb{Z}}$ be a complete $\Add(M)$-resolution of $X$. Then $X\simeq\Coker(d_X^{-1})$ and
$\Hom_A(M, \cpx{M}_X)$ is acyclic. Clearly, there is a canonical ring homomorphism $\Lambda\twoheadrightarrow\Gamma$ and
the functor $\StHom_A(M,-):A\Modcat\to \Gamma\Modcat$ is naturally isomorphic to the composition of $G$ with
$\Gamma\otimes_\Lambda-$. Then $\StHom_A(M, \cpx{M}_X)\simeq\Gamma\otimes_\Lambda\Hom_A(M, \cpx{M}_X)$ as complexes.
Since $\Gamma\otimes_\Lambda-$ is right exact, the sequence $$\StHom_A(M, M^{n-1}_X)\lraf{(d_X^{n-1})^*}\StHom_A(M, M^n_X)\lra \StHom_A(M, \Omega_M^{-n}(X))\lra 0$$ is exact for all $n$.  As the inclusion $\lambda_n: \Omega_M^{-n}(X)\to M^{n+1}_X$ is a left $\Add(M)$-approximation of $\Omega_M^{-n}(X)$,
it follows from Lemma \ref{Left exact} that the sequence
$$
0\lra \StHom_A(N, \Omega_M^{-n}(X))\lraf{(\tau_M)_*} \StHom_A(M,\Omega_M^{-n}(X))\lraf{(\lambda_n)^*} \StHom_A(M, M_X^{n+1})
$$
is exact. Consequently,  $H^n(\StHom_A(M, \cpx{M}_X))\simeq H^n(\Gamma\otimes_\Lambda\Hom_A(M, \cpx{M}_X))\simeq\StHom_A(N, \Omega_M^{-n}(X)).$
Since $\StHom_A(N, M')=0$ for any $M'\in\Add(M)$, we have
$$\StHom_A(N, \Omega_M^{-n}(X))=\StHom_M(N, \Omega_M^{-n}(X))\simeq \StHom_M(\Omega_M^n(N), X).$$
This shows $(2)$.

\subsection{Compact objects from left approximations}\label{Compact-approximation}

In this subsection, we  characterize compact objects in $\mathscr{C}$ in terms of $M$-filtered modules (see Proposition \ref{CP}) and show that all objects of $(\mathscr{E}\cap\mathscr{S})/[M]$ are compact in $\mathscr{C}$. We then establish a connection between $(\mathscr{G}\cap\mathscr{S})/[M]$ and $(\mathscr{E}\cap\mathscr{S})/[M]$ by employing strong generators (see Corollary \ref{Compact recollement}).

 As a preparation, we recall a construction of ${\bf L}(X)$ and ${\bf R}(X)$ from the proof of \cite[Theorem III.2.3]{BI}.

Let $\mathfrak{A}:=\{\underline{M}[n]\mid n\geq 0\}\subseteq\mathscr{X}$. Denote by $\Add(\mathfrak{A})$ the full subcategory of $\Stmc{A}$ consisting of direct summands of arbitrary direct sums of objects of $\mathfrak{A}$. For a full subcategory $\mathscr{U}$ of $\Stmc{A}$, we denote by $\mathscr{U}^{\star n}=\mathscr{U}\star\mathscr{U}\star\cdots\star\mathscr{U}$ ($n$-factors) the category of $n$-extensions of $\mathscr{U}$ by $\mathscr{U}$ in $\Stmc{A}$.

Let $X\in A\Modcat$. We construct a right $\Add(\mathfrak{A})$-approximation $f_1: Q_1\to X$ of $X$ as follows:
Consider the set $I_X$ of the union of $\underline{\Hom}_{A}(P,X)$ with $P$ running over $\mathfrak{A}$,
define $Q_1=\bigoplus_{\lambda\in I_X}P_\lambda$ and take $f_1$ to be the morphism induced by $I_X$, where $\lambda$ is morphism from $P_{\lambda}$ to $X$. Then we extend $f_1$ to a triangle $Q_1\lraf{f_1}X\lraf{g_1} X_1\lra Q_1[1]$ in $\Stmc{A}$.
Now, we can repeat this construction by replacing $X$ by $X_1\in A\Modcat$. In general, for each $n\geq 0$, we can inductively construct a triangle $$Q_{n+1}\lraf{f_{n+1}} X_n\lraf{g_{n+1}}X_{n+1}\lra Q_{n+1}[1]$$ in $\Stmc{A}$ such that $f_{n+1}$ is a right $\Add(\mathfrak{A})$-approximation  of $X_n$ with $X_0:=X$. Setting $T_1:=Q_1$ and $h_1:=f_1$, we then construct inductively a tower of objects
$T_1\lraf{\tau_1} T_2\lraf{\tau_2}T_3\lra\cdots$ which is embedded into the following tower of triangles in $\Stmc{A}$:
$$
\xymatrix{
&T_1\ar[r]^-{h_1}\ar[d]_-{\tau_1}& X\ar[r]^-{g_1}\ar@{=}[d]& X_1\ar[r]\ar[d]_-{g_2}& T_1[1]\ar[d]_-{\tau_1[1]}\\
(\star\star)&T_2\ar[r]^-{h_2}\ar[d]_-{\tau_2}& X\ar[r]^-{g_1g_2}\ar@{=}[d]& X_2\ar[r]\ar[d]_-{g_3}& T_2[1]\ar[d]_-{\tau_2[1]}\\
&T_3\ar[r]^-{h_3}\ar[d]& X\ar[r]^-{g_1g_2g_3}\ar@{=}[d]& X_3\ar[r]\ar[d]& T_3[1]\ar[d]\\
&\vdots &\vdots&\vdots &\vdots
}
$$
Applying the Octahedral Axiom of triangulated categories for each $n$ yields a series of triangles
$$(\diamond\diamond) \quad T_n\lraf{\tau_n} T_{n+1}\lraf{\sigma_n} Q_{n+1}\lra T_n[1].$$ This implies $T_n\in \Add(\mathfrak{A})^{\star n}\subseteq \mathscr{X}$ for $n\geq 1$.

Let $\underrightarrow{\hocolim}(T_n)$ be the \emph{homotopy colimit}
in $\Stmc{A}$ of the tower of objects
$$
T_1\lraf{\tau_1}T_2\lraf{\tau_2}T_3\lraf{\tau_3}\cdots\lra T_n\lraf{\tau_n}T_{n+1}\lra\cdots
$$
defined by the triangle
$$
(\sharp)\quad\quad \bigoplus_{n\geq 1}T_n\lraf{(1-\tau_*)} \bigoplus_{n\geq 1}T_n\lra \underrightarrow{\hocolim}(T_n)\lra \bigoplus_{n\geq 1}T_n[1]
$$
where the morphism $(1-\tau_*)$  is induced by  $({\rm Id}_{T_n}, -\tau_n): T_n\to T_n\oplus T_{n+1\hookrightarrow }\bigoplus_{n\geq 1}T_n$.
Now, we choose $\tau_n$ as a representative in $A\Modcat$ and denote by $\lim\limits_{\longrightarrow}T_n$ the colimit of the direct system $\{(T_n, \tau_n)\mid n\geq 1\}$ of $A$-modules.
Then there is a short exact sequence
$$
0\to \bigoplus_{n\geq 1}T_n\lraf{(1-\tau_*)} \bigoplus_{n\geq 1}T_n\lra \lim\limits_{\longrightarrow}T_n\lra 0$$
which induces a canonical triangle in $\Stmc{A}$:
$$
\bigoplus_{n\geq 1}T_n\lraf{(1-\tau_*)} \bigoplus_{n\geq 1}T_n\lra \lim\limits_{\longrightarrow}T_n\lra \bigoplus_{n\geq 1}T_n[1].
$$
This implies that  $\underrightarrow{\hocolim}(T_n)\simeq \lim\limits_{\longrightarrow}T_n$ in $\Stmc{A}$.
Since the homotopy colimit of a tower of triangles in $\Stmc{A}$ is a triangle (for example, see \cite[Lemma 4.2]{Rk}), there is a triangle in $\Stmc{A}$:
$$(\ddag)\quad\quad
\underrightarrow{\hocolim}(T_n)\lra X \lra \underrightarrow{\hocolim}(X_n)\lra \underrightarrow{\hocolim}(T_n)[1].
$$
By the proof of \cite[Theorem III.2.3 and Remark III.2.7]{BI}, we can show $\underrightarrow{\hocolim}(T_n)\in\mathscr{X}$ and $\underrightarrow{\hocolim}(X_n)\in\mathscr{Y}$.
Since $(\mathscr{X}, \mathscr{Y})$ is a torsion pair in $\Stmc{A}$,  there are isomorphisms
$${\bf R}(X)\simeq\underrightarrow{\hocolim}(T_n)\;\;\mbox{and}\;\;{\bf L}(X)\simeq\underrightarrow{\hocolim}(X_n).$$
Recall that $T: \mathscr{D}\to\mathscr{C}$ stands for the left adjoint of the inclusion $\mathscr{C}\to \mathscr{D}$ (see Lemma \ref{Cotorsion pair}(3)) and $N$ is defined to be $H_*^0(M)$.

\begin{Lem}\label{R-L}
$(1)$ If $X\in {^{\bot >0}M}$, then ${\bf R}(X)$ lies in $\underline{\mathscr{G}\cap\mathscr{S}}$ and is isomorphic in $\mathscr{C}$ to $T(T_3)$.

$(2)$ If $X\in\mathscr{G}$, then $ {\bf L}(X)\in\underline{M}^{\bot}$. If further $X\in\mathscr{E}$, then
$T_3\in \stmc{A}$.

$(3)$  $\Omega_M^-(N)\simeq\Omega_A^-(N)\simeq T(\Omega_A^-(M))$ in $\mathscr{C}$.

\end{Lem}
{\it Proof.}  $(1)$ Recall from Section \ref{Torsion pair} that $\mathscr{X}=\underline{M}^{\bot >0}\cap \underline{\mathscr{S}}$ and $\mathscr{Y}={^{\bot \geq -1}}\underline{M}\subseteq {^{\bot >0}}\underline{M}=\underline{^{\bot >0}M}$.  Let $X\in {^{\bot >0}M}$. Since there is a triangle ${\bf R}(X)\to X\to {\bf L}(X)\to {\bf R}(X)[1]$ in $\Stmc{A}$ with ${\bf R}(X)\in\mathscr{X}$ and ${\bf L}(X)\in\mathscr{Y}$, we see that ${\bf R}(X)$ lies in  $\underline{^{\bot >0}M}$, and therefore in $\underline{\mathscr{G}\cap\mathscr{S}}$.

Set $\Lambda_n:=\{i\in\mathbb{N}\mid \StHom_A(M[i], X_n)\neq 0\}$ for $n\ge 0$. Since $\add({_A}M)=\add(\nu_A(M))$, it is clear that $\underline{{^{\bot >0}}M}=\underline{M}^{\bot \leq -2}$ by $(\lozenge)$.  This implies $\Lambda_0\subseteq \{0,1\}$. So we can choose $Q_1=M^{(I_{1, 0})}\oplus M[1]^{(I_{1,1})}$ in $\Stmc{A}$ for some index sets $I_{1, 0}$ and $I_{1,1}$. For a natural number $n$, we apply $\StHom_A(M[i], -)$ for $i\geq 0$ to the triangle
$$Q_{n+1}\lraf{f_{n+1}} X_n\lraf{g_{n+1}}X_{n+1}\lra Q_{n+1}[1]$$ in $\Stmc{A}.$  This yields an exact sequence of abelian groups:
$$
\StHom_A(M[i], Q_{n+1})\lraf{(f_{n+1})^*}\StHom_A(M[i], X_n)\lraf{(g_{n+1})^*} \StHom_A(M[i], X_{n+1})\lra \StHom_A(M[i], Q_{n+1}[1]).
$$
Since $f_{n+1}$ is a right $\Add(\mathfrak{A})$-approximation  of $X_n$, the map $(f_{n+1})^*$ is always surjective, and therefore there is an injection
$\StHom_A(M[i], X_{n+1})\hookrightarrow\StHom_A(M[i], Q_{n+1}[1])$ for $i\geq 0$. Using the fact $\Add(\underline{M})\subseteq\underline{\mathscr{G}}=\underline{M}^{\bot\neq 0,-1}$, we then can show $\Lambda_n\subseteq \{i\in\mathbb{N}\mid n\leq i\leq 2n+1\}$ by induction on $n$ and choose $Q_{n+1}=\bigoplus_{j\in\Lambda_n}M[j]^{(I_{n+1,j})}$ for some index sets $I_{n+1,j}$.

Consider $n\geq 3$. Then $\StHom_A(Q_{n+1}, M[1])=0$. It follows from ($\diamond\diamond$) and Lemma \ref{M-stable category}(2) that there is a triangle $T_n\to T_{n+1}\to Q_{n+1}\to \Omega_M^-(T_n)$ in $\mathscr{D}$.  We apply the functor $T$ to this triangle and produce another triangle $T(T_n)\to T(T_{n+1})\to T(Q_{n+1})\to \Omega_M^-(T(T_n))$ in $\mathscr{C}$. Since $\StHom_A(M[m], \mathscr{G})=0$ for any $m\geq 2$, the minimal injective envelope $M[m]\to I^m$ of $M[i]$ is a left $\mathscr{G}$-approximation of $M[m]$.  This means $T(M[m])\simeq I^m=0$ in $\mathscr{C}$ for any $m\geq 2$. Since $T$ commutes with direct sums, we have $T(Q_{n+1})\simeq\bigoplus_{j\in\Lambda_n}T(M[j])^{(I_{n+1,j})}=0$. Thus $T(T_n)\simeq T(T_{n+1})$. Consequently, the homotopy colimit  $\underrightarrow{\hocolim}(T(T_n))$ in $\mathscr{C}$ of  $\{(T(T_n), T(\tau_n)\mid n\geq 1\}$ is isomorphic  to $T(T_3)$.  Moreover, due to ${\bf R}(X)\in {^{\bot 1}M}$,  we see from Lemma \ref{M-stable category}(2) that ${\bf R}(X)$ is also the homotopy colimit in $\mathscr{D}$ of $\{T_n\}$.  Since ${\bf R}(X)\in\mathscr{G}$ and $T$ commutes with homotopy colimits by Lemma \ref{Cotorsion pair}(3), there are isomorphisms ${\bf R}(X)\simeq T({\bf R}(X))\simeq \underrightarrow{\hocolim}(T(T_n))\simeq T(T_3)$ in $\mathscr{C}$.

$(2)$  If $X\in\mathscr{G}$,  then ${\bf L}(X)\in \underline{M}^{\bot >0}$. This implies $ {\bf L}(X)\in\underline{M}^{\bot}$ because ${\bf L}(X)\in\mathscr{Y}=\underline{M}^{\bot\leq 0}$.
Let $$\mathscr{A}:=\{U\in \Stmc{A}\mid \bigoplus_{n=0}^{\infty}\StHom_A(M[n], U)\in\Gamma\modcat\}.$$  Since $\Gamma$ is an Artin algebra
over a commutative Artin ring $R$, a $\Gamma$-module $N$ is finitely generated if and only if $N$ is a finitely generated  $R$-module.
This implies that an $A$-module $U$ lies in $\mathscr{A}$ if and only if $\bigoplus_{n=0}^{\infty}\StHom_A(M[n], U)\in R\modcat$. The latter is
equivalent to saying that there is a non-negative integer $\delta_U$ such that $\StHom_A(M[n], U)=0$ for $n>\delta_U$ and $\StHom_A(M[n], U)\in R\modcat$ for $0\leq n\leq \delta_U$. Moreover, $\mathscr{A}$ is closed under direct summands, finite direct sums, the shift $[-1]$ and extensions of triangles in $\Stmc{A}$. Now, it follows from $M\in\underline{M}^{\bot\neq 0,-1}$ and $\Hom_A(M[1], M)\simeq D\StHom_A(M, \nu_A(M))$ that $M[j]\in\mathscr{A}$ for any $j\in\mathbb{Z}$. Consequently, for $U\in\mathscr{A}$, there is a right $\Add(\mathfrak{A})$-approximation $f_U: Q_U\to U$ of $U$ such that $Q_U\in\add(\bigoplus_{n=0}^m M[n])$ for some $m\geq 0$ and the third term $C_U$
of the triangle $Q_U\lraf{f_U} U\to C_U\to Q_U[1]$ in $\Stmc{A}$ still lies in $\mathscr{A}$.

By Lemma \ref{Equality}, there holds $\underline{\mathscr{E}}\subseteq \mathscr{A}$. If $X\in\mathscr{E}$, then $I_{n+1,j}$ can be chosen to be finite sets and $Q_{n+1}\in\add(\bigoplus_{j\in\Lambda_n}M[j])\subseteq\stmc{A}$.

$(3)$ In the proof of $(1)$, we take $X=N$. Then $X\in\mathscr{E}$ by Theorem \ref{Compact}(1). Recall that $X\in\mathscr{H}\subseteq\underline{M}^{\bot\neq 0}$ by Proposition \ref{Self-orthogonal}, and $\StHom_A(M, M)\simeq\StHom_A(M, X)$ by Lemma \ref{Isomorphism}(1). Further, we even have

$(a)$ $\Lambda_0\subseteq\{0\}$, $Q_1=M$ and $X_1[-2]\in\mathscr{H}$ (see the proof of Lemma \ref{CR}); and

$(b)$ $\Lambda_n\subseteq \{i\in\mathbb{N}\mid n+1\leq i\leq 2n\}$ and $Q_{n+1}\in\add(\bigoplus_{i=n+1}^{2n}M[i])\in\stmc{A}$ for $n\geq 1$.

\medskip
Since $M$ lies in $\underline{\mathscr{G}}=\underline{M}^{\bot\neq 0,-1}$, there holds $\StHom_A(Q_3, M[1])=0$. By Lemma \ref{M-stable category}(2), there exists a right triangle $T_2\lraf{\tau_2} T_3\lraf{\sigma_2} Q_3\ra \Omega_M^-(T_2)$ in $\mathscr{D}$. This gives rise to another triangle $T(T_2)\to T(T_3)\to T(Q_3)\to \Omega_M^-(T(T_2))$ in $\mathscr{C}$. Clearly, $T(Q_3)=0$ by $Q_3\in\add(M[3]\oplus M[4])$. Thus
$T(T_2)\simeq T(T_3)$. It then follows from (1) and $X\in\mathscr{H}\subseteq \mathscr{X}$ that $X \simeq {\bf R}(X)\simeq T(T_2)$ in $\mathscr{C}$.

By Lemma \ref{CR}, $X[1]=\Omega_A^-(X)\simeq \Omega_M^-(X)$ in $\mathscr{C}$. Further, we will show $T_2[1]=\Omega_A^-(T_2)\simeq \Omega_M^-(T_2)$ in $\mathscr{C}$. Actually, it suffices to prove that any homomorphism from $T_2$ to a module in $\Add(M)$ factorizes through the injective envelope of $T_2$, or equivalently, $\StHom_A(T_2, M)=0$. Thanks to $\add(M)=\add(\nu_A(M))$, we will show $\StHom_A(T_2, \nu_A(M))=0$.

Indeed, $\StHom_A(T_2, \nu_A(M))\simeq D\StHom_A(M[1], T_2)$  by $(\lozenge)$.  We apply $\StHom_A(M[1],-)$ to the triangle
$X_2[-1]\to T_2\to X\to X_2$ (see ($\star\star$)) and obtain an isomorphism $\StHom_A(M[1], X_2[-1]) \simeq \StHom_A(M[1], T_2)$. Here we use the fact $X\in\mathscr{H}\subseteq\underline{M}^{\bot\neq 0}$. So, it is enough to show $\StHom_A(M[1], X_2[-1])=0$. Recall that the triangle $Q_2\lraf{f_2} X_1\to X_2\to Q_2[1]$ has the following properties:  $Q_2\in\add(M[2])$ by $(b)$, $X_1\in\mathscr{H}[2]$ by $(a)$, and there exists an injection $\StHom_A(M[2], X_2)\hookrightarrow\StHom_A(M[2], Q_2[1])$ by the proof of $(1)$. It follows from $\StHom_A(M, M[1])=0$ that $\StHom_A(M[1], X_2[-1])\simeq \StHom_A(M[2], X_2)=0$. Thus  $\StHom_A(M[1], T_2)=0$, and therefore $\StHom_A(T_2, \nu_A(M))=0$.

By the facts that $Q_2\in\add(M[2])$ and $\StHom_A(M[2], M)=0$, we have $$\Ext_A^1(Q_2[1], M) \simeq \StHom_A(Q_2[1], M[1]) \simeq \StHom_A(Q_2, M)=0.$$By Lemma \ref{M-stable category}(2), the triangle $M[1]\lraf{\tau_1[1]} T_2[1]\lraf{\sigma_1[1]} Q_2[1]\lra M[2]$ in $\Stmc{A}$ can be extended to a right triangle $M[1]\lraf{\tau_1[1]} T_2[1]\lraf{\sigma_1[1]} Q_2[1]\lra \Omega_M^{-}(M[1])$ in $\mathscr{D}$. This yields the triangle in $\mathscr{C}$
$$T(M[1])\lraf{T(\tau_1[1])} T(T_2[1])\lraf{T(\sigma_1[1])} T(Q_2[1])\lra \Omega_M^{-}(T(M[1]))$$ by Lemma \ref{Cotorsion pair}(3).
It then follows from $T(M[3])=0$ that $T(Q_2[1])=0$  and  $T(M[1])\simeq T(T_2[1])$ in $\mathscr{C}$. Moreover, $T(T_2[1])\simeq T(\Omega_M^-(T_2))\simeq \Omega_M^-(T(T_2))\simeq \Omega_M^-(X)$ in $\mathscr{C}$, where the second isomorphism is due to the fact that $T$ commutes with $\Omega_M^-$. Thus $T(\Omega_A^-(M))=T(M[1])\simeq\Omega_M^-(X)$ in $\mathscr{C}$. $\square$

\medskip
Now, we state a property of finitely $M$-filtered $A$-modules introduced in Definition \ref{Countably filtered}.

\begin{Lem}\label{FMF}
If $X$ is a finitely $M$-filtered $A$-module, then $X\in \mathscr{G}^{\bot >0}$ and has an $\add(M)$-resolution
of finite length. Moreover, if the $A$-module $X$ lies in $\mathscr{G}$, then $X\in\add(M)$.
\end{Lem}
{\it Proof.}
Since $\mathscr{G}^{\bot >0}$ contains $M$ and is closed under cosyzygies and extensions in $A\Modcat$, we have $X\in \mathscr{G}^{\bot >0}$.
Now, let $\mathscr{M}$ be the full subcategory of $A\modcat$ consisting of all those modules having an $\add(M)$-resolution of finite length. Then $Y\in\mathscr{M}$ if and only if $\Hom_A(M, Y)\in \mathscr{P}^{<\infty}(\Lambda)$. Since ${_A}M$ is a self-orthogonal generator, $Y\in\mathscr{M}$ if and only if there is an exact sequence $0\to M_n\to M_{n-1}\to\cdots\to M_0\to Y\to 0$ in $A\modcat$ for some $n\in\mathbb{N}$ such that $M_j\in\add(M)$ for all $0\leq j\leq n$. Clearly, $\mathscr{M}\subseteq M^{\bot >0}\subseteq M^{\bot 1}$.  As $\mathscr{P}^{<\infty}(\Lambda)$ is always closed under extensions in $\Lambda\modcat$, $\mathscr{M}$ is closed under extensions in $A\modcat$. So, to show that the finitely $M$-filtered module $X$ belongs to $\mathscr{M}$, it suffices to show $\Omega_A^{-i}(M)\in\mathscr{M}$ for each $i\geq 0$. However, this follows from the exact sequence $0\to M\to I^0\to\cdots\to I^{i-1}\to \Omega_A^{-i}(M)\to 0$, where $I^j$ is injective and therefore in $\add(M)$ for $0\leq j\leq i-1$.
As to the last statement in the lemma, we notice $\mathscr{M}\cap{^{\bot >0}}M=\add(M)$, $\mathscr{G}\subseteq {^{\bot >0}}M$ and $\mathscr{G}\cap\mathscr{M}=\add(M)$. $\square$

\medskip
In the following, we describe zero objects and compact objects in $\mathscr{C}$.
This is related to pure-projective modules. It is known that a module over an Artin algebra is \emph{pure-projective} if and only if it is a direct sum of finitely generated, indecomposable modules.

\begin{Prop}\label{CP}
The following hold for $X\in\mathscr{E}\cap\mathscr{S}$.

$(1)$ $X$ is compact in $\mathscr{C}$ and isomorphic in $\Stmc{A}$ to an $M$-filtered module.

$(2)$ If $X^{\bot 1}$ is closed under direct sums of countably many finitely $M$-filtered $A$-modules in $A\Modcat$, then $X\in\Add(M)$.
In particular, if $X$ is pure-projective, then $X\in\Add(M)$.

$(3)$ If $X\in\mathscr{H}\cap\stmc{A}$, then $X$ is projective.
\end{Prop}

{\it Proof.} Let $X\in\mathscr{E}\cap\mathscr{S}$. We keep all notations in the proof of Lemmas \ref{R-L}(1)-(2).

$(1)$ By Lemma \ref{Equality}, $X\in\mathscr{X}$ and $X\simeq {\bf R}(X)$ in $\Stmc{A}$.
Moreover, by Lemmas \ref{R-L}(1)-(2), ${\bf R}(X)\simeq T(T_3)$ in $\mathscr{C}$ with $T_3\in\stmc{A}$.  Each object of $\stmc{A}$ is compact in $\Stmc{A}$ and the quotient functor $\Stmc{A}\to \mathscr{D}$ preserves compact objects by Lemma \ref{AJCP}(2). Hence $T_3$ is compact in $\mathscr{D}$. As $T$ preserves compact objects by Lemma \ref{Cotorsion pair}(3), both $T(T_3)$ and $X$ are compact in $\mathscr{C}$.

It follows from $X\in\mathscr{E}$ that the proofs of Lemmas \ref{R-L}(1)-(2) yield $Q_{n+1}=\bigoplus_{j\in\Lambda_n}M[j]^{(I_{n+1,j})}$, where $\Lambda_n\subseteq \{i\in\mathbb{N}\mid n\leq i\leq 2n+1\}$, $M[j]=\Omega_A^{-j}(M)$ and $I_{n+1,j}$ are finite sets. This implies $T_n\in \stmc{A}$ for $n\geq 1$. Since triangles in $\stmc{A}$ (up to isomorphism) are induced from short exact sequences, we can add finitely generated projective modules to $T_n$ and assume that the associated map $\tau_n: T_n\to T_{n+1}$ is injective and $\Coker(\tau_n)$ is isomorphic to a finite direct sum of modules in $\{A\}\cup\{\Omega_A^{-i}(M)\mid i\in\mathbb{N}\}$. Now, let $X'$ be the colimit of the direct system $\{(T_n, \tau_n)\mid n\geq 1\}$ in $A\Modcat$. Then the canonical map $T_n\to X'$ is injective. Thus $T_n$ can be regarded as a submodule of $X'$ and $T_n\subseteq T_{n+1}$. Moreover, $X'=\bigcup_{n=0}^{\infty}T_n$ and is $M$-filtered. Now (1) follows from $X\simeq {\bf R}(X)\simeq X'$ in $\Stmc{A}$.

$(2)$ Since all $T_n$ are finitely $M$-filtered and $X\in\mathscr{G}$, we have $T_n\in X^{\bot 1}$ by Lemma \ref{FMF}. Let $Z:=\bigoplus_{n\geq 1}T_n$.
The assumption of $(2)$ implies $\Ext_A^1(X, Z)=0$. So the exact sequence
$0\to Z\lraf{(1-\tau_*)}Z\to X'\to 0$ induced by $\{\tau_n\mid n\geq 1\}$ splits, and therefore $X'$ is isomorphic to a direct summand of
$Z$. Since $T_n$ is a finitely generated $A$-module, it is a direct sum of finitely many indecomposable submodules with the local endomorphism rings. It follows from \cite[Corollary 26.6]{AF} that $X'\simeq \bigoplus_{n\geq 1}T_n'$, where $T_n'$ is a direct summand of $T_n$.
Moreover, by Lemma \ref{FMF}, $T_n$ has an $\add(M)$-resolution of finite length, and so does $T_n'$. It follows from $X\simeq X'$ in $\Stmc{A}$ and $X\in\mathscr{G}$ that $X'\in\mathscr{G}$, and therefore $T_n'\in\mathscr{G}$. Consequently, $T_n'\in\add(M)$ by the proof of Lemma \ref{FMF}. Hence $X'\in\Add(M)$ and $X\in\Add(M)$.

For any $Z\in A\modcat$, $Z$ is clearly finite presented. So $Z^{\bot 1}$ is always closed under arbitrary direct sums in $A\Modcat$. This implies that if $X$ is pure-projective, then $X^{\bot 1}$ is closed under arbitrary direct sums in $A\Modcat$, and therefore $X\in \Add(M)$.

$(3)$ By Lemmas \ref{Equality}(1) and \ref{E-property}(1), $\mathscr{H}\cap\stmc{A}\subseteq \underline{\mathscr{E}\cap\mathscr{S}}$. It follows from $X\in\stmc{A}$ that $X^{\bot 1}$ is closed under arbitrary direct sums in $A\Modcat$. Further, $X\in\Add(M)$ by $(2)$, and $\StHom_A(N, M)=0$ by Lemma \ref{CR}. Due to $\Add(M)=\Prod(M)$, we have  $\StHom_A(N, X)=0$. It follows from Lemma \ref{Isomorphism}(3) and $X\in\mathscr{H}$ that $X\simeq 0$ in $\Stmc{A}$. Thus $X$ is projective. $\square$

\medskip
Now, we describe compact generators of the right-hand side in the recollement of Theorem \ref{Gorenstein-recollement}.

Let $\mathcal{T}$ be a triangulated category.  For a full subcategory $\mathcal{U}$ of $\mathcal{T}$, we denote by $\langle\mathcal{U}\rangle$ the smallest full subcategory of $\mathcal{T}$ containing $\mathcal{U}$ and being closed under finite coproducts, direct summands and shifts.
Let $\mathcal{V}$ be another full subcategory of $\mathcal{T}$. Define $\mathcal{U}\star \mathcal{V}$ to be the full subcategory of $\mathcal{T}$ consisting of objects $X$ such that there is a triangle $U\to X\to V\to U[1]$ with $U\in\mathcal{U}$ and $V\in\mathcal{V}$.
Following \cite[3.1]{RR}, we set $\mathcal{U}\diamond\mathcal{V}:=\langle\mathcal{U}\star \mathcal{V}\rangle$, and then
define $\langle\mathcal{U}\rangle_0:=0$ and $\langle\mathcal{U}\rangle_{n+1}:=\langle\mathcal{U}\rangle_n\diamond\langle\mathcal{U}\rangle$ for $n\ge 0$ inductively. Clearly, the objects of $\langle\mathcal{U}\rangle_{n+1}$ are the \emph{direct summands} of the objects obtained by taking $(n+1)$-fold extensions of finite direct sums of shifts of objects of $\mathcal{U}$. Following \cite[Definition 3.1]{RR}, the \emph{dimension} of $\mathcal{T}$, denoted by $\dim(\mathcal{T})$, is defined to be the minimal natural number $n$
such that there is an object $X\in\mathcal{T}$ such that $\mathcal{T}=\langle X\rangle_{n+1}$. If no such $X$ exists, one defines $\dim(\mathcal{T})=\infty$. If $\dim(\mathcal{T})$ is finite, then such an object $X$ is called a \emph{strong generator} of $\mathcal{T}$.

For a pair of integers $i\leq j$, we  denote by $\langle\mathcal{U}\rangle_{n+1}^{[i,j]}$ the full subcategory of $\langle\mathcal{U}\rangle_{n+1}$ consisting of all objects which are obtained by taking $(n+1)$-fold extensions of finite direct sums of objects in the class
$\{U[-s]\mid U\in\mathcal{U}, s\in\mathbb{Z}, i\leq s\leq j\}$. Here, we do not require taking both direct summands
and arbitrary shifts in $\langle\mathcal{U}\rangle_{n+1}^{[i,j]}$.

\begin{Koro} \label{Compact recollement}
$(1)$ Let $\mathcal{S}$ be the finite set of isomorphism classes of simple objects of $\mathscr{H}$, and let $n$ be the Loewy length of the algebra $\Gamma$.
Then $(\mathscr{G}\cap\mathscr{S})/[M]$ is a compactly generated triangulated category:
$$\big(\mathscr{G}\cap\mathscr{S}\big)/[M]=\langle\Add(\mathcal{S})\rangle_{2n}^{[-1,0]}\;\;\mbox{and}\;\;
\big((\mathscr{G}\cap\mathscr{S}\big)/[M])^{\rm c}=(\mathscr{E}\cap\mathscr{S})/[M]=\langle\mathcal{S}\rangle_{2n}^{[-1,0]}.$$

$(2)$ If $A$ is symmetric and $\Gamma$ is semisimple, then there are triangle equivalences
$$(\mathscr{G}\cap\mathscr{S})/[M]=\Add\big(N\oplus \Omega_A^-(N)\big)\lraf{\simeq}(\Gamma\times\Gamma)\Modcat,$$
$$(\mathscr{E}\cap\mathscr{S})/[M]=\add\big(N\oplus \Omega_A^-(N)\big)\lraf{\simeq}(\Gamma\times\Gamma)\modcat,$$
where $(\Gamma\times\Gamma)\Modcat$, as a triangulated category, has the shift functor induced from the automorphism of the algebra
$\Gamma\times\Gamma$ by permutating the first and second coordinate.

$(3)$ If $T(X)\in\mathscr{E}/[M]$ for any $X\in A\modcat$, then the recollement in Proposition \ref{Restriction} restricts to
a half recollement of triangulated categories:
$$
\xymatrix@C=1.2cm{(\underline{M}^{\bot})^{\rm c}\ar[r]^-{{\pi\circ \rm inc}}
&\mathscr{C}^{\rm c}\ar[r]^-{\widetilde{\Phi}}\ar_-{\widetilde{\Psi}}@/_1.3pc/[l]
&(\mathscr{E}\cap\mathscr{S})/[M].\ar@/_1.2pc/[l]_{{\rm inc}}}
$$
\end{Koro}

{\it Proof.} $(1)$ Let $\mathscr{R}:=(\mathscr{G}\cap\mathscr{S})/[M]$. By Lemma \ref{Equality}(2)(b), we can identify $\mathscr{H}$ with its (essential) image in $\mathscr{R}$ under the quotient functor $\underline{\mathscr{G}\cap\mathscr{S}}\to\mathscr{R}$. It follows from the first part of Corollary \ref{Gamma-mod}(3) that $\mathscr{R}=\mathscr{H}[1]\star \mathscr{H}$, where  $[1]$ denotes the functor $\Omega_M^{-}$.  To characterize $\mathscr{R}$ in terms of $\mathcal{S}$, we use the functor $\StHom_A(N, -): \mathscr{H}\to \Gamma\Modcat$ which is an equivalence of abelian categories sending $N$ to $\Gamma$, due to Lemma \ref{Isomorphism}(2)(3). Recall that $\mathscr{H}^{\rm fg}:=\{X\in\mathscr{H}\mid \StHom_A(M, X)\in\Gamma\modcat\}=\mathscr{H}\cap\underline{\mathscr{E}}$. By the proof of Corollary \ref{Gamma-mod}(1), the functor $\StHom_A(N, -)$ restricts to an equivalence: $\mathscr{H}^{\rm fg}\lraf{\simeq}\Gamma\modcat$. This equivalence clearly sends simple objects of $\mathscr{H}$ to simple $\Gamma$-modules. Since $\Gamma$ is an Artin algebra, it has only finitely many isomorphism classes of simple modules and each (respectively, finitely generated)
$\Gamma$-module is generated by simple modules under arbitrary (respectively, finite) direct sums and taking $n$-fold extensions. Consequently, $\mathcal{S}$ is a finite set and each object of $\mathscr{H}$ (respectively, $\mathscr{H}^{\rm fg}$) is generated by $\mathcal{S}$ under arbitrary (respectively, finite) direct sums and taking $n$-fold extensions. Note that each short exact sequence in $\mathscr{H}$ induces a triangle in $\mathscr{R}$ by Lemma \ref{M-stable category}(2). Now, by $\mathscr{R}=\mathscr{H}[1]\star \mathscr{H}$, the first equality in $(1)$ holds, and therefore $\mathscr{R}=\Loc(\mathcal{S})$. Similarly, from the equality $(\mathscr{E}\cap\mathscr{S})/[M]=\mathscr{H}^{\rm fg}[1]\star \mathscr{H}^{\rm fg}$ by the second part of Corollary \ref{Gamma-mod}(3), we see that $(\mathscr{E}\cap\mathscr{S})/[M]$ is generated by $\mathcal{S}$ under taking the shifts $[i]$ for $i=0,1$ and $(2n)$-fold extensions.  This implies the third equality in $(1)$. By Proposition \ref{CP}(1), each object of $\mathscr{H}^{\rm fg}$ is compact in $\mathscr{C}$ and thus also in $\mathscr{R}$. Hence $\mathscr{R}$ is compactly generated. By \cite[Theorem 4.4.9]{Neemanbook}, we have $\mathscr{R}^{\rm c}=\thick(\mathcal{S})$.
Clearly, $\thick(\mathcal{S})\subseteq(\mathscr{E}\cap\mathscr{S})/[M]=\langle\mathcal{S}\rangle_{2n}^{[-1,0]}\subseteq \thick(\mathcal{S})$. Thus $\mathscr{R}^{\rm c}=\langle\mathcal{S}\rangle_{2n}^{[-1,0]}$.

$(2)$ Suppose that $A$ is symmetric and $\Gamma$ is semisimple. Then $\Gamma$ is symmetric. By Corollary \ref{2-periodic}, $\Omega_M^2(N)\simeq N$ in $\mathscr{C}$.
Moreover, by Lemma \ref{CR}, $\Omega_M^-(N)\simeq \Omega_A^-(N)$ in $\mathscr{C}$. It follows that $\Omega_M^{2i}(N)\simeq N$ and $\Omega_M^{2i+1}(N)\simeq \Omega_A^-(N)$ for any $i\in\mathbb{Z}$. Let $W:=N\oplus \Omega_A^-(N)$ and $\mathscr{W}:=\Add(W)$. Then $\mathscr{W}=\Omega_M(\mathscr{W})$. Clearly, $\mathscr{W}$ contains $N$ and is closed under direct sums in $\mathscr{C}$. As $\Gamma$ is semisimple, there holds $\add(N)=\add(\mathcal{S})\subseteq\mathscr{R}$. To show $\mathscr{R}=\mathscr{W}$, it suffices to show that $\mathscr{W}$ is a triangulated subcategory of $\mathscr{C}$.
Since $\mathscr{W}$ is closed under $\Omega_M$, we only need to show that, for any morphism $f: X_1\to X_2$ in $\mathscr{W}$ and  triangle $X_1\stackrel{f}{\ra} X_2\ra X_3\ra \Omega_M^-(X_1)$ in $\mathscr{C}$, the term $X_3$ belongs to $\mathscr{W}$.

Since $\mathscr{H}$ is the heart of the torsion pair $(\mathscr{X},\mathscr{Y})$ in $\Stmc{A}$, $\Ext_{\mathscr{H}}^j(U, V)\simeq \StHom_A(U, V[j])$ for any $U, V\in \mathscr{H}$ and $j=0, 1$. Since $N$ is a projective object in $\mathscr{H}$ by Lemma \ref{Isomorphism}(3), we have $\StHom_A(N, N[1])\simeq\Ext_{\mathscr{H}}^1(N, N)=0$. Further, by Lemma \ref{Gorenstein} and \cite[Lemma 3.5]{cs}, if $V_1$ and $V_2$ lie in $\mathscr{G}$, then $\StHom_A(V_1, V_2[n])$ $\simeq \StHom_M(V_1, \Omega_M^{-n}(V_2))$ for all $n\geq 1$. This implies $\StHom_A(N, N[1])\simeq \StHom_M(N, \Omega_M^-(N))$, and therefore $\StHom_M(N, \Omega_M^-(N))=0=\StHom_M(\Omega_M^-(N), N)$. Let $B:=\StEnd_A(W)$. Then $B\simeq \Gamma \oplus \Gamma$ as algebras.
Since $B$ is semisimple and $W$ is compact in $\mathscr{C}$ by Lemma \ref{Compact}(1), the functor $\StHom_M(W,-): \mathscr{W}\to B\Modcat$ is an additive equivalence. By this equivalence and the fact that $B\Modcat$ is a semisimple abelian category, it can be proved that $f$ as a morphism in $\mathscr{W}$ is isomorphic to a direct sum of the identity map of $Z_1$ with the zero map $Z_2\to Z_3$, where $Z_i\in \mathscr{W}$ for $1\leq i\leq 3$.
Consequently, $X_3\simeq \Omega_M^{-}(Z_2)\oplus Z_3$. It then follows from $\mathscr{W}=\Omega_M(\mathscr{W})$ that $X_3\in\mathscr{W}$. Thus $\mathscr{W}$ is a triangulated subcategory of $\mathscr{C}$, and therefore $\mathscr{R}=\mathscr{W}$.
By $(1)$, $\mathscr{R}^{\rm c}=\add(W)$.

$(3)$ By Lemma \ref{Cotorsion pair}(4) and \cite[Theorem 4.4.9]{Neemanbook}, the assumption of Corollary \ref{Compact recollement}(3) implies
$\mathscr{C}^{\rm c}\subseteq \mathscr{E}/[M]$. Combining Corollary \ref{Compact-reco} with Proposition \ref{CP}(1), we see that $\widetilde{\Phi}$ sends compact objects of $\mathscr{C}$ to objects of $(\mathscr{E}\cap\mathscr{S})/[M]$ which are also compact in $\mathscr{C}$. Note that
$\widetilde{\Psi}:\mathscr{C}\to \underline{M}^{\bot}$ always preserves compact generating sets. Thus the first two lines of functors in the recollement of Corollary \ref{Compact-reco} restrict to the half recollement in Corollary \ref{Compact recollement}(3). $\square$

\begin{Koro}\label{dimension}
$(1)$ $\dim\big((\mathscr{E}\cap\mathscr{S})/[M])\big)\leq 2\; LL(\Gamma)-1$, where $LL(\Gamma)$ is the Loewy length of $\Gamma$.

$(2)$ $\dim\big((\mathscr{E}\cap\mathscr{S})/[M])\big)\leq 2\,\gd(\Gamma)+1$, where $\gd(\Gamma)$ is the global dimension of $\Gamma$ .
\end{Koro}

{\it Proof.} $(1)$ follows from Corollary \ref{Compact recollement}(1).

$(2)$ If $\gd(\Gamma)$ is infinite, then the inequality in $(2)$ holds trivially. Now, let $m:= \gd(\Gamma)<\infty$. In the following, we follow the notation in the proof of Corollary \ref{Compact recollement}(1) and show that $(\mathscr{E}\cap\mathscr{S})/[M]=\langle\add(N)\rangle_{2m+2}^{[-m-1,0]}$. This implies
$\dim\big((\mathscr{E}\cap\mathscr{S})/[M])\big)\leq 2m+1$.

By the second part of Corollary \ref{Gamma-mod}(3), $(\mathscr{E}\cap\mathscr{S})/[M]=\mathscr{H}^{\rm fg}[1]\star \mathscr{H}^{\rm fg}$. So, it suffices to control the objects of
$\mathscr{H}^{\rm fg}$ by $N$. We take an object $X\in\mathscr{H}^{\rm fg}$. Since
the functor $\StHom_A(N, -):\mathscr{H}^{\rm fg}\to\Gamma\modcat$ is an equivalence of abelian categories sending $N$ to $\Gamma$ and since each finitely generated $\Gamma$-module has a projective resolution of length $m$ by finitely generated projective $\Gamma$-modules, there is a long exact sequence
$(\ast):\;0\to N_m\to \cdots\to N_1\to N_0\to X\to 0$ in $\mathscr{H}^{\rm fg}$ with $N_i\in\add(N)$ for $0\leq i\leq m$. Note that each short exact sequence
$0\to X\to Y\to Z\to 0$ in $\mathscr{H}$ gives rise to a triangle
$X\to Y\to Z\to \Omega_A^{-1}(X)$ in $\Stmc{A}$ with the terms $X,Y,Z\in\mathscr{H}$.
Since $\mathscr{H}\subseteq \underline{\mathscr{G}\cap\mathscr{S}}\subseteq\underline{{^{\bot 1}}M}$ by Lemma \ref{Equality}(1), this triangle induces a triangle $X\to Y\to Z\to X[1]$ in $\mathscr{C}$ by Lemma \ref{M-stable category}(2). So, we divide the sequence $(\ast)$ into a series of short exact sequences in $\mathscr{H}^{\rm fg}$
and then obtain $X\in\langle\add(N)\rangle_{m+1}^{[-m,0]}$.
Thus $(\mathscr{E}\cap\mathscr{S})/[M]=\langle\add(N)\rangle_{2m+2}^{[-m-1,0]}$. $\square$

\medskip
{\bf Proof of Theorem \ref{Gorenstein-recollement}.} The existence of the recollement in $(1)$ follows from Proposition \ref{Restriction}, while $(2)$ and (3) are exactly Corollaries \ref{Compact-reco} and  \ref{dimension}, respectively. $\square$

\medskip
{\bf Proof of Proposition \ref{Compact objects in recollement}.}  $(2)$ and $(3)$ follow from Corollary \ref{Compact recollement}(1), while the first part of $(1)$ is Proposition \ref{CP}(1). Let $\mathscr{E}_0:=\mathscr{E}\cap\mathscr{S}\cap A\modcat$. Clearly, $\add(M)\subseteq \mathscr{E}_0$. If $X\in\mathscr{E}_0$, then $X\in\Add(M)$ by Proposition \ref{CP}(2) because finitely generated $A$-modules are pure-projective.  As $X$ is finitely generated, it lies in $\add(M)$. Thus $\mathscr{E}_0=\add(M)$.
This implies the second part of $(1)$.  $\square$

\section{Tachikawa's second conjecture}\label{Characterization}

In this section we prove Theorem \ref{Conjecture} and Corollary \ref{Gendo-symmetric}. As a consequence, we show that the Nakayama conjecture is true for Gorenstein-Morita algebras (see Definition \ref{Gorenstein-Morita}(ii)). Moreover, we introduce two homological conditions ${\rm (CI)}$  and ${\rm (CII)}$). They are connected with finitistic dimension (see Lemma \ref{Three conditions}). Also, the invariance of ${\rm (CII)}$ under different types of equivalences between algebras is discussed in Corollary \ref{DM}.

Let $B$ be an Artin algebra. The \emph {dominant dimension} of  $B$, denoted by
$\dd(B)$, is by definition the largest natural number $n$ or $\infty$ such that, in a minimal injective coresolution
$0\to {_B}B\to I^0\to I^1\to\cdots\to I^n\to\cdots$,
all these module $I^i$ are projective for $0\leq i < n$.
The unsolved Nakayama Conjecture says that an Artin algebra is self-injective whenever its dominant dimension is infinite. Related to  this conjecture, Tachikawa proposed two conjectures in \cite[p.\! 115-116]{Tac}.

{\bf (TC1)}:
If an Artin algebra $B$ satisfies $\Ext^n_B(D(B),B)=0$ for all $n\geq 1$, then $B$ is self-injective.

{\bf (TC2)}:
Let $B$ be a self-injective algebra and $Y$ a finitely
generated $B$-module. If $\Ext^n_B(Y,Y)=0$ for all $n\geq 1$, then $Y$ is projective.

As pointed in the introduction, the two conjectures \textbf{(TC1)} and \textbf{(TC2)} hold true for all algebras if and only if so does the Nakayam conjecture for all algebras. Moreover, it was shown in \cite{Muller} that, given a pair $(B, Y)$ with $B$ a self-injective algebra and $Y$ a finitely generated, self-orthogonal  $B$-module, the algebra $\End_B(B\oplus Y)$ satisfies the Nakayama conjecture if and only if $Y$ is projective.

By \cite{Muller}, algebras of dominant dimension at least $2$ are exactly endomorphism algebras of
generator-cogenerators over algebras. Recall that a finitely generated module $X$ over an algebra $C$ is called
a \emph{generator-cogenerator} if $C\oplus D(C)\in\add({_C}X)$.
If $B$ is the endomorphism algebra of a generator-cogenerator $X$ over an algebra $C$, then, by M\"{u}ller's theorem (see \cite[Lemma 3]{Muller}), $\dd(B)=n>1$ if and only if $\Ext_C^i(X, X)=0$ for $1\leq i<n-1$.
In particular, $\dd(B)=\infty$ if and only if $X$ is self-orthogonal, that is, $\Ext_C^i(X, X)=0$ for all $i\geq 1$.

\medskip
{\bf Proof of Theorem \ref{Conjecture}.} $(1)\Rightarrow (2)\mbox{-}(5)$.
Suppose ${_A}M$ is projective. Then $\mathscr{G}=A\Modcat$ and $W=\Omega_A^-(M)=0$.
This implies $(4)$ and $(5)$. Since each $A$-module is always a filtered colimit of finitely generated $A$-modules, $(2)$ holds. Note that $A$-filtered modules are projective. Thus $(3)$ holds.

$(2)\Rightarrow (4)$. This is clear since $W\in\mathscr{G}$.

$(4)\Rightarrow (5)$.  Let $f:\Omega_A^-(M)\to W$ be the minimal left $\mathscr{G}$-approximation of $\Omega_A^-(M)$ with $W\in \mathscr{G}$. Assume that $W$ is a filtered colimit of $\{W_i\mid i\in I\}$ with $I$ an essentially small, filtered category and with all $W_i\in\mathscr{G}_0$. We show that $W^{\bot 1}$ is closed under countable direct sums in $A\Modcat$. This implies $(4)$.

Let $\lambda_i: W_i\to W$ be the canonical homomorphism of the colimit $W$. Since $\Omega_A^-(M)$ is finitely generated and even finitely presented, the canonical map $\lim\limits_{\longrightarrow} \Hom_A(\Omega_A^-(M), W_i)\to \Hom_A(\Omega_A^-(M), W)$ induced by
$\{\Hom_A(\Omega_A^-(M), \lambda_i)\mid i\in I\}$ is an isomorphism. As $I$ is a filtered category, the colimt $\lim\limits_{\longrightarrow} \Hom_A(\Omega_A^-(M), W_i)$ of abelian groups is a quotient group of the direct sum $\bigoplus_{i\in I}\Hom_A(\Omega_A^-(M), W_i)$ and each of its elements can be represented by a homomorphism of $\Hom_A(\Omega_A^-(M), W_i)$ for some index $i\in I$. Consequently, there is an index $i\in I$ and a homomorphism $f_i: \Omega_A^-(M)\to W_i$ such that $f=f_i\lambda_i$. Further, due to $W_i\in\mathscr{G}$, the approximation implies that there is a homomorphism $g_i: W\to W_i$ such that $f_i=fg_i$. It follows that $f=fg_if_i$. Since $f$ is left minimal, $g_if_i$ is an isomorphism. Thus $g_i$ is split-injective. Since $W_i$ is finitely generated, $W$ is also finitely generated (presented). Hence, $W^{\bot 1}$ is closed under arbitrary direct sums in $A\Modcat$.

$(5)\Rightarrow (1)$.  By Lemma \ref{R-L}(3), $T(\Omega_A^-(M))\simeq \Omega_M^-(N)$ in $\mathscr{C}$. Thus $W\simeq\Omega_M^-(N)$ in $\mathscr{C}$. By Lemma \ref{M-stable category}(1), there are $M_1, M_2\in\Add(M)$ such that $W\oplus M_1\simeq \Omega_M^-(N)\oplus M_2$. Let $\mathscr{M}$ be the full subcategory of $A\modcat$ consisting of modules with finite $\add(M)$-resolutions. By $M\in \mathscr{G}^{\bot >0}$ and $W\in\mathscr{G}$, we have $\mathscr{M}\subseteq \mathscr{G}^{\bot >0}\subseteq W^{\bot 1}$. For $X\in\Add(M)$, since ${_A}M$ is self-orthogonal and finitely presented, $X^{\bot 1}$ contains $\mathscr{M}$ and is closed under arbitrary direct sums in $A\Modcat$. Thus $(3)$ implies that $\Omega_M^-(N)^{\bot 1}$ is closed under countable direct sums in $A\Modcat$ of modules in $\mathscr{M}$. By Lemma \ref{FMF},  $\Omega_M^-(N)^{\bot 1}$ is closed under countable direct sums in $A\Modcat$ of finitely $M$-filtered $A$-modules. Further, it follows from $N\in \mathscr{E}\cap\mathscr{S}$, Corollary \ref{localizing} and Lemma \ref{E-property}(2) that $\Omega_M^-(N)\in\mathscr{E}\cap\mathscr{S}$. By Proposition \ref{CP}(2), $\Omega_M^-(N)\in\Add(M)$. This shows $\Omega_M^-(N)=0$ (and thus also $N=0$) in $\mathscr{C}$. By $\StEnd_A(N)=\StEnd_M(N)$, we have $N=0$ in $\Stmc{A}$. Since $\StEnd_A(M)\simeq \StEnd_A(N)$ by Lemma \ref{Isomorphism}(2), $M=0$ in $\Stmc{A}$. In other words, ${_A}M$ is projective.

$(3)\Rightarrow (1)$. By Proposition \ref{Compact objects in recollement}(1),
the module $N$ is $M$-compact (that is, compact in $\mathscr{C}$)
and isomorphic in $\Stmc{A}$ to an $M$-filtered module $X$. Then there are projective $A$-modules $P$ and $Q$ with $N\oplus P\simeq X\oplus Q$. Since projective $A$-modules are zero in $\mathscr{C}$, $X$ is $M$-compact. By $(3)$, $X$ lies in $\Add(M)$, and therefore $N=0$ in $\mathscr{C}$. By Lemma \ref{Isomorphism}(2), we have $M=0$ in $\Stmc{A}$. Thus ${_A}M$ is projective. $\square$

\medskip
Let $\mathscr{P}^{<\infty}(B)$ be the category of finitely generated $B$-modules with finite projective dimension, and let $\GP{B}_\omega$ be the category of \emph{countably generated, compactly Gorenstein-projective} $B$-modules. Clearly, finitely generated, Gorenstein-projective $B$-modules are in $\GP{B}_\omega$. Note that the category of countably generated $B$-modules is
a Serre subcategory of $B\Modcat$. This is due to the fact:  A ring $R$ has the property that each submodule of each countably generated left $R$-module is countably generated if and only if each left ideal of $R$ is countably generated.

We consider the following two homological conditions:

\smallskip
${\rm (CI)}$ The direct sum of countably many $B$-modules from $\mathscr{P}^{<\infty}(B)$ belongs to $\GP{B}_\omega^{\,\,\bot >0}$.

${\rm (CII)}$ Any compactly Gorenstein-projective, compactly filtered $B$-module is projective.

\smallskip
The \emph{finitistic dimension} of an Artin algebra $B$ is the supremum of projective dimensions of all $B$-modules in $\mathscr{P}^{<\infty}(B)$. The well-known finitistic dimension conjecture says that an Artin algebra $B$ should always have finite finitistic dimension.
The validity of this conjecture for $B$ implies the one of the Nakayama conjecture for $B$. However, the finitistic dimension conjecture is still open.

\begin{Lem}\label{Three conditions}
$(1)$ If ${\rm (CI)}$ holds, then so does ${\rm (CII)}$.

$(2)$ If an Artin algebra $B$ has finite finitistic dimension or is a virtually Gorenstein algebra, then ${\rm (CI)}$ holds.

\end{Lem}

{\it Proof.}
$(1)$
Let $Y$ be a compactly Gorenstein-projective $B$-module which is compactly filtered by a sequence  $\{Y_i\mid i\in\mathbb{N}\}$ of submodules of $Y$. Set $X:=\bigoplus_{i\in\mathbb{N}}Y_i$. Then there is a canonical exact sequence $0\to X\to X\to Y\to 0$ in $B\Modcat$. Since $Y_i$ are finitely generated for all $i$, the module $Y$ is countably generated. This implies $Y\in \GP{B}_\omega$. Moreover, since $Y_i\in\mathscr{P}^{<\infty}(B)$ for all $i\in\mathbb{N}$, the condition (CI) implies that
the sequence splits and $Y$ is isomorphic to a direct summand of $X$. Note that $Y_i$ is a direct sum of finitely many indecomposable submodules with the local endomorphism rings. By \cite[Corollary 26.6]{AF}, $Y\simeq \bigoplus_{i\in \mathbb{N}}Y_i'$, where $Y_i'$ is a direct summand of $Y_i$ for each $i$.
Since $\mathscr{P}^{<\infty}(B)\cap\GP{B}=\add(B)$, $Y_i'$ is projective for all $i$. Thus $Y$ is projective.  This shows that (CII) holds.

$(2)$ Clearly, $\GP{B}^{\bot >0}$ contains all $B$-modules of finite projective dimension. For a virtually Gorenstein algebra $B$, $\GP{B}^{\bot >0}={^{\bot >0}}\GI{B}$, and therefore $\GP{B}^{\bot >0}$ is closed under arbitrary direct sums in $B\Modcat$. Thus (2) holds.
$\square$

\begin{Lem}\label{EH}
Let $B$ and $C$ be Artin algebras and $X$ a finitely generated $C$-$B$-bimodule. Suppose that ${_C}X$ and $X_B$ are projective and the tensor functor $X\otimes_B-: B\Modcat\to C\Modcat$ induces a triangle equivalence $B\mbox{-}\underline{\rm GProj}\to C\mbox{-}\underline{\rm GProj}$.
If $C$ satisfies ${\rm (CII)}$, then so does $B$.
\end{Lem}

{\it Proof.} Let $F:={_C}X\otimes_B-$. Since both ${_C}X$ and $X_B$ are projective, the functor $F:B\Modcat\to C\Modcat$ is exact and preserves projective modules.
This yields $F(\mathscr{P}^{<\infty}(B))\subseteq \mathscr{P}^{<\infty}(C)$. Since $F$ commutes with filtered colimits, it sends compactly filtered $B$-modules
to compactly filtered $C$-modules. Further, since $F$ induces a triangle equivalence $B\mbox{-}\underline{\rm GProj}\to C\mbox{-}\underline{\rm GProj}$,
it reflects projective modules and sends compactly Gorenstein-projective $B$-modules to compactly Gorenstein-projective $C$-modules.
This implies Lemma \ref{EH}. $\square$

\medskip
Next, we point out that the condition ${\rm (CII)}$ is preserved by several classes of equivalences between algebras.
For the unexplained notions below of stable equivalences of adjoint type and singular equivalences of Morita type with level,
we refer to \cite{xc} and \cite{Wang}, respectively. Given a finitely generated $B$-module $N$,
we denote by $\thick({_B}N)$ the smallest thick subcategory of $B\modcat$ which contains $N$.

\begin{Koro}\label{DM}
Let $B$ and $C$ be finite-dimensional algebras over a field. Suppose that

\smallskip
$(a)$ $B$ and $C$ are derived equivalent, or

$(b)$ $B$ and $C$ are stably equivalent of adjoint type, or

$(c)$ $B$ and $C$ are singularly equivalent of Morita type with level defined by a pair of bimodules $({_C}X{_B}, {_B}Y_C)$ such that
$\Hom_C(X, C)\in\thick({_B}B\oplus D(B))$ and $\Hom_B(Y, B)\in\thick({_C}C\oplus D(C))$.

Then $B$ satisfies ${\rm (CII)}$ if and only if so does $C$.
\end{Koro}

{\it Proof.} $(a)$ follows from Lemma \ref{EH} and \cite[Example 4.7 and Corollary 5.4]{HP}. The case $(b)$ is a special case of $(c)$.
In the case $(c)$, the modules ${_C}X$ and $X_B$ are finitely generated and projective by the definition of singular equivalences of Morita type with level.
By Lemma \ref{EH}, it suffices to show that the functor $F:={_C}X\otimes_B-: B\Modcat\to C\Modcat$ induces a triangle equivalence $B\mbox{-}\underline{\rm GProj}\to C\mbox{-}\underline{\rm GProj}$. Note that for finitely generated Gorenstein-projective modules, the triangle equivalence $B\mbox{-}\underline{\rm Gproj}\to C\mbox{-}\underline{\rm Gproj}$ is known in \cite[Proposition 4.5]{Wang} under a stronger assumption that $\Hom_C(X, C)\in\mathscr{P}^{<\infty}(B)$ and $\Hom_B(Y, B)\in\mathscr{P}^{<\infty}(C)$. Since our discussions involve infinitely generated Gorenstin-projective modules, we have to check whether $F$ and $G$ can be regraded as functors between  $\GP{B}$ and $\GP{C}$.

Let $U\in\GP{B}$ with a complete projective resolution $\cpx{P}$. Then $F(\cpx{P}):=(F(P^i)_{i\in\mathbb{Z}})$ is an exact complex of projective $C$-modules. Moreover, $\cpx{\Hom}_C(F(\cpx{P}), C)\simeq\cpx{\Hom}_B(\cpx{P}, \Hom_C(X,C))$ as complexes.
By \cite[Theorem 2]{BK}, $\thick({_B}B\oplus D(B))=\GP{B}^{\bot >0}\cap B\modcat$. Since each cocycle of $\cpx{P}$ lies in $\GP{B}$ and
$\Hom_C(X, C)\in\thick({_B}B\oplus D(B))$, the complex $\cpx{\Hom}_B(\cpx{P}, \Hom_C(X,C))$ is exact. Consequently, $F(\cpx{P})$ is a compete projective resolution of $F(U)$, and therefore $F(U)\in \GP{C}$. In other words, $F$ restricts to a functor $\GP{B}\to\GP{C}$. Similarly, the functor $G:={_B}Y\otimes_C-: C\Modcat\to B\Modcat$ also restricts to a functor $\GP{C}\to\GP{B}$.

Recall that $B\mbox{-}\underline{\rm GProj}$ is a triangulated category with the shift functor $\Sigma_B$ which is a quasi-inverse of the syzygy functor $\Omega_B$. This implies $Y\otimes_CX \otimes_B-\simeq \Omega^n_B(-)\simeq \Sigma_B^{-n}: B\mbox{-}\underline{\rm GProj}\lra B\mbox{-}\underline{\rm GProj},$ where $n$ is the level of the singular equivalence in $(c)$ and $\Sigma_B^{n}$ denotes the $n$-th shift functor of $B\mbox{-}\underline{\rm GProj}$. Similarly, there are equivalences $X\otimes_BY \otimes_C-\simeq \Omega^n_C(-)\simeq \Sigma_C^{-n}: C\mbox{-}\underline{\rm GProj}\ra C\mbox{-}\underline{\rm GProj}.$ Thus
$F$ induces a triangle equivalences $B\mbox{-}\underline{\rm GProj}\lra C\mbox{-}\underline{\rm GProj}$ with a quasi-inverse
$\Sigma_B^{n}\circ G$: $C\mbox{-}\underline{\rm GProj}\to B\mbox{-}\underline{\rm GProj}$. $\square$

\medskip
Recall from Definition \ref{Gorenstein-Morita}(ii) that an algebra $B$ is said to be \emph{compactly Gorenstein} if ${\rm (CII)}$ holds; and \emph{Gorenstein-Morita} if $B$ is both strongly Morita and compactly Gorenstein. By Lemma \ref{Three conditions}, a strongly Morita algebra is Gorenstein-Morita if it is virtually Gorenstein or has finite finitistic dimension.

\medskip
{\bf Proof of Corollary \ref{Gendo-symmetric}.}
Suppose that $\Lambda$ is a Gorenstin-Morita algebra. Then there is a self-injective algebra $A$ and a finitely generated $A$-module $M$ which is a generator such that $\Lambda =\End_A(M)$
and $\add(M)=\add(\nu_A(M))$. Further, suppose  $\dd(\Lambda)=\infty$. Then ${_A}M$ is self-orthogonal by \cite[Lemma 3]{Muller}.

In the following, we show that $N:=H_*^0(M)=0$ in $\mathscr{C}$ (see Section \ref{Torsion pair}). If it is the case, then $\StEnd_A(M)\simeq \StEnd_A(N)$ by Lemma \ref{Isomorphism}(2), and therefore $M=0$ in $\Stmc{A}$, that is, ${_A}M$ is projective. Thus $\Lambda$ is Morita equivalent to $A$ and must be self-injective, and therefore the first part of Corollary \ref{Gendo-symmetric} is proved.

Now, let $M=A\oplus \bigoplus_{i=1}^{m}M_i$ with $m\in\mathbb{N}$ and $M_i$ indecomposable and non-projective for $1\le i\le m$.
By Lemma \ref{Gamma-mod}(2), $N\simeq \bigoplus_{i=1}^{m}H_*^0(M_i)$
in $\mathscr{C}$ and $H_*^0(M_i)$ are indecomposable in $\mathscr{C}$ for $1\le i\le m$.  Clearly, $H_*^0(M_i)$ as an $A$-module does not contain non-projective direct summands in $\Add(M)$, due to $\StHom_A(N, M)=0$. So we can assume that $H_*^0(M_i)=P_i\oplus N_i$ as $A$-modules, where $P_i$ is  a projective $A$-module, and $N_i$ is either zero or an indecomposable $A$-module that does not lie in $\Add(M)$. Now, we turn to show $N_i=0$ for all $i$. This implies that $N$ is projective, and thus
$N=0$ in $\mathscr{C}$.

By Theorem \ref{Compact}(1) and Proposition \ref{CP}(1), $N_i$ is isomorphic in $\Stmc{A}$ to an $M$-filtered $A$-module $X$. Now, suppose that $X$ is filtered by a sequence $\{X_n\mid n\in\mathbb{N}\}$ of submodules of $X$ (see Definition \ref{Countably filtered}). Then $X=\lim\limits_{n\rightarrow\infty}X_n=\bigcup_{n=0}^{\infty}X_n$. By Lemma \ref{FMF}, $X_n\in M^{\bot >0}$ and has an $\add(M)$-resolution of finite length.
Applying $G:=\Hom_A(M,-)$ to the inclusions $X_n\to X_{n+1}$, we obtain $G(X_n)\in\mathscr{P}^{<\infty}(\Lambda)$ and injective homomorphisms $G(X_n)\to G(X_{n+1})$. Since ${_A}M$ is finitely presented, $G$ commutes with filtered colimits. It follows that $G(X)=\lim\limits_{n\rightarrow\infty} G(X_n)=\bigcup_{n=0}^{\infty}G(X_n)$. In other words, $G(X)$ is compactly filtered.

Assume $N_i\neq 0$. Then $N_i$ is indecomposable and does not belong to $\Add(M)$. It follows from $X\simeq N_i$ in $\Stmc{A}$ and $N_i\in\mathscr{G}$ that $X\in\mathscr{G}$ and $\Omega_M^-\Omega_M(X)\simeq\Omega_M^-\Omega_M(N_i)$. By Remark \ref{Decomposition}, $\Omega_M^-\Omega_M(N_i)\simeq N_i$, and therefore $X\simeq N_i\oplus M_0$ for some $M_0\in\Add(M)$. As $N$ is compact in $\mathscr{C}$ by Theorem \ref{Compact}(1), $X$ is also compact in $\mathscr{C}$. It follows from Lemma \ref{Gorenstein}(2) that $G(X)$ is compactly Gorenstein-projective, but not projective. Since $\Lambda$ is compactly Gorenstein, $G(X)$ must be projective. This leads to a contradiction. Thus $N_i=0$.

If $\Lambda$ is gendo-symmetric, then $A$ can be assumed to be symmetric. Since gendo-symmetric, virtually Gorenstein algebras are Gorenstein-Morita, the second part of Corollary \ref{Gendo-symmetric}
follows from the first part of Corollary \ref{Gendo-symmetric}. $\square$

\medskip
To end this section, we mention the following questions related to the results in this paper.

\smallskip
{\bf Question} 1. Find necessary and sufficient conditions for Artin algebras to be compactly Gorenstein.

Finitely generated Gorenstein-projective modules over an Artin algebra are compactly Gorenstein-projective. The next question is about countably generated Gorenstein-projective modules.

{\bf Question} 2. Describe the class of countably generated, compactly Gorenstein-projective modules over an Artin algebra.

{\bf Question} 3. Find more new examples of Gorenstein-Morita algebras.

\medskip
{\bf Acknowledgements.} Both authors would like to thank Ming Fang for discussions on gendo-symmetric algebras and stable categories. The research work was supported partially by the National Natural Science Foundation of China (Grant 12031014, 12122112 ).

{\footnotesize
\smallskip
Hongxing Chen,

School of Mathematical Sciences  \&  Academy for Multidisciplinary Studies, Capital Normal University, 100048
Beijing,

P. R. China

{\tt Email: chenhx@cnu.edu.cn (H.X.Chen)}

\smallskip
Changchang Xi,

School of Mathematical Sciences, Capital Normal University, 100048
Beijing, P. R. China; and
School of Mathematics and Statistics, Shaanxi Normal University, 710119 Xi'an, P. R. China

{\tt Email: xicc@cnu.edu.cn (C.C.Xi)}
}


\medskip
{\footnotesize
\begin{thebibliography}{99}
\bibitem{AF} {{\sc F. W. Anderson} and {\sc K. R. Fuller}, \emph{Rings and categories of modules}, Second Edition.
\textit{Graduate Texts in Mathematics} \textbf{13}, Springer-Verlag 1992.}

\bibitem{Asashiba}{{\sc H. Asashiba,} The selfinjectivity of a local algebra $A$ and the condition $\Ext^1_A(DA, A)=0,$ The 5th Int. Conf. Represent. Algebras (Tsukuba,1990), \emph{Canad. Math. Soc. Conf. Proc.} \textbf{11} (1991) 9-23.}

\bibitem{APT}{{\sc M. Auslander}, {\sc I.M. Platzeck} and {\sc G. Todorov}, Homological theory of idempotent ideals, {\it Trans. Amer. Math. Soc.}
{\bf 332} (1992) 667-692.}

\bibitem{abs}{{\sc L. L. Avramov, R.-O. Buchweitz} and {\sc L. M. Sega}, Extensions of a dualizing complex by its ring: Commutative versions of a conjecture of Tachikawa, \emph{J. Pure Appl. Algebra} \textbf{201} (2005) 218-239.}

\bibitem{AR}{{\sc M. Auslander} and {\sc I. Reiten}, On a generalized version of the Nakayama conjecture, {\it Proc. Amer. Math. Soc.} 52 (1975) 69-74.}

\bibitem{AR1} {{\sc M. Auslander} and {\sc I. Reiten}, Represenbtation theory of Artin algebras III. Almost split sequences,
{\it Comm. Algebra} {\bf 3} (1975) 239-294.}

\bibitem{ar1990}{{\sc M. Auslander} and {\sc I. Reiten}, Applications of contravariantly finite subcategories, {\it Adv. Math.} 86 (1991) 111-152.}
\bibitem{ars}{{\sc M. Auslander, I. Reiten} and {\sc S. Smal{\o}}, \emph{Representation theory of Artin algebras}, Cambridge Studies in Advanced Mathematics
\textbf{36}. Cambridge University Press, Cambridge, 1995.}

\bibitem{AINS}{{\sc L. L. Avramov, S. B. Iyengar, S. Nasseh} and {\sc K. Sather-Wagstaff}, Persistence of homology over commutative noetherian rings, \emph{J. Algebra} \textbf{610} (2022) 463-490.}

\bibitem{BBD}{{\sc A. A. Beilinson, J. Bernstein} and {\sc P.
Deligne},  Faisceaux pervers, \emph{Asterisque} \textbf{100} (1982) 5-171.}

\bibitem{Bel} {{\sc A. Beligiannis}, Cohen-Macaulay modules, (co)torsion pairs and virtually Gorenstein algerbas, {\it J. Algebra}. {\bf 288} (2005) 137-211.}

\bibitem{Bel2} {{\sc A. Beligiannis}, On algebras of finite Cohen-Macaulay type, {\it Adv. Math.} {\bf 226} (2011) 1973-2019.}

\bibitem{BI}{{\sc A. Beligiannis} and {\sc I. Reiten}, Homological and homotopical aspects of torsion theories, {\it Mem. Amer. Math. Soc}.
{\bf 188} (2007), no. 883, 1-207.}

\bibitem{BK} {{\sc A. Beligiannis} and {\sc H. Krause}, Thick subcategories and virtually Gorenstein algerbas,
{\it Illinois J. Math}. {\bf 52} (2008) 551-562.}

\bibitem{c}{{\sc H. X. Chen}, Applications of hyperhomology to adjoint functors, {\it Comm. Algebra} \textbf{50} (1) (2022) 19-32.}

\bibitem{cs}{{\sc H. X. Chen} and {\sc S. Koenig}, Ortho-symmetric modules, Gorenstein algebras, and derived equivalences,
\emph{Int. Math. Res. Not.} IMRN\textbf{2016} (2016) 6979-7037.}

\bibitem{xc1}{{\sc H. X. Chen} and {\sc C. C. Xi}, Good tilting
modules and recollements of derived module categories, \emph{Proc.
Lond. Math. Soc.} \textbf{104} (2012) 959-996.}

\bibitem{xc6}{{\sc H. X. Chen} and {\sc C. C. Xi}, Dominant dimensions, derived equivalences and tilting modules, \emph{Isr. J. Math.}
\textbf{215} (2016) 349-395.}

\bibitem{xcf}{{\sc H. X. Chen}, {\sc M. Fang} and {\sc C. C. Xi}, Tachikawa's second conjecture, derived
recollements, and gendo-symmetric algebras, \emph{Compos. Math.} (to appear), arXiv: 2211.08037.}

\bibitem{DH}{{\sc P. Dr\"axler} and {D. Happel,}
A proof of the generalized Nakayama conjecture for algebras with $J^{2l+1}=0$ and $ A/J^l$ representation finite, \emph{J. Pure Appl. Algebra} \textbf{78} (2) (1992) 161-164.}

\bibitem{FK} {{\sc M. Fang} and {\sc S. Koenig},
 Endomorphism algebras of generators over symmetric algebras,
 \textit{J. Algebra} \textbf{332} (2011) 428-433.}

\bibitem{GoTa}{{\sc S. Goto} and {\sc R. Takahashi}, On the Auslander-Reiten conjecture for Cohen-Macaulay local rings, \emph{Proc. Amer. Math. Soc.} \textbf{145} (8) (2017) 3289-3296.}

\bibitem{Hoshino}{{\sc M. Hoshino}, Modules without self-extensions and Nakayama's conjecture, \emph{ Arch. Math. (Basel)}
\textbf{43} (1984) 493-500.}

\bibitem{HKM}{{\sc M. Hoshino}, {\sc Y. Kato} and {\sc J-I. Miyachi}, On $t$-structures and torsion theories induced by
compact objects, \emph{J. Pure Appl. Algebra}
\textbf{167} (2002) 15-35.}

\bibitem{HP}{{\sc W. Hu} and {\sc S. Y. Pan}, Stable functors of derived equivalences and Gorenstein projective modules,
\emph{Math. Nachr.} \textbf{290} (2017) 1512-1530.}

\bibitem{Huang}{{\sc Z. Y. Huang}, Proper resolutions and Gorenstein categories, \emph{J. Algebra} \textbf{393} (2013) 142-169.}

\bibitem{keller}{{\sc B. Keller}, Deriving DG categories, \emph{Ann. Sci.
\'Ecole Norm. Sup.} {\bf 27} (1994) 63-102.}

\bibitem{KY} {{\sc O. Kerner} and {\sc K. Yamagata}, Morita algebras,
\emph{J. Algebra} \textbf{382} (2013) 185-202.}

\bibitem{Kr0}{{\sc H. Krause}, A Brown representability theorem via coherent functors, \emph{Topology} \textbf{41} (2002) 853-861.}

\bibitem{Kr1}{{\sc H. Krause}, The stable derived category of a noetherian scheme, \emph{Compos. Math.}
\textbf{141} (2005) 1128-1162.}

\bibitem{Muller} {{\sc B. M\"{u}ller}, The classification of algebras by dominant dimension, \emph{Canad. J. Math.} \textbf{20} (1968) 398-409.}

\bibitem{Nakayama}{{\sc T. Nakayama}, On algebras with complete homology, \emph{Abh. Math. Sem. Univ.
Hamburg} \textbf{22} (1958) 300-307.}

\bibitem{Neemanbook}{{\sc A. Neeman}, \emph{Triangulated categories}, Annals of Mathematics Studies {\bf 148},
Princeton University Press, Princeton, 2001.}

\bibitem{Neeman1}{{\sc A. Neeman}, The homotopy category of injectives, \emph{Algebra Number Theory} \textbf{8} (2014) 429-456.}

\bibitem{Rk}{{\sc J. Rickard}, Idempotent modules in the stable category,  \emph{J. Lond. Math. Soc}. \textbf{56} (1997) 149-170.}

\bibitem{Rk2}{{\sc J. Rickard}, Unbounded derived categories and the finitistic dimension conjecture,  \emph{Adv. Math}. \textbf{354} (2019) 106735, 21pp.}

\bibitem{RR}{{\sc R. Rouquiver}, Dimensions of triangulated categories,  \emph{J. K-Theory} \textbf{1} (2008) 193-256.}

\bibitem{SSW}{{\sc S. Sather-Wagstaff}, {\sc T. Sharif} and {\sc D. White}, Stability of Gorenstein categories, \emph{J. Lond. Math. Soc.}
\textbf{77} (2008) 481-502.}

\bibitem{Sch}{{\sc R. Schulz}, Boundedness and periodicity of modules over QF rings, {\it J. Algebra} \textbf{101} (1986) 450-469.}

\bibitem{Tac} {{\sc H. Tachikawa}, \emph{Quasi-Frobenius rings and generalizations},
Lecture Notes in Mathematics \textbf{351}, Berlin-Heidelberg New York 1973.}

\bibitem{Wang}{{\sc Z. F. Wang}, Singular equivalence of Morita type with level, \emph{J. Algebra} \textbf{439} (2015) 245-269.}

\bibitem{Weibel1}{{\sc Ch. Weibel}, \emph{An introduction to homological algebra}, Cambridge University Press. Cambridge, 1994.}

\bibitem{xc}{{\sc C. C. Xi}, Stable equivalences of adjoint type, \emph{Forum Math.}
\textbf{20} (2008) 81-97.}

\bibitem{Yamagata}{{\sc K. Yamagata}, Frobenius algebras, in: \emph{Handbook of Algebra}, vol. 1, North-Holland, Amsterdam, 1996, pp. 841-887. }
\end{thebibliography}}
\end{document}